\documentclass[twoside,11pt]{article}

\usepackage{multirow}
\usepackage{booktabs}
\usepackage{jmlr2e}
\usepackage{float}
\usepackage{natbib}
\usepackage{graphicx}
\usepackage{subfigure}
 \usepackage{algorithm}
 \usepackage{algorithmic}
 \usepackage{cite}
 \usepackage{boxedminipage}
 \usepackage{float}
 \usepackage{amsmath}
\RequirePackage{hyperref}
\RequirePackage{nameref}
\hypersetup{colorlinks, linkcolor=blue, citecolor=blue, urlcolor=magenta, linktocpage, plainpages=false}
\usepackage{bm}
\usepackage{footnote}
\usepackage{tablefootnote}
\allowdisplaybreaks[3]
 \newtheorem{lem}{{Lemma}}
 \newtheorem{assumption}{Assumption}
 \newtheorem{cor}{{Corollary}}
 
 \newcommand\argmin{\mathop{\textrm{argmin}}}
 \newcommand\crule[1][5cm]{%
 	\par
 	\nointerlineskip
 	\centerline{\hbox to #1{\hrulefill}}%
 	\nointerlineskip}

\newtheorem{thm}{{Theorem}}
 \newtheorem{rem}{Remark}

 \usepackage[misc]{ifsym}

 \newcommand{\N}{\mathbb{N}}
 \newcommand{\R}{\mathbb{R}}
 \newcommand{\E}{\mathbb{E}}

\def\bigbl{\left\lbrace}
\def\bigbr{\right\rbrace}

\begin{document}


\firstpageno{1}

\title{On the Convergence of Stochastic Gradient Descent with Bandwidth-based Step Size}

\author{\name Xiaoyu Wang \email wxy@lsec.cc.ac.cn \\
       \addr Academy of Mathematics and Systems Science \\Chinese Academy of Sciences\\ Beijing 100190, China\\ 
       University of Chinese Academy of Sciences \\
       No.19A Yuquan Road, Beijing 100049, China
       \AND Ya-xiang Yuan 
       \email yyx@lsec.cc.ac.cn\\
       \addr State Key Laboratory of Scientific/Engineering Computing, Institute of Computational
       Mathematics and Scientific/Engineering Computing, Academy of Mathematics and Systems Science\\Chinese Academy of Sciences\\Beijing 100190, China}

\editor{}

\maketitle

\begin{abstract}
We investigate the stochastic gradient descent (SGD) method where the step size lies within a banded region instead of being given by a fixed formula. The optimal convergence rate under mild conditions and large initial step size is proved. Our analysis provides comparable theoretical error bounds for SGD associated with a variety of step sizes. In addition, the convergence rates for some existing step size strategies, e.g., triangular policy and cosine-wave, can be revealed by our analytical framework under the boundary constraints. The bandwidth-based step size provides efficient and flexible step size selection in optimization. We also propose a $1/t$ up-down policy and give several non-monotonic step sizes. Numerical experiments demonstrate the efficiency and significant potential of the bandwidth-based step-size in many applications.

\end{abstract}
\vspace{0.03in}

\section{Introduction}\label{sec:1}

In this paper, we consider the stochastic optimization problem as follows 
\begin{equation}\label{P1}
\min_{x\in \R^d} f(x) = \E_{\xi \sim \Xi }[ f(x; \xi)],
\end{equation} 
where $\xi$ is a random variable drawn from some source distribution $\Xi$. This problem is often encountered in machine learning and statistics, and attracts much attention along with the big data and artificial intelligence.  The corresponding empirical risk problem is to minimize $f(x) = \frac{1}{n}\sum_{i=1}^{n} f(x; \xi_i)$,
where each $\xi_i$ ($ i \in \left\lbrace 1, 2,\dots, n \right\rbrace $) denotes a realization of $\xi$.

The stochastic gradient descent (SGD) algorithm~\citep{SGD-1951} is widely used to solve the above problems in which the iterative scheme is
\begin{equation}\label{SGD}
x_{t+1} = x_{t} - \eta(t) \nabla f(x_t;\xi_{i_t}),
\end{equation}
where $\eta(t) > 0$ is the step size and $\xi_{i_t}$ is a realization of $\xi$ at iteration $t$ (or $i_t$ is chosen uniformly at random from  $\left\lbrace 1, 2, \cdots, n\right\rbrace $). The stochastic gradient $g_t := \nabla f(x_t;\xi_{i_t})$ satisfies
\begin{equation*}
\E[ \nabla f(x_t;\xi_{i_t}) \mid \mathcal{F}_{t}\footnote{We use $\mathcal{F}_t$ to denote $\sigma$-algebra of the random information at  iteration $t$.}] = \nabla f(x_t).
\end{equation*}

SGD is often
preferable in large-scale machine learning because of its simplicity and low-cost computation per iteration. However, determining the step size is of the key importance and challenging due to the gradient noise of SGD. In this paper we focus on the mini-batch version of SGD where its gradient is estimated on a small sample set $\Omega_t$ ($|\Omega_t|=b$), i.e., $g_t = \frac{1}{b}\sum_{i\in \Omega_t}\nabla f(x_t, \xi_i)$.
\subsection{Theoretical Analysis of SGD for Various Step Sizes}\label{literature_review}

The analysis of SGD at beginning  lies on the asymptotic results \citep{chung1954,leen1994optimal,leen1998two}. \citet{leen1994optimal} analyzed the asymptotic properties around the locally optimal solution $x^{\ast}$ with $\eta(t) = \eta_0/t$ and show that if $ \eta_0 > 1/(2\lambda_{\min})$ ($\lambda_{\min}$ is the smallest eigenvalue of $\nabla^2 f(x^{\ast})$), the error $\E[\left\|x_t-x^{\ast} \right\|^2 ] $ has order $ \mathcal{O}(1/t)$, which is an optimal rate (minimax rate) \citep{polyak1992, agarwal2009information, ghadimi2012optimal}. 

Recently, the focus has been shifted to study the non-asymptotic results.
\citet{Moulines-Bach2011} established the convergence rate of SGD for a class of step size $\eta(t) = \eta_0/t^p$ for $p\in (0,1]$. For strongly convex and $L$-smooth functions, SGD exhibits an optimal error bound $\mathcal{O}(1/T)$ ($T$ is the number of iterations) with $\eta(t) = \eta_0/t$~\citep{Moulines-Bach2011,Rakhlin-Shamir-Sridharan2011,nguyen2019tight}. 
However, the situation becomes complicated if the function is not $L$-smooth. The best known result on the last iterate is $\E[f(x_{T}) - f(x^{\ast})] \leq O(\log(T)/T)$ with $\eta(t) = 1/(\mu t)$~\citep{Shamir-Zhang2013}, which is proved to be tight~\citep{harvey2018tight}. Many averaging techniques such as suffix averaging \citep{Rakhlin-Shamir-Sridharan2011} and polynomial-decay averaging~\citep{Shamir-Zhang2013,Lacoste-Schmidt-Bach2012} are incorporated into SGD and obtain an optimal $\mathcal{O}(1/T)$ rate. \citet{hazan2014beyond} achieved an $\mathcal{O}(1/T)$ convergence rate by exponentially decreasing the step size after a consecutive period which grows exponentially,  and adopting a simple modification where the inner iterations are averaged as an output. \citet{jain2019making} designed the piece-wise decay step size with the form of $\mathcal{O}(1/t)$ per period and obtained  an optimal error bound $\E[f(x_{T}) - f(x^{\ast})] \leq O(1/T)$ on the last iterate. \citet{ge2019step} analyzed a step decay step size which decays exponentially after $T/\log(T)$ epochs  and achieved a near-optimal $\mathcal{O}(\log(T)/T)$ convergence rate for the least squares problems. Very recently, \citet{li2020exponential} proposed the continuous version of the step decay step size and proved a near-optimal convergence rate when the Polyak-L\'ojasiewicz condition holds.

To the best of our knowledge, there are many other efficient step sizes preferred in deep learning, e.g., adaptive methods \citep{AdaGrad,RMSProp,zeiler2012adadelta,Adam,S_polyak},  Barzilai-Borwein based~\citep{tan2016barzilai, yang2018random}, line-search based~\citep{keskar2015nonmonotone, vaswani2019painless},
cyclical learning rate (step size) \citep{smith2017cyclical, loshchilov2016sgdr, sinewave}.\\

\subsection{Motivation}
In this paper, we focus on the non-asymptotic convergence rate of the SGD method in which the step size $\left\lbrace \eta(t)\right\rbrace $ varies in a bounded region rather than a fixed schedule. The lower and upper bounds of the region are defined by two functions $\delta_1(t)$ and $\delta_2(t)$ w.r.t. the iteration number $t$. More specifically, we assume there exist two positive constants $ m \leq M$ such that 
\begin{equation}\label{eta_condition}
m \delta_1(t) \leq \eta(t) \leq  M \delta_2(t), \,\forall \, t \geq 1,\tag{BD}
\end{equation}
and $d\delta_1(t)/dt \leq 0$ and $d\delta_2(t)/dt \leq 0$.  When $\delta_1(t) = \delta_2(t) = 1/t$, we call it $1/t$-band. Such an idea is motivated by the piece-wise decay step size~\citep{hazan2014beyond, jain2019making, ge2019step}, which is a step function whose graph consists of some line segments lying within two curves (i.e., their lower and upper bounds).  The diminishing step size $\eta(t) = \eta_0/t$ and
 piece-wise decay step sizes proposed by \citet{hazan2014beyond} and \citet{jain2019making} can be regarded as the special cases of $1/t$-band.
 
 \citet{dauphin2014identifying} pointed out that a great obstacle to minimize deep neural networks with high possibility arose from saddle points instead of poor local minima.  The proposed non-monotonic scheduling (\ref{eta_condition}), admitting some intermediate increasing in step size, might help
rapidly traverse the saddle points and find flat minima. \citet{smith2017cyclical} described a type of cyclical learning rate (step size) which varies within a band of minimum and maximum values and showed the potential benefits to train deep neural networks. Similarly, \citet{sinewave} proposed a sine-wave learning rate framework.  Their boundaries decay exponentially after a few fixed epochs. The policy lets the step size locally vary within a reasonable band. Although their mechanisms might have a short-term negative effect, it is beneficial overall.

 We are interested in the bandwidth-based step size described in (\ref{eta_condition}), which has more freedom to be used to design more efficient  step sizes in practice. Although there are many specific and effective schedules mentioned in Section \ref{literature_review}, it is still an interesting and challenging topic to analyze the convergence properties of the SGD method based on such a generic step size. Moreover, some step sizes, e.g., exponentially decaying step size after a few fixed epochs and cyclical learning rate \citep{smith2017cyclical, loshchilov2016sgdr, sinewave}, lack non-asymptotic convergence guarantees. To overcome the limitations, we explore their connections in theory and practice using the bandwidth-based step size (\ref{eta_condition}).

\subsection{Main Contributions}
The main contributions of this work include: (1) the proof of the optimal convergence rate of SGD w.r.t the bandwidth-based step size (\ref{eta_condition}), (2) the error bound analysis for bandwidth-based step size with the same boundary order, (3) the error bound in terms of the different boundary orders, and (4) development of the $1/t$ up-down policy and four non-monotonic step sizes with demonstrated efficiency. 
More precisely, we make the following contributions.

First, we investigate the step sizes lying in a bounded region to achieve an $\mathcal{O}(1/T)$ convergence rate for strongly convex problems. 
The main results are summarized in Table \ref{tab:11}. The commonly used step size $\eta(t) = \eta_0/t$ is extended to $1/t$-band 
which allows the step size to vary locally like parallel line, triangular, cosine-wave or other ways. We further relax its lower and upper bounds, which can provide theoretical guarantees for large step sizes on the initial iterations. 
 Our analysis also provides the error bounds $\E[\left\|x_{T+1}-x^{\ast}\right\|^2]\, \text{or}\, \E[f(\hat{x}\footnote{Here $\hat{x}_{T}$ is a type of averaging of the previous iterations $x_t$ from $t=1,2,\cdots, T$.}) - f(x^{\ast})]\leq \mathcal{O}(1/T)$ for some step sizes, e.g., $\eta(t) = \eta_0/t$ and piece-wise decay step size \citep{hazan2014beyond, jain2019making}, which are comparable to those of \citet{Moulines-Bach2011}, \citet{Shamir-Zhang2013} and \citet{Lacoste-Schmidt-Bach2012} but slightly worse than that of \citet{jain2019making}. Actually, the convergence rate proposed in \citet{jain2019making} is hardly achieved under the general framework (\ref{eta_condition}).  In particular, the cyclical step sizes developed in \citet{smith2017cyclical} and \citet{ sinewave}, which lack convergence guarantees, can achieve an $\mathcal{O}(1/T)$ convergence rate if their boundaries satisfy the cases discussed in Table \ref{tab:11}. 
\begin{table}[h]
	\centering
	\renewcommand{\arraystretch}{1.2}	
	\begin{tabular}{|l||c|l||c|c|}
\hline	
 \multicolumn{2}{|c|}{ $\delta_1(t)$}  & \multicolumn{2}{|c|}{$\delta_2(t)$} &  Theorem\\
\hline
{\bf $(A) $}& $1/t$ & $(B)$&  $1/t$ &  1 \& 2 \\ \hline 
$(A_1)$  &  $	\sum_{t=t^{\ast}}^{T} \eta(t) \geq C\ln\left((T+1)/t^{\ast}\right)$ & $(B)$ & $1/t$ & 3 \\ \hline
\multirow{2}{*}{$(A)$}&   \multirow{2}{*}{ $1/t$} & \multirow{2}{*}{$(B_1)$} & $1/t^{r}$, $t \in [C_1T^p]$\tablefootnote{We use $[C_1T^p]$ to denote a positive integer set from 1 to $C_1T^p$, which is also suited for $[T]$.} & \multirow{2}{*}{4 }	\\
 & & & $ 1/t$, $t \in [T]\backslash [C_1T^p]$\tablefootnote{Let $[T]\backslash [C_1T^p] $ denote a positive integer set from $C_1T^p+1$ to $T$.}  &   \\\hline
\multirow{2}{*}{$(A_2)$}& 1, $t \in [C_1T^p]$ & \multirow{2}{*}{$(B_2)$}& 1, $t \in [C_1T^p]$ & \multirow{2}{*}{5} \\
&  $1/t$, $t \in  [T] \backslash [C_1T^p]$ & & $ 1/t$, $ t \in [T]\backslash [C_1T^p]$ & \\\hline
	\end{tabular}
\caption{The bandwidth-based step sizes in Section \ref{sec:3} to achieve $\mathcal{O}(\frac{1}{T})$}	
	\label{tab:11}
\end{table}

Second, we propose the error bound analysis for the cases that $\delta_1(t) = \delta_2(t) = \delta(t)$ and results are shown in Table \ref{tab:12}, where $\delta(t)$ satisfies (\ref{stepsize_cond3}) proposed by~\citet{nguyen2018new}. The results are comparable to those in the prior literature for $\eta(t) = 1/t^p$ ($p\in (0,1]$)~\citep{Moulines-Bach2011}. When $\lim_{t\rightarrow \infty} \delta(t)t =0$, the result is novel.  In particular, we add a condition $-d\delta(t)/dt \leq c_1 \delta(t)^2$ which clarifies ``in the most general case'' mentioned in~\citet{nguyen2018new} and give a more rigorous proof. Moreover, our analysis provides better upper bounds in some cases such as $\eta(t)=1/\sqrt{t}$ and $1/(t\log(t))$ than those of theorem 10 in~\citet{nguyen2018new}.
 \begin{table}[h]
	\centering
	\renewcommand{\arraystretch}{1.1}	
	\begin{tabular}{|c||c|c|c|}	
	\hline
		\multicolumn{2}{|c|}{Conditions} & $\E[\left\|x_{T+1}- x^{\ast}\right\|^2]$& Theorem   \\ \hline
\multirow{3}{*}{$\delta_1(t) = \delta_2(t)$ }	  & $\lim_{t\rightarrow \infty} t \delta(t) = 1$  &  $\mathcal{O}(1/T^{\tau\mu m}) + \mathcal{O}(1/T)$ &  1 \\		
& $\lim_{t\rightarrow \infty} t \delta(t) = 0$ & $\mathcal{O}(\exp(-\tau\mu m \sum_{t=1}^{T}\delta(t)))$&  \multirow{2}{*}{ 6}\\
\,\,\,\,\, \,\,$=\delta(t)$	& $\lim_{t\rightarrow \infty} t \delta(t) = \infty$ &  $\mathcal{O}(\delta(t) )+ \mathcal{O}(\exp(-\tau\mu m \sum_{t=1}^{T}\delta(t)))$  & \\	\hline\hline
		\multirow{6}{*}{$\delta_1(t) \neq \delta_2(t)$}  & $\delta_1(t)= 1/t$ & \multirow{2}{*}{$\mathcal{O}(\log^2(T)/T)$}  
		&     \multirow{2}{*}{7} \\ 
		& $\delta_2(t)= \log(t)/t$ & &   \\ \cline{2-4}		
    		&  $\delta_1(t)= 1/t$ & \multirow{2}{*}{$\mathcal{O}(1/T^{2\alpha-1})$ }&    \multirow{2}{*}{ 8}\\
		& $\delta_2(t)= 1/t^{\alpha}$ &  &   \\\cline{2-4}	
		& $\delta_1(t)= 1/(t\log(t))$ &\multirow{2}{*}{ $\mathcal{O}(1/\log(T)^{\tau\mu m})$ } &   \multirow{2}{*}{9} \\
		& $\delta_2(t)= 1/t^{\alpha}$ &  &    \\ \hline
	\end{tabular}
	\caption{Summary of the bandwidth-based step sizes discussed in Section \ref{sec:4} and \ref{sec:5}}	
	\label{tab:12}
\end{table}
 
Third, we discuss the cases of the  lower and upper bounds being in different orders (i.e., $\delta_1(t) \neq \delta_2(t)$), listed in Table \ref{tab:12}.  The theoretical results explore the connections between the band and its boundaries and broaden the boundaries of the step size for analyzing the convergence behaviors of SGD.

Finally, we propose a $1/t$ up-down policy and design four non-monotonic step sizes including $1/t$ Fix-period band, $1/t$ Grow-period band, $1/t$ Grow-Exp and $1/t$ Fix-Exp. We test regularized logistic regression and some nonconvex problems (e.g., deep neural networks, VGG-16 \citep{Simonyan2015}  and ResNet-18 \citep{he2016deep}) on several real datasets (e.g., MNIST, CIAFR-10 and CIFAR-100). The numerical experiments verify the efficiency compared to their baselines $\eta(t) = \eta_0/t$, exponentially decaying step size \citep{hazan2014beyond} respectively.  We also implement other default algorithms in deep learning, e.g., averaged SGD \citep{polyak1992}, SGD with momentum \citep{Momentum_Polyak, Momentum_DNN} and Adam~\citep{Adam}. The results show that  $1/t$ up-down policy also works for averaged SGD and momentum. Moreover, we compare the proposed step-size strategies to other popular step sizes, e.g., triangular policy \citep{smith2017cyclical} and cosine annealing \citep{loshchilov2016sgdr}. A great potential is shown when the step size is created based on the bandwidth, especially for nonconvex optimization.  

\vspace{0.1cm}

\noindent{\bf Organization:}  in Section \ref{sec:2}, we present some necessary definitions and lemmas used in the downstream analysis. In Section \ref{sec:3}, we investigate the conditions for bandwidth-based step size of SGD to achieve the $\mathcal{O}(1/T)$ convergence rate. Section 
\ref{sec:4} discusses the scenario where the ending points of the bandwidth step size being in the same order which covers the most cases we met. Section \ref{sec:5} considers the situation where the bands have different lower and upper boundaries. In Section \ref{sec:6}, we perform some numerical experiments for the proposed step sizes based on bandwidth. Then we make a conclusion in Section \ref{sec:7}. 

\noindent{\bf Notation.}
Let $x^{\ast}$ be the unique minimum of $f$, that is $x^{\ast} = \argmin_{x\in \R^d} f(x)$. We use $\mathcal{F}_t$ to denote $\sigma$-algebra of the random information at $t$-th iteration.  In default, the expectation is taken  with respect to  the source distribution
$\Xi$, that is $\E[\cdot] = \E_{\Xi}[\cdot] := \E_{\xi \sim \Xi}[\cdot]$. Other notations include: $\left\|\cdot\right\|$:=  $\left\|\cdot \right\|_2 $; $[n] = \left\lbrace 1,2,\dots, n \right\rbrace $; $[n] \backslash [n_1]=\left\lbrace n_1+1, n_1+2, \dots, n \right\rbrace $ for any $n_1 < n \in \N$.

\section{Preliminaries}\label{sec:2}

 In this part, we will give some definitions and basic lemmas used in the later sections.

\begin{assumption}[$\mu$-strongly convex]\label{strongly_convex}
	The objective function $f(\cdot):\R^d \longmapsto \R$ is $\mu$-strongly convex if there exists a constant $\mu > 0$ such that 
	\begin{equation}\label{stronglyconvex}
	f(x) - f(\hat{x}) \geq \left\langle \nabla f(\hat{x}), x-\hat{x} \right\rangle + \frac{\mu}{2} \left\|x-\hat{x} \right\|^2, 
	\end{equation} 
	for all $x, \hat{x} \in \R^d$.
\end{assumption}
Note that  $f(x;\xi)$ for each $\xi$ is not guaranteed to be convex even we assume that $f(x)$ is $\mu$-strongly convex.
.

\begin{assumption}[Expected smoothness]\label{expected_smooth}
	There exists a constant $L_f > 0$ such that
	\[
	\E[ \left\| \nabla f(x; \xi) -\nabla f(x^{\ast};\xi) \right\|^2] \leq 2L_f(f(x) - f(x^\ast)).
	\]
	Let $\E[\left\| \nabla f(x^{\ast} ;\xi)\right\|^2 ] = \sigma^2$, where $\sigma$ is a finite constant, we have 
\begin{equation}
\E[ \left\| \nabla f(x; \xi) \right\|^2] \leq  4L_f(f(x) - f(x^\ast)) + 2\sigma^2.
\end{equation}
\end{assumption}

\noindent{\bf Uniformly bounded gradient.} The assumption of uniformly bounded gradient (i.e., $\E[\left\|g_t \right\|^2 ] \leq G^2$ for some fixed $G > 0$) is used in some recent papers~\citep{Moulines-Bach2011,Shamir-Zhang2013,Rakhlin-Shamir-Sridharan2011,hazan2014beyond,jain2019making}. However,  this is clearly false if
$f$ is strongly convex, which has been pointed out by \citet{nguyen2018new}. If $f$ is $\mu$-strongly convex and $\E[\left\|g_t \right\|^2 ] \leq G^2$, by \emph{Jensen inequality} in expectation that $\left\|\E[X]\right\|^2 \leq \E[\left\|X\right\|^2]$, we have
\begin{equation*}
\mu^2 \left\|x-x^{\ast}\right\|^2 \leq 2\mu (f(x) - f(x^{\ast})) \leq \left\|\nabla f(x) \right\|^2 = \left\|\E[\nabla f(x; \xi)] \right\|^2 \leq \E[\left\|\nabla f(x;\xi) \right\|^2] \leq G^2.
\end{equation*}
In this case,  $f(x) - f(x^\ast)$ and  $\left\|x-x^{\ast}\right\|^2$ should be bounded on the whole space $\R^d$. However,  this  leads to a contradiction when $\left\|x_t -x^{\ast} \right\| $ is sufficiently large. Thus we assume expected smoothness \citep{gower2018stochastic} rather than uniformly bounded gradient.

\noindent{\bf $L$-smooth property vs expected smoothness.} Suppose that $f$ is $\mu$-strongly convex. By (\ref{dis_L}), the $L$-smooth property used in \citet{Moulines-Bach2011} 
\begin{equation}\label{dis_L}
\left\|\nabla f(x; \xi) - \nabla f(x^{\ast};\xi)\right\|^2 \leq L^2\left\|x-x^{\ast} \right\|^2 \leq \frac{2L^2}{\mu}[f(x) - f^{\ast}],
\end{equation}
implies expected smoothness  with $L_f =L^2/\mu$, but the opposite does not hold (see \citet{nguyen2018new}). Moreover, if $f$ is convex and $L$-smooth, the expected smoothness assumption can be satisfied with $L_f = 2L$ but the opposite is not true. Indeed, the example 2.2 of  \citet{Gower-etal.2019} shows that Assumption \ref{expected_smooth} holds even when $f(x;\xi)$ or $f$ is not convex.

\begin{lem}\label{sec:lem1}
Let a constant $\tau \in [1, 2)$. We assume that $f$ is $\mu$-strongly convex, then
	\begin{equation}\label{lem1:inequ1}
	\left\langle \nabla f(x), x - x^{\ast} \right\rangle 
	\geq (2-\tau)(f(x) - f(x^{\ast}))  + \frac{\tau\mu}{2}\left\|x -x^{\ast} \right\|^2, \,\text{for}\,\, x \in \R^d.
	\end{equation} 
\end{lem}

In Lemma \ref{sec:lem1}, the constant $\tau \in [1,2)$ is introduced to balance the weights of $f(x) - f(x^{\ast})$ and $\left\|x -x^{\ast} \right\|^2$. Let $x=x_t$ in (\ref{lem1:inequ1}). The term $(2-\tau)(f(x_t) - f(x^{\ast}))$ is used to eliminate $4L_f[f(x_t) - f(x^{\ast})]$ which is introduced by the expected smoothness assumption. All proofs of the lemmas in this section are provided in Appendix {\bf A}.

\begin{lem}\label{sec:lem2}
	Suppose that  Assumptions \ref{strongly_convex} and  \ref{expected_smooth} hold. Considering the mini-batch version of the SGD method, we have $\E[  \left\|x_{t+1}-x^{\ast} \right\|^2 \mid \mathcal{F}_t] $ is at most
\begin{equation}\label{lem2_equ1}
 (1- \tau\mu\eta(t))\left\|x_t -x^{\ast} \right\|^2 + 2\eta(t)^2 \sigma^2 + (4L_f\eta(t)^2 - 2(2-\tau)\eta(t))[f(x_t) - f(x^{\ast})].
\end{equation}
Besides, let $n_0 := \sup \left\lbrace t \in \N^{+}: \eta(t) >  \frac{(2-\tau)}{2L_f} \right\rbrace $ and $f_{n_0} := \mathop{\max}\limits_{1 \leq t \leq n_0}\left\lbrace f(x_t) - f(x^{\ast})\right\rbrace $. If  $n_0$ is a constant independent of $T$ (the budget of  the iteration $t$), then for $t > n_0$, we have $\E[  \left\|x_{t+1}-x^{\ast} \right\|^2] $ is at most
\begin{equation}\label{lem2_equ2}
\exp\left( -\tau\mu\sum_{l=1}^{t}\eta(l)\right) \Delta_{n_0}^0   + 2\sigma^2\sum_{l=1}^{t}\eta(l)^2\exp\left( -\tau\mu\sum_{u>l}^{t}\eta(u)\right),
\end{equation}
where $\Delta_{n_0}^0 = \left\|x_1- x^{\ast} \right\|^2 + \frac{n_0 \chi_{n_0}f_{n_0}}{\exp(-\tau\mu \sum_{l=1}^{n_0}\eta(l))}$ and $\chi_{n_0}=\mathop{\max}\limits_{1 \leq t \leq n_0}\left\lbrace 4L_f\eta(t)^2 - 2(2-\tau)\eta(t)\right\rbrace $.
\end{lem}

In order to estimate the upper bound (\ref{lem2_equ2}) of $\E[  \left\|x_{T+1}-x^{\ast} \right\|^2] $, $n_0$ needs to be  independent of the iterates budget $T$.  This can be achieved by requiring the step size $\eta(t)$ smaller than $(2-\tau)/2L_f$ after $n_0$ iterations. For the commonly used step size $\eta(t) = \eta_0/t^p$  ($p\in (0,1]$) which finally decreases to zero, this is also easily satisfied when $n_0 \leq (2\eta_0L_f/(2-\tau))^{1/p}$.

\begin{rem}
For simplicity, we might as well let
\begin{subequations}
\begin{equation}\label{gamma1} 
\Gamma_{T}^1 := \exp\left( -\tau\mu\sum_{l=1}^{T}\eta(l)\right)\Delta_{n_0}^0, 
\end{equation}
\begin{equation}\label{gamma2}
 \Gamma_{T}^2   :=2 \sigma^2\sum_{l=1}^{T}\eta(l)^2\exp\left( -\tau\mu\sum_{u>l}^{T}\eta(u)\right). 
\end{equation}
\end{subequations}
From Lemma \ref{sec:lem2}, let $t = T$, we have
\begin{equation}\label{main_inequ1}
\E[  \left\|x_{T+1}-x^{\ast} \right\|^2 ] \leq  \Gamma_{T}^1 +  \Gamma_{T}^2.
\end{equation}
Based on (\ref{main_inequ1}), the upper bound of $\E[  \left\|x_{T+1}-x^{\ast} \right\|^2 ] $ is divided into two parts $ \Gamma_{T}^1$ and $\Gamma_{T}^2$. Once the summation $\sum_{l=1}^{T}\eta(l)$ is determined, $\Gamma_{T}^1$ can be estimated by (\ref{gamma1}). Therefore, the challenge of the following analysis falls on the evaluation of $\Gamma_{T}^2$.
\end{rem}

\section{Non-Asymptotic Analysis of SGD for An Optimal Convergence Rate }\label{sec:3}

In this section, we will analyze the non-asymptotic convergence rate of the classical SGD algorithm where the step size $\eta(t)$ satisfies the following conditions
\begin{enumerate}
	\item[${(A)}$] there exists a constant $m > 0 $ such that $\eta(t) \geq \frac{m}{t}$,
\item[${(B)}$] there exists a constant $ M \geq m $ such that $\eta(t) \leq \frac{M}{t}$.
	\end{enumerate}
This is a special case of (\ref{eta_condition}) with $\delta_1(t) = \delta_2(t) = 1/t$.
The step size under these conditions is more general and possibly non-monotonic compared with the common choice $\eta(t) = \eta_0/(a + t)$ \citep{Rakhlin-Shamir-Sridharan2011,Moulines-Bach2011, Shamir-Zhang2013, Lacoste-Schmidt-Bach2012, Bottou-Curtis-Nocedal2018, Gower-etal.2019}.

The natural questions arising in this setting are the convergence of  SGD  and, if the convergence holds, the corresponding convergence rate (e.g., $\mathcal{O}(1/T)$ rate).  It is easy to see that SGD converges under $(A)$ and $(B)$ since they satisfy the well-known conditions ${(1')}$ $\sum_{t=1}^{\infty} \eta(t) = \infty$ and ${(2')}$ $\sum_{t=1}^{\infty} \eta(t)^2 < \infty$ given by \citet{SGD-1951}. The remaining question is which cases can ensure that SGD obtains the optimal $\mathcal{O}(1/T)$ convergence rate under the bandwidth condition (\ref{eta_condition}). All 
proofs in this section are given in Appendix {\bf B}.

\begin{thm}\label{sec:thm1}
	Let Assumptions \ref{strongly_convex} and \ref{expected_smooth} hold. We consider the step size $\eta(t)$ satisfy the conditions (A) and (B) for all $ 1 \leq t\leq T$ and let $n_0 := \sup \left\lbrace t \in \N^{+}: \eta(t) > \frac{(2-\tau)}{2L_f} \right\rbrace $. After at most $T >  n_0$ iterates, we have
	\begin{equation*}
	\E[\left\|x_{T+1} -x^{\ast} \right\|^2 ] \leq 
	\left\{\begin{array}{lr}
	\frac{\Delta_{n_0}^0 }{(T+1)^{\tau\mu m}} + 2\sigma^2\exp(1) M^2 \frac{\ln(T) +1}{T+1} & \,\text{if}\,\, m = \frac{1}{\tau\mu};\\
	\frac{\Delta_{n_0}^0 }{(T+1)^{\tau\mu m}} + \frac{\varepsilon_1M^2}{(\tau\mu m -1)} \frac{(T+1)^{(\tau\mu m -1)} + \tau\mu m -2 }{(T+1)^{\tau\mu m}} &  \,\text{else}\,\, m \neq \frac{1}{\tau\mu},
	\end{array}
	\right.
	\end{equation*}
	where $\tau \in (1,2]$ is a constant, $\varepsilon_1=2\sigma^2\exp(\tau\mu m)$ and $\Delta_{n_0}^0 $ has the same definition as Lemma \ref{sec:lem2}. 
\end{thm}

From condition $(B)$,  it is easy to see that $n_0 \leq 2ML_f/(2-\tau)+1$ which is independent of the iteration budget $T$.

\begin{rem}[The parameter $\tau$]
The parameter $\tau$ equals 2 in \citet{leen1994optimal} and \citet{leen1998two}, which reveals the connections of the initial step size $\eta_0$ and the strongly convex parameter $\mu$ with its convergence rate around the local minimizer $x^{\ast}$. However, in the non-asymptotic setting such as \citet{nguyen2018new} and \citet{Gower-etal.2019}, the value of $\tau$ is 1 and it  is hardly improved to be 2 under Assumption \ref{expected_smooth}.  By introducing $\tau \in [1,2)$ in Lemma  \ref{sec:lem1}, we observe the similarly asymptotic state around $x^{\ast}$ in Theorem \ref{sec:thm1}  when $\tau$  approaches $2$.

\end{rem}

\begin{cor}\label{sec:coro1}
	Under the same conditions as Theorem \ref{sec:thm1}, we can achieve the following results
	\begin{equation*}
	\E[\left\|x_{T+1} -x^{\ast} \right\|^2 ] \leq 
	\left\{\begin{array}{lr}
	\vspace{0.5em}
	\frac{\Delta_{n_0}^0  + (\frac{2-\tau\mu m}{1-\tau\mu m}) \varepsilon_1M^2}{(T+1)^{(\tau\mu m)}}  &  \,\text{if}\,\, m < \frac{1}{\tau\mu};\\
	\vspace{0.5em}
	\frac{\Delta_{n_0}^0  + 2\sigma^2M^2\exp(1)}{T+1} + \frac{2\sigma^2M^2\exp(1)\ln(T)}{T+1} & \,\text{else if}\,\, m = \frac{1}{\tau\mu};\\
	\vspace{0.5em}
	\frac{\Delta_{n_0}^0  + \varepsilon_1M^2}{(T+1)^{(\tau\mu m)}} +\frac{\varepsilon_1M^2}{(\tau\mu m -1)} \cdot\frac{1}{T+1}	&  \,\,\text{else}\, m > \frac{1}{\tau \mu}. \\
	\end{array}
	\right.
	\end{equation*}
\end{cor}

The corollary reveals the variation of the convergence rates with the coefficient $m$ of the lower bound $\delta_1(t)$. When $ m > 1/(\tau\mu)$, an optimal $O(1/T)$ convergence rate is attained measured by $\E[\left\|x_{T+1}-x^{\ast} \right\|^2 ]$. Note that $m = 1/(\tau\mu) $ is a special case which achieves a nearly-optimal $\mathcal{O}(\ln(T)/T)$ convergence rate. Besides, we can see that if $m$ is very small, it greatly slows down the convergence of SGD with the rate $\mathcal{O}(1/T^{(\tau\mu m)})$. Thus the value of $m$ is critical. The similar behaviors are observed in \citet{leen1994optimal}, \citet{SGD-complex} and \citet{Moulines-Bach2011}  for $\eta(t) = \eta_0 / t$.

\begin{thm}\label{sec:thm2_average}
	Let Assumptions \ref{strongly_convex} and \ref{expected_smooth} hold. We consider the step size $\eta(t)$ satisfy the conditions (A) and (B) for all $ 1 \leq t\leq T$.  Let $n_1 := \sup \left\lbrace t \in \N^{+}: \eta(t) >  \frac{2-\tau}{4L_f} \right\rbrace $ and $f_{n_1} = \max_{1 \leq t \leq n_1}\left\lbrace f(x_t) - f(x^{\ast})\right\rbrace$. If $\tau\mu m \geq 1$, for $ T > n_1$, we have that $\E[f(\hat{x}_T) - f(x^{\ast})] $ is bounded by	
\begin{equation}\label{func_inequ1}
\frac{1}{(2-\tau)mS_1}\left[\upsilon_1\Delta_{n_1}^0 + \upsilon_2(1-\frac{\tau}{2})mf_{n_1}+2 \sigma^2M^2(T-n_1+t_0\ln(T/n_1))\right],
\end{equation}
where $\hat{x}_T = \frac{\sum_{t=1}^{T}(t + t_0)x_t}{S_1}$, $S_1 = \frac{T(T+t_0)(t_0+1)}{2}$, $\Delta_{n_1}^0 = \frac{\left\|x_1 -x^{\ast} \right\|^2}{(n_1+1)^{\tau\mu m} } + 4\sigma^2M^2 + n_1\chi_{n_1}f_{n_1}$,   $\upsilon_1=(n_1+t_0+1)\left(n_1+1 - \tau\mu m\right)$  and  $\upsilon_2=(1+t_0)(n_1+t_0)$.  
\end{thm}

Moreover, we derive the error bound on the functions values of order $\mathcal{O}(1/T + \ln(T)/T^2)$. The convergence rate is comparable to those in the existing literature \citep{Rakhlin-Shamir-Sridharan2011,Lacoste-Schmidt-Bach2012,Shamir-Zhang2013,hazan2014beyond}. From (\ref{func_inequ1}), we know that  the convergence rate depends on $M^2/m$. If $M \approx m$, compared to Theorem \ref{sec:thm2_average}, the averaging technique weakens the effects of $m$ and $M$. 

\begin{rem}[Other averaging techniques]
In (\ref{func_inequ1}), for any $T > 0$, let $\hat{x}_T = \sum_{t=1}^{T}\alpha(t) x_t$, where $\alpha(t) = (t+t_0)/S_1$, we have $$\frac{\alpha(t)}{\alpha(t+1)} = \frac{t+t_0}{t+t_0+1}.$$ If $t_0 = 1$, the weight scheme in (\ref{func_inequ1}) is exactly the same as \citet{Lacoste-Schmidt-Bach2012}. For different $t_0 > 1$,  $\hat{x}_T$ produces a generalized weighted average iterate, different from those in \citet{Lacoste-Schmidt-Bach2012} and \citet{Shamir-Zhang2013}. We can see that for fixed $ 0 < t < T$, the ratio between the weights $\alpha(t)/\alpha(t+1) = t/(t+\eta) $ \citep{Shamir-Zhang2013} is smaller than $(t+t_0)/(t+t_0+1)$ if $\eta \geq  1$ and $t_0 \geq 1$. This means that the weight of (\ref{func_inequ1}) from $t$ to $t-1$ decays slower than that in \citet{Shamir-Zhang2013}. Moreover, if $\alpha(t) = (t+t_0)^k/\sum_{t}(t+t_0)^k$ for some $k\in \N^{+}$,  we have 
$$\frac{\alpha(t)}{\alpha(t+1)} = \frac{(t+t_0)^k}{(t+t_0+1)^k}.$$
This form is actually equivalent to that of \citet{Shamir-Zhang2013} and the integer $k$  corresponds to $\eta$. 
\end{rem}

We further relax the lower or upper bound of $\eta(t)$ and figure out in which cases  the optimal $\mathcal{O}(1/T)$ convergence rate can also be obtained. To better understand how the lower or upper bound affects the convergence rate, we only change one of them at a time. In general, if we fix the upper bound $\delta_2(t)$, the lower bound of $
\eta(t)$ can be extended to ($A_1$) (see Theorem \ref{sec:thm3}). Moreover,  in Remark \ref{rem:2} we reveal that the condition $(A_1)$ is essential to reach the optimal $\mathcal{O}(1/T)$ convergence rate.

\begin{thm}\label{sec:thm3}
	Suppose that Assumptions \ref{strongly_convex} and \ref{expected_smooth} hold. We consider the step size $\eta(t)$ satisfy the following conditions
	\begin{enumerate}
		\item[${(A_1)}$]\label{condition_31} there exists a constant $C > 0$ such that for all $t^{\ast} \in \left\lbrace 1,2, \cdots, T \right\rbrace $, we have 
		\begin{equation}\label{sum_C}
		\sum_{t=t^{\ast}}^{T} \eta(t) \geq C\ln\left(\frac{T+1}{t^{\ast}}\right);
		\end{equation}
		\item[${(B)}$]\label{condtion_32} there exists a constant $M > 0$ such that $ \eta(t) \leq \frac{M}{t}$ for all $ 1 \leq  t \leq T$.		
	\end{enumerate}
Let $n_0 := \sup \left\lbrace t \in \N^{+}: \eta(t) > \frac{(2-\tau)}{2L_f} \right\rbrace $.	If $C > \frac{1}{\tau\mu}$, for $t > n_0$, we have $\E[\left\|x_{T+1} -x^{\ast} \right\|^2 ]$ is at most
	\[
	\frac{ \Delta_{n_0}^0 + 8\sigma^2M^2}{(T+1)^{(\tau\mu C)}} +  \frac{8\sigma^2M^2 \exp(1)}{(\tau\mu C-1)}\cdot\frac{1}{T+1}.
	\]		
\end{thm}

The theorem shows that if the upper bound $\delta_2(t)$ is of order $1/t$, the lower bound of $\eta(t)$ can be extended to be of order $1/t$ on average to obtain an $\mathcal{O}(1/T)$ rate. Note that condition $(A_1)$ does not require $\eta(t)$ to be larger than $C/t$ for all $1 \leq t \leq T$.  For example, if $\eta(t) $ is larger than $m/t$ for $t \in [1, \alpha T]$ where $\alpha \in (0,1)$ and satisfies condition $(B)$,  we still can derive an $\mathcal{O}(1/T)$ bound. 

\noindent{\bf Compared to the step size proposed by \citet{jain2019making}.}
The following piece-wise decay step size which is modified by \citet{jain2019making} for strongly convex problems (see (4) of \citet{jain2019making}) 
\begin{equation*}
\eta(t) = 2^{-i} \cdot\frac{1}{\mu t }, \,\text{for }\, T_i <t  \leq T_{i+1}, \, T_{i} = T - \lceil T\cdot 2^{-i}\rceil,
\end{equation*}
 satisfies $(A_1)$ and $(B)$. From Theorem \ref{sec:thm3}, we are able to achieve a slightly weaker result than that of \citet{jain2019making} measured by functions values on the last iterate. However, it hardly gives the error bound of \citet{jain2019making} under such general conditions $(A_1)$ and $(B)$ without any modification. To keep our focus, we will not give further analysis.

\begin{rem}\label{rem:2}	
	In order to analyze the convergence rate of SGD, the key step is to estimate $\Gamma_T^2 $ defined by (\ref{gamma2}). If $\eta(t)$ has an upper bound $M /t $ for all $1\leq t \leq T$, we have
	\[
	\Gamma_T^2 = 2\sigma^2\sum_{t=1}^{T} \eta(t)^2 \exp\left( -\tau\mu \sum_{u >t}^{T}\eta(u)\right)  \leq 2\sigma^2\sum_{t=1}^{T}\frac{1}{t^2}\cdot \exp\left( -\tau\mu \sum_{u >t}^{T}\eta(u)\right) .
	\]
Considering the partial summation of $\frac{1}{t^2} \exp\left( -\tau\mu \sum_{u >t}^{T}\eta(u)\right) $ from $t^{\ast}$ to T,  for all $ 1 \leq t^{\ast} \leq T$, we have
	\[
	\sum_{t=t^{\ast}}^{T} \frac{1}{t^2}\cdot \exp\left( -\tau\mu \sum_{u >t}^{T}\eta(u)\right)  \geq \sum_{t=t^{\ast}}^{T}\frac{1}{t^2}\cdot \exp\left( -\tau\mu \sum_{u=t^{\ast}}^{T}\eta(u)\right).
	\]
	In order to achieve the convergence rate such that $\E[\left\|x_{T+1} - x^{\ast} \right\|^2 ] \leq \mathcal{O}(1/T)$, we have to require that 
	\[
	\sum_{t=t^{\ast}}^{T}\frac{1}{t^2}\cdot \exp\left( -\tau\mu \sum_{u=t^{\ast}}^{T}\eta(u)\right)  \leq \mathcal{O}\left( \frac{1}{T}\right) .
	\]
Then
	\[
	\exp\left( -\tau\mu \sum_{u=t^{\ast}}^{T}\eta(u)\right)  \left( \frac{1}{t^{\ast}} - \frac{1}{T}\right)  \leq \mathcal{O}\left( \frac{1}{T}\right) 
	\Longrightarrow	\sum_{u=t^{\ast}}^{T}\eta(u) \geq \frac{1}{\tau\mu}\ln\left( \frac{T}{t^{\ast}} -1\right)  + \mathcal{O}(1).
	\]
	Thus we see that condition $(A_1)$ in Theorem \ref{sec:thm3} is essential to achieve the optimal $\mathcal{O}(1/T)$ convergence rate under condition $(B)$.
\end{rem}

As we know, large step sizes are often preferred in practice, especially at the initial training~\citep{huang2017snapshot,loshchilov2016sgdr}. In the following parts, we extend the upper bound of $\eta(t)$ and discuss the convergence rates of SGD with larger step sizes. We firstly study the case that $\delta_2(t)$ is of order $1/t^r$ ($r \in (\frac{1}{2}, 1)$) at the initial $C_1T^p$ iterations ($p \in (0,1)$). For simplicity, we assume that $C_1T^{p}$ is an integer. 

\begin{thm}\label{sec:opt_thm}
	Let Assumptions \ref{strongly_convex} and \ref{expected_smooth} hold. Let  $\eta(t)$ satisfies the following conditions
	\begin{enumerate}
			\item[${ (A)}$] there exists a constant $m > 0$ such that $\eta(t) \geq \frac{m}{t}$ for all $ 1 \leq t \leq T$;
		\item[${ (B_1)}$] There exist constants $p \in (0,1)$, $r \in (\frac{1}{2}, 1)$ and $C_1, M_1, M_2 > 0$ such that $\eta(t) \leq \frac{M_1}{t^{r}}$ for $ t \in  [C_1T^{p}]$ and $\eta(t) \leq \frac{M_2}{t}$ for $ t \in  [T]\backslash [C_1T^{p}] $.
\end{enumerate}
Let $n_0 := \sup \left\lbrace t \in \N^{+}: \eta(t) >  \frac{(2-\tau)}{2L_f} \right\rbrace $.	If $ m  > \frac{1}{\tau\mu}$, then for $t > n_0$, we have $\E[\left\|x_{T+1} - x^{\ast} \right\|^2 ]$ is at most
	\begin{align*}
 \frac{\Delta_{n_0}^0  + \varepsilon_1(M_1^2+M_2^2)}{(T+1)^{\tau\mu m}} + \frac{\varepsilon_1 M_1^2(C_1+1)^{\varsigma_2}}{\varsigma_2T^{\varsigma_1}}
	 + \frac{\varepsilon_1M_2^2}{(\tau\mu m -1)(T+1)},
	\end{align*}
	where $\varepsilon_1=2\sigma^2\exp(\tau\mu m)$, $ \varsigma_1 =  (1-p)\tau\mu m + p(2r-1)$ and $\varsigma_2 = 1-2r+\tau\mu m > 0$.
\end{thm}

From Theorem \ref{sec:opt_thm}, if $\varsigma_1 \geq 1$ and $\tau\mu m >1$, an $\mathcal{O}(1/T)$ convergence rate can be obtained. To ensure that $\varsigma_1 \geq 1$,  we have $0 < p \leq (\tau\mu m -1)/(\tau\mu m -2r +1) < 1 $ for $m > 1/(\tau\mu)$ and $r\in (1/2,1)$. So the value of $p$ is reasonable.  When $r$ is near to $1/2$, $(\tau\mu m -1)/(\tau\mu m -2r +1) \approx (\tau\mu m -1)/(\tau\mu m)$. In this case if $\tau\mu m$ is very large, $p$ is very close to 1. Alternatively, we can require $\tau\mu m \geq (1-p(2r-1))/(1-p) > 1$ for $p\in (0,1)$ and $r\in (1/2,1)$,  which is stronger than the condition  of Theorem \ref{sec:thm1}.

In the following theorem, the step size $\eta(t)$ is allowed to vary within a band whose lower and upper bounds consist of two positive constants for $1\leq t \leq C_1T^p$ ($p \in (0,1)$) for an optimal $\mathcal{O}(1/T)$ convergence rate. 
\begin{thm}\label{sec2:thm6}
	We assume that Assumptions \ref{strongly_convex} and \ref{expected_smooth} hold. If the step size $\eta(t)$ satisfies the following conditions: there are some constants $p \in (0,1)$, $C_1 > 0$, $0 < m_1 \leq M_1$, $m_2 \leq M_2$ such that 
	\begin{enumerate}
			\item[${ (A_2)}$]  $ \eta(t) \geq  m_1$ for $t \in  [C_1T^{p}]$ and $\eta(t) \geq \frac{m_2}{t}$ for $t \in [T]\backslash [C_1T^{p}]$;
		\item[${ (B_2)}$]  $ \eta(t) \leq M_1$ for $ t \in  [C_1T^{p}]$ and $\eta(t) \leq \frac{M_2}{t}$ for $t \in [T]\backslash [C_1T^{p}]$.
\end{enumerate}
Let $n_0 := \sup \left\lbrace t \in \N^{+}: \eta(t) >  \frac{(2-\tau)}{2L_f} \right\rbrace $.	If $\kappa = (\tau\mu m_2)(1-p) \geq 1$ and $n_0$ is a constant that is independent of $T$, then for $t > n_0$, $\E[\left\|x_{T+1} - x^{\ast} \right\|^2 ]$ is at most
	\[
 \frac{\exp(\tau\mu m_2)}{\tau\mu m_1C_1}\cdot\frac{\Delta_{n_0}^0  }{T^{(\kappa+p)}} + 2\sigma^2\exp(\tau\mu m_2)\left(\frac{M_1^2 C_1^{(\tau\mu m_2)}}{\tau\mu m_1T^{\kappa}} + \frac{M_2^2}{(\tau\mu m_2 - 1)}\cdot \frac{1}{T+1}\right).
	\]
\end{thm}

Let $\kappa = (\tau\mu m_2)(1-p) \geq 1$, we have $\tau\mu m_2 \geq 1/(1-p) > 1$ ($p \in (0,1)$), which is stronger than the requirement on $\tau\mu m$ in Theorem \ref{sec:opt_thm}. 
When the iteration budget $T$ is very large, we see that $C_1T^p$ can also be very large. For $T  \gg  K$ where $K$ is the condition number,  our result possibly enlarges the existing result of \citet{Gower-etal.2019} equipped with a constant step size on the initial $4K$ steps. Note that \citet{allen2018make} proposes an algorithm $\text{SGD}^{sc}$(a.k.a. SGD after SGD), in which the step size $\eta(t) = 1/(2L)$ for the initial $\lfloor T/2 \rfloor$ iterates, where $L$ is the parameter of smoothness. However, the output of each inner loop is an average of all inner iterates, which is different from the classic SGD algorithm we focus on in this paper. Thus we do not give any further comparison.

\subsection{Discussions on Other Popular Step Sizes}\label{sec:3:1}
\citet{hazan2014beyond} proposed the following piece-wise decay step size within the $i$-th run
\begin{equation}\label{exp_decay1}
\eta(t) = \eta_i = \frac{\eta_{i-1}}{2}, \,t \in [T_i, T_{i+1}), \,T_{i+1} = 2T_i,
\end{equation}
where $\sum_{i} T_i = T$.  The above step size exponentially decays per cycle but the period of each cycle increases by a given factor. Clearly, it satisfies the conditions $(A)$ and $(B)$. So its convergence rates are guaranteed by our analysis (see Theorem \ref{sec:thm1} and \ref{sec:thm2_average}).  

The exponentially decaying step size is popular and often used  in deep learning, that is
\begin{equation}\label{exp_fix1}
\eta(t) = \eta_0 \alpha^{\lfloor t/T_0 \rfloor}, 
\end{equation}
where $\alpha \in (0,1)$ is a constant which is independent of $T$ and $T_0$ accounts for how many epochs have been performed since the last run. For simplicity, we let $\alpha =1/2$. If the period $T_0$ is the same at each cycle, we consider the following cases that
\begin{itemize}
	\item $T_0 = 1$, or  a constant (independent of $T$). Its non-asymptotic convergence can not be guaranteed because $\sum_{t=1}^{\infty}\eta(t) = + \infty$ is not satisfied. 
	\item $T_0 = \lfloor T^{r}\rfloor$, $ r \in (0,1)$. When $ k_0=\lfloor t/T_0 \rfloor = \lfloor r\log_2(T)\rfloor$,  the partial summation $\sum_{t=k_0T_0}^{(k_0+1)T_0}\eta(t)^2\exp(-\tau\mu \sum_{u=(k_0 +1)T_0}^{T}\eta(u)) \geq  \exp(-2\tau\mu \eta_0)\sum_{t=k_0T_0}^{(k_0+1)T_0}\eta(t)^2 = \mathcal{O}(1/T^r)$. From Lemma \ref{sec:lem2},   it hardly obtains the non-asymptotical  $\mathcal{O}(1/T)$ convergence rate.	
 \item $T_0 = \lfloor T/\log(T)\rfloor$. Let $k_0 = \lfloor t/T_0 \rfloor  = \lfloor\log_2(T) - \log_2\log_2(T)\rfloor$. In this case we have $\sum_{t=k_0T_0}^{(k_0+1)T_0}\eta(t)^2\exp(-\tau\mu \sum_{u=k_0T_0}^{T}\eta(u)) \geq \exp(-2\tau\mu \eta_0)\sum_{t=k_0T_0}^{(k_0+1)T_0}\eta(t)^2 = \mathcal{O}(\log_2(T)/T)$. We can see that the best result will not exceed $\mathcal{O}(\log_2(T)/T)$ from Lemma \ref{sec:lem2}. Such a rate has been demonstrated by \citet{ge2019step} for the least squares problems.
 \item $T_0 = \lfloor T/m \rfloor$, where $m \in \N^{+}$ is a constant which is independent of $T$. This case is more often met in practice. When $m$ is a positive constant,  the final step size is $2^{-m} \gg 1/T$. It is impossible to achieve the non-asymptotical  $\mathcal{O}(1/T)$ rate.	
\end{itemize}
From the above discussions, if $T_0$ is fixed and is the same per cycle, it hardly achieves the ideal $\mathcal{O}(1/T)$  convergence rate for strongly convex problems. 


In the cyclical step size schedules, for example: a triangular policy was proposed where the step size is locally  increased and then  decreased linearly within a band and  the schedule (\ref{exp_fix1}) is adopted as the baseline \citep{smith2017cyclical}; a sine-wave learning rate was proposed where the step size decays exponentially and  locally oscillations within a range of values \citep{sinewave}.  If the boundaries of the step size are based on (\ref{exp_fix1}), from the above analysis we know that the optimal $
\mathcal{O}(1/T)$ convergence rate can not be guaranteed in theory. However, to say the least, once the boundaries are taken as (\ref{exp_decay1}) or satisfy the situations we discussed, the optimal $\mathcal{O}(1/T)$ convergence rate can be obtained by our analysis.

\section{Convergence Analysis Under the Same Boundary Order}\label{sec:4}
In this section, we investigate the convergence rate of SGD where the bandwidth-based step size (\ref{eta_condition})  has  the same boundary order, i.e., $\delta_1(t) = \delta_2(t)$. 

The convergence conditions on step size for standard SGD were proposed by \citet{SGD-1951} 
\begin{equation}\label{stepsize_cond1}
 {(1')} \,\,\sum_{t=1}^{\infty} \eta(t) = +\infty; \quad \quad\quad\quad\quad\quad 
{(2')} \,\,\sum_{t=1}^{\infty} \eta(t)^2 < +\infty. \tag{H1}
\end{equation}
Obviously,  the common choice $\eta(t) = 1/t^p$ for $p\in (\frac{1}{2}, 1]$ satisfies (\ref{stepsize_cond1}). However, (\ref{stepsize_cond1}) does not hold for $\eta(t) = 1/t^p$ with $0 < p \leq 1/2$, which has been proven to converge \citep{leen1998two, ljung1977analysis, Moulines-Bach2011}. Moreover, the step size under (\ref{stepsize_cond1}) is possibly non-monotonic. For example, it can oscillate between two boundaries $\eta(t) = 1/t$ and $\eta(t) = 1/\sqrt{t}$.

\citet{ljung1977analysis} proposed the following convergence conditions for recursive stochastic algorithms
\begin{equation}\label{stepsize_cond2}
\begin{array}{l}
{(1')}\,\,\sum_{t=1}^{\infty} \eta(t)  = +\infty;\,\,\quad\quad \quad\quad\quad \quad \,\,\,  { (2')}\,\,\sum_{t=1}^{\infty} \eta(t)^p < +\infty, \, \text{for some}\,\,  p > 0; \quad \vspace{0.2cm}\\
{(3')} \,\,\eta(\cdot) \, \text{is a decreasing sequence};   \quad \quad {(4')} \,\,\lim_{t \rightarrow \infty} \sup[1/\eta(t) - 1/\eta(t-1)] < \infty. \tag{H2}
\end{array}
\end{equation}
Compared to (\ref{stepsize_cond1}), (\ref{stepsize_cond2}) can cover more generic cases, e.g., $\eta(t) = 1/t^p $ for all $p \in (0,1]$. However, for example $\eta(t) = 1/(t\log(t+1))$, which satisfies (\ref{stepsize_cond1}), is not admitted by (\ref{stepsize_cond2}). Moreover,  the step size $\eta(t)$ in (\ref{stepsize_cond2}) is assumed to decrease which is not essential for (\ref{stepsize_cond1}).

Recently, \citet{nguyen2018new} extended (\ref{stepsize_cond1}) and (\ref{stepsize_cond2}) to the following cases
\begin{equation}\label{stepsize_cond3}
 {(1')}\,\, \sum_{t=1}^{\infty}\eta(t) = +\infty; \quad \quad
\quad   {(2')} \,\,\lim_{t\rightarrow +\infty} \eta(t) = 0;  \quad\quad
\quad  {(3')}\,\,\frac{d\eta(t)}{dt} \leq 0. 
\tag{H3}
\end{equation}
The common choices $\eta(t) = 1/t^p$  for $p\in (0, 1]$ and $ 1/(t\ln(t))$ satisfy (\ref{stepsize_cond3}). In addition, we see that $\eta(t) = 1/\ln(t)$ which goes down slower than any polynomial, satisfies the above conditions.
They proved the convergence of SGD and derived a uniform formula to describe the convergence rates for the step sizes satisfying (\ref{stepsize_cond3}) (see theorem 9 and 10 in \citet{nguyen2018new}). 

In this section, we focus on the sequence of step size $\left\lbrace \eta(t) \right\rbrace $ that satisfies 
\begin{equation}\label{stepsize_cond4}
m \delta(t) \leq \eta(t) \leq  M \delta(t), \tag{BD-S}
\end{equation}
where $ m \leq M$ are two positive constants and the boundary function $\delta(t)$ satisfies (\ref{stepsize_cond3}).  The main theorem is presented as follows, which covers the general cases mentioned above. The proofs in this section are provided in Appendix {\bf C}.
\begin{thm}\label{sec:thm4}
	Suppose  Assumptions \ref{strongly_convex} and \ref{expected_smooth} hold. The step size sequence $\left\lbrace \eta(t)\right\rbrace $
	satisfies condition (\ref{stepsize_cond4}) and the boundary function $\delta(t)$ is differentiable and satisfies (\ref{stepsize_cond3}). Let $n_0 := \sup \left\lbrace t \in \N^{+}: \eta(t) >  \frac{(2-\tau)}{2L_f} \right\rbrace $ and we assume that  $n_0$ is a constant which is independent of $T$. For $t > n_0$, 
	\begin{enumerate}
		\item[{\bf 1.}] if $\lim_{t\rightarrow \infty} t\delta(t)=0$,  we have that $\E[ \left\|x_{t+1} - x^{\ast} \right\|^2] $ is at most
		\begin{align*}
\left(\Delta_{n_0}^0 + \varepsilon_2 \frac{\delta(1)^2(t_{\epsilon}-1) + 2\epsilon^2}{\exp\left( -\tau\mu m  \int_{u=1}^{t_{\epsilon}}\delta(u)du\right)}\right) \exp\left( -\tau\mu m\int_{u=1}^{t+1}\delta(u)du\right), 			
\end{align*}
where  $\epsilon$ and $t_{\epsilon}$ are constants appeared in the proof, $\varepsilon_2 =2\sigma^2M^2 \exp(\tau\mu m \delta(1))$.
\item[{\bf 2.}] If $\lim_{t\rightarrow \infty} t\delta(t)=1$, the results of Theorem \ref{sec:thm1} can be applied.

\item[{\bf 3.}] If $\lim_{t\rightarrow \infty} t\delta(t) = + \infty$ and there exist constants $c_1 \leq \frac{\tau\mu m}{2}$ and $T_M \in \N$ such that $- \frac{d\delta(t)}{dt} \leq c_1 \delta(t)^2$ for all $t\geq T_M$, then $\E[\left\|x_{t+1} -x^{\ast} \right\|^2]$ is at most
\begin{align*}
 \frac{\varepsilon_2}{\tau\mu m-c_1}
\delta(t+1) + \left( \Delta_{n_0}^0  +\frac{\varepsilon_2\delta(1)^2T_M}{\exp\left( -\tau\mu m  \int_{u=1}^{t_{M}}\delta(u)du\right)} \right)\exp\left( -\tau\mu m \int_{u=1}^{t+1}\delta(u) du\right),
\end{align*}
where $\varepsilon_2$ is the same as the first case.
\end{enumerate}	
\end{thm}

Corresponding to the limit of $\delta(t)t$, three different cases are analyzed in Theorem \ref{sec:thm4}. When $\lim_{t\rightarrow \infty} t\delta(t) = 0$, the result is new. It covers the case which drops faster than $1/t$, e.g., $\delta(t) = 1/(t\ln(t))$.  In the third case, we add a condition that $- d\delta(t)/dt \leq c_1 \delta(t)^2 $ ($\forall t \geq T_M$) for some $c_1$ and $T_M$. The common choices, e.g., $\delta(t) = 1/t^p$ for all $p\in (0, 1]$ and $\delta(t) = 1/\ln(t)$, all satisfy the condition. More cases such as $\delta(t) = \ln(t+1)/t^p$ for all $p\in (0,1]$ are also included in the discussions. In addition, we see that when $t$ is continuous, $(4')$ of (\ref{stepsize_cond2}) can be reformulated as $$\lim_{t \rightarrow \infty} \sup[1/\eta(t) - 1/\eta(t-1)] = \lim_{t \rightarrow \infty} \sup\left[\frac{\eta(t-1) - \eta(t)}{\eta(t)\eta(t-1)}\right] = \lim_{t \rightarrow \infty}\sup \frac{\frac{- d\eta(t)}{dt}}{\eta(t)^2} < +\infty.$$  
This implies that there exists a constant $c_1 > 0$ such that $ - d\eta(t)/dt \leq c_1 \eta(t)^2$ for sufficiently large $t$. However, $c_1$ is supposed to be smaller than $ \tau\mu m/2$ in the third case. The following lemma reveals that as long as such $c_1 > 0$ exists, there must be a constant $c_1 > 0$ such that $c_1 \leq \tau\mu m/2$.
\begin{lem}\label{sec:lem4}
	We suppose that $\lim_{t\rightarrow \infty} t\delta(t) = + \infty$. If there exist constants $c_1 > 0$ and $T_{M} \in \N^{+}$ such that $- \frac{d\delta(t)}{dt} \leq c_1 \delta(t)^2$ for all $t \geq T_M$, there must be such a constant $c_1$ that satisfies $c_1 \leq \frac{\tau\mu m}{2}$. 
\end{lem}
\begin{rem}
	Theorem \ref{sec:thm4} shows the convergence rate of SGD where the bandwidth-based step size satisfies (\ref{stepsize_cond4}). We emphasize that 
	\begin{enumerate}		
		\item 
In the proof of the third case,  an important step is to use integral $\int_{l=1}^{t}P(l)dl$ to evaluate the summation $\sum_{l=1}^{t}P(l)$ where $P(l) $ is the product of $\delta(l)^2$ and $\exp(-\tau\mu m \int_{u=l}^{t+1} \delta(u)du) $. Even though $\delta(l)$ is decreasing and $\exp(-\tau\mu m \int_{u=l}^{t+1} \delta(u)du)$ is increasing, there can be many possibilities for their product. \citet{nguyen2018new} considered three cases for the product that, e.g., decreases and then increases,  keeps on increasing or decreasing (see the proof of theorem 9 in \citet{nguyen2018new}). However, as we know  the product of $\delta(l)^2$ and $\exp(-\tau\mu m \int_{u=l}^{t+1} \delta(u)du) $ increases and then decreases in \citet{ge2019step}. In Theorem \ref{sec:thm4}, we add a condition $ -d\delta(t)/dt \leq c_1 \delta(t)^2$ to describe ``most general cases" mentioned in \citet{nguyen2018new} and make the proof more rigorous.
		\item Theorem \ref{sec:thm4} reveals  the convergence rate of SGD, which is totally determined by $\delta(t+1)$ or $\exp(-\tau\mu m \int_{u=1}^{t+1}\delta(u) du)$. Compared to that of \citet{nguyen2018new}, our result provides better upper bounds in many cases. For example, when $\eta(t) = 1/(t\ln(t))$, theorem 10 of \citet{nguyen2018new} no longer gives an upper bound but Theorem \ref{sec:thm4} shows that it is bounded  by $\exp(-\tau\mu m \int_{u=1}^{t+1}\delta(u) du)$. In the case that $\eta(t) = 1/\sqrt{t}$, the first term of the upper bound  in theorem 10 \citep{nguyen2018new} is actually larger than $\eta(t+1)$, which is worse than that of Theorem \ref{sec:thm4}.		
\item The step size $\eta(t)$ in Theorem \ref{sec:thm4} can be non-monotonic, rather than monotonic or given in monotonic forms (e.g., $\eta_0/T$ or $\eta_0/t^p$ for $p\in(0,1]$) in most of the literature analyzing the convergence rate of SGD \citep{Rakhlin-Shamir-Sridharan2011,Moulines-Bach2011, Shamir-Zhang2013, Lacoste-Schmidt-Bach2012, Bottou-Curtis-Nocedal2018, Gower-etal.2019, jain2019making}.
\end{enumerate}
\end{rem}

\section{Convergence Analysis based on the Different Boundary Orders }\label{sec:5}
In this section, we will analyze the convergence rate of SGD where the lower bound function $\delta_1(t) $ and the upper bound function $ \delta_2(t)$  are  in different orders. From Section \ref{sec:4}, if the lower and upper bounds of the step size $\eta(t)$ are in the same order, their convergence rate is actually consistent with their boundaries. In the following part,  we want to find out the convergence behaviors of SGD when the boundaries of the step size are in different orders.

Firstly, we are interested in the case $\delta_2(t) = \ln(t+1)/(t+1)$ which decays slower than the lower upper $\delta_1(t) = 1/(t+1)$. 
\begin{thm}\label{sec4:thm1}
Suppose that Assumptions \ref{strongly_convex} and \ref{expected_smooth} hold. Let the step size sequence $\left\lbrace \eta(t)\right\rbrace $  satisfy that
	\begin{equation*}
	\frac{m}{t+1}\leq \eta(t) \leq \frac{M\ln(t+1)}{t+1}, \,t \geq 1,
	\end{equation*}
	for $ 0 <  m \leq M$. Let $n_0 := \sup \left\lbrace t \in \N^{+}: \eta(t) >  \frac{(2-\tau)}{2L_f} \right\rbrace $. For $t > n_0$, we have
	\[ \E[ \left\|x_{T+1} -x^{\ast} \right\|^2] \leq \left\{
	\begin{array}{lr}
	\frac{2^{(\tau\mu m)}\Delta_{n_0}^0  + \sigma^2M^2 \exp(1)\ln2}{(T+2)} + \frac{\varepsilon_1M^2 \ln^3(T+2)}{3(T+2)}& \text{if\,} \,m = \frac{1}{\tau\mu}; \\ 
	\frac{2^{(\tau\mu m)}\Delta_{n_0}^0  + 2\varepsilon_1\nu_1M^2 }{(T+2)^{(\tau\mu m)}}& \text{else if} \,\,m <\frac{1}{\tau\mu};	\\
	\frac{2^{(\tau\mu m)}\Delta_{n_0}^0  + \varepsilon_1\nu_2M^2}{(T+2)^{(\tau\mu m)}} + \left[ \frac{\ln^2(T+2)}{(\tau\mu m -1)} + \frac{2}{(\tau\mu m-1)^3}\right]   \frac{\varepsilon_1M^2}{(T+2)}& \text{else\,}\, m > \frac{1}{\tau\mu},
	\end{array}
	\right.
	\]
	where $\varepsilon_1=2\sigma^2\exp(\tau\mu m)$, $\nu_1= \frac{\ln2}{2} + \frac{2+2\ln2+\ln^22}{(1-\tau\mu m)^3}$ and $ \nu_2=\frac{\ln2}{2} + \frac{2^{(\tau\mu m)}\ln2 }{(\tau\mu m -1)^2}$.	
\end{thm}

The theorem reveals that when $m > 1/(\tau\mu)$, SGD can achieve an $\mathcal{O}(\ln^2(T)/T)$ bound which is nearly optimal. The proofs in this section are given in Appendix {\bf D}. 

As we know, (\ref{stepsize_cond1}) is sufficient for the convergence of SGD, but the convergence rate under (\ref{stepsize_cond1}) is unknown yet. If we keep the lower bound $\delta_1(t) = 1/t$ and continue to extend the upper bound $\delta_2(t)$, what kinds of results will we get?  The following result gives an answer to this interesting question.
\begin{thm}\label{sec4:thm2}
	We assume that Assumptions \ref{strongly_convex} and \ref{expected_smooth} hold. If the step size $\eta(t)$ satisfies that
	\begin{equation}\label{L-U_inequ1}
	\frac{m}{t}\leq \eta(t) \leq \frac{M}{t^{\alpha}}
	\end{equation} 
	for $\alpha \in (1/2, 1]$. Let $n_0 := \sup \left\lbrace t \in \N^{+}: \eta(t) >  \frac{(2-\tau)}{2L_f} \right\rbrace $.
For $t > n_0$, we have 
	\begin{equation*}
	\E[ \left\|x_{T+1} -x^{\ast} \right\|^2]  \leq \left\{
	\begin{array}{lr}	 \frac{\Delta_{n_0}^0  + 2\sigma^2M^2\exp(2\alpha-1)(1+\ln(T+1))}{(T+1)^{(2\alpha-1)}} & \,\,\text{if} \,\, \tau\mu m = 2\alpha -1;\\
	\vspace{0.3em}
	\frac{\Delta_{n_0}^0  + \frac{\varepsilon_1M^2(\tau\mu m-2\alpha)}{\tau\mu m -2\alpha +1}}{(T+1)^{(\tau\mu m)}} + \frac{\varepsilon_1M^2}{(\tau\mu m - 2\alpha +1) }\frac{1}{(T+1)^{(2\alpha-1)} } 	& \text{else}\,\, \tau\mu m \neq  2\alpha -1,
	\end{array}
	\right.
	\end{equation*}	
	where $\varepsilon_1=2\sigma^2\exp(\tau\mu m)$ is the same as that of Theorem \ref{sec:thm1}.		
\end{thm}

In Theorem \ref{sec4:thm2}, the upper bound $\delta_2(t)$ in (\ref{L-U_inequ1}) is extended to $1/t^{\alpha}$ for any $\alpha \in (1/2, 1]$. It is easy to see that (\ref{stepsize_cond1}) holds for the step size $\eta(t)$ which satisfies (\ref{L-U_inequ1}). The corresponding convergence rate is $\mathcal{O}(1/(T+1)^{2\alpha-1})$ which is relied on $\alpha$ when $ \tau\mu m > 1$. Obviously, this result is worse than those achieved at its lower and upper bounds. Unfortunately,  at present we are not able to improve Theorem \ref{sec4:thm2}. On the other direction, we reduce the lower bound $\delta_1(t)$ to $1/((t+1)\ln(t+1))$, which decreases faster than the case $\delta_1(t) = 1/t$ in Theorem \ref{sec4:thm2}.

\begin{thm}\label{sec4:thm3}
	Suppose that Assumptions \ref{strongly_convex} and \ref{expected_smooth} hold. Let the step size $\eta(t)$ satisfy
	\begin{equation}\label{step_cond9}
	\frac{m}{(t+1)\ln(t+1)}\leq \eta(t) \leq \frac{M}{(t+1)^{\alpha}}, \,\, t \geq 1,
	\end{equation} 
	for $\alpha \in (1/2, 1]$. 
	Then for sufficiently large $t$, there must be a constant $C_2 > 0$ such that
	\begin{equation*}
	\E[ \left\|x_{T+1} -x^{\ast} \right\|^2]   \leq  \frac{C_2}{(\ln(T+2))^{(\tau\mu m)}}.
	\end{equation*}

\end{thm}

 Theorem \ref{sec4:thm3} shows that the convergence rate of SGD where the step size satisfies (\ref{step_cond9}) is consistent with the result achieved at the lower boundary $\eta(t) = m/((t+1)\ln(t+1))$.

\section{Numerical Experiments}
\label{sec:6}%
In this section, we propose several non-monotonic step sizes within $1/t$-band to show the effectiveness compared to their baselines, e.g., $\eta(t) = \eta_0/t$ (called ${1/t}${\bf-stepsize}) and exponentially decaying step size.  The $1/t$ step size decays very fast, so we update all the step sizes after one epoch shown as Algorithm \ref{alg:1} (called Epoch-SGD). 
	\begin{algorithm}[H]
	\caption{Epoch-SGD}\label{alg:1}
	\begin{algorithmic}[1]			
	\STATE \textbf{Initialization}: initial point $x_0 = x_1^1$, \# inner loop $m^{'}$, \# outer loop $N$
	\FOR{ $t =  1: N$ }
	\STATE Update the step size $\eta(t)$	
	\FOR{ $i = 1: m^{'}$}
	\STATE Choose a subset $\Omega_i \subseteq [n]$ randomly, where $|\Omega_i|=b$ \\
	\STATE Compute $g_i^t = \frac{1}{b}\sum_{l \in \Omega_i}\nabla f(x_i^t;\xi_l)$\\ 
	\STATE $x_{i+1}^t = x_i^t - \eta(t) g_i^t$
	\ENDFOR
	\STATE $x_1^{t+1} = x_{m^{'}+1}^t$
	\ENDFOR
	\STATE {\bf Return} $x_{m^{'}+1}^{N}$
	\end{algorithmic}
	\end{algorithm}

  
  \subsection{The $1/t$-band Step Sizes}\label{sec:step}
We formulate some non-monotonic step sizes  $\eta(t)$ which belongs to a banded region $[\eta_0/t, s\eta_0/t]$ (named $\bm{1/t}${\bf -band}), where $s>1$. The boundary function $\eta(t) = \eta_0/t$ is called $1/t$-stepsize. Let $ t_i \,   (i=1,2,\cdots, 1 \leq t_1 < t_2 < t_3 <\cdots)$ be the nodes where the step size might be non-monotonic or non-differentiable. For $t \in [t_i, t_{i+1})$, let
\begin{equation}\label{t1t2_eta}
\eta(t) = \frac{\hat{A}_i}{\hat{B}_it+1},
\end{equation}
 where $\hat{A}_i, \hat{B}_i$ are constants such that $\eta(t_i) = s\eta_0 / t_i$ and $\eta(t_{i+1}) = \eta_0/t_{i+1}$. In reality, there exist other forms of $\eta(t)$, e.g., linear decay and concave decay. In the paper we are interested in the case that $\eta(t)$ has the form of (\ref{t1t2_eta}). 
We consider the two cases: {\bf (1)} $t_{i+1}-t_{i}$ is fixed. We call this {$\bm{1/t}$ {\bf Fix-period band}; {\bf (2)} $t_{i+1}-t_{i}$ grows exponentially. We call this $\bm{1/t}$ {\bf Grow-period band}. For intuitive explanation, we plot the two cases and their boundaries $1/t$-stepsize ($s=3$) in Figure  \ref{stepsize:1}.  
 \begin{figure}[h]
\centering
	\subfigure[]{\label{stepsize:1}
	\includegraphics[width=0.47\textwidth]{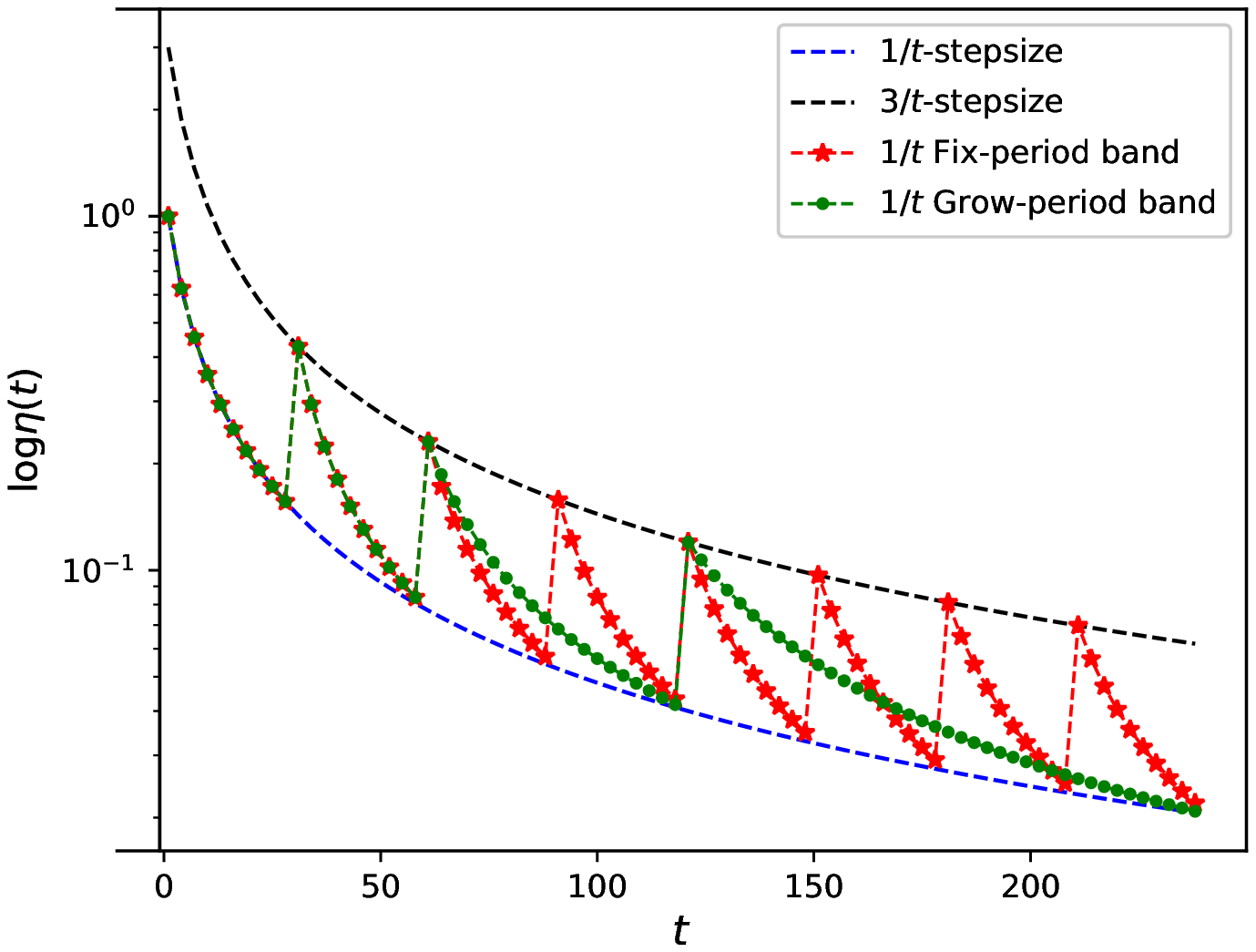}
	}
\subfigure[]{\label{stepsize:2}
\includegraphics[width=0.47\textwidth]{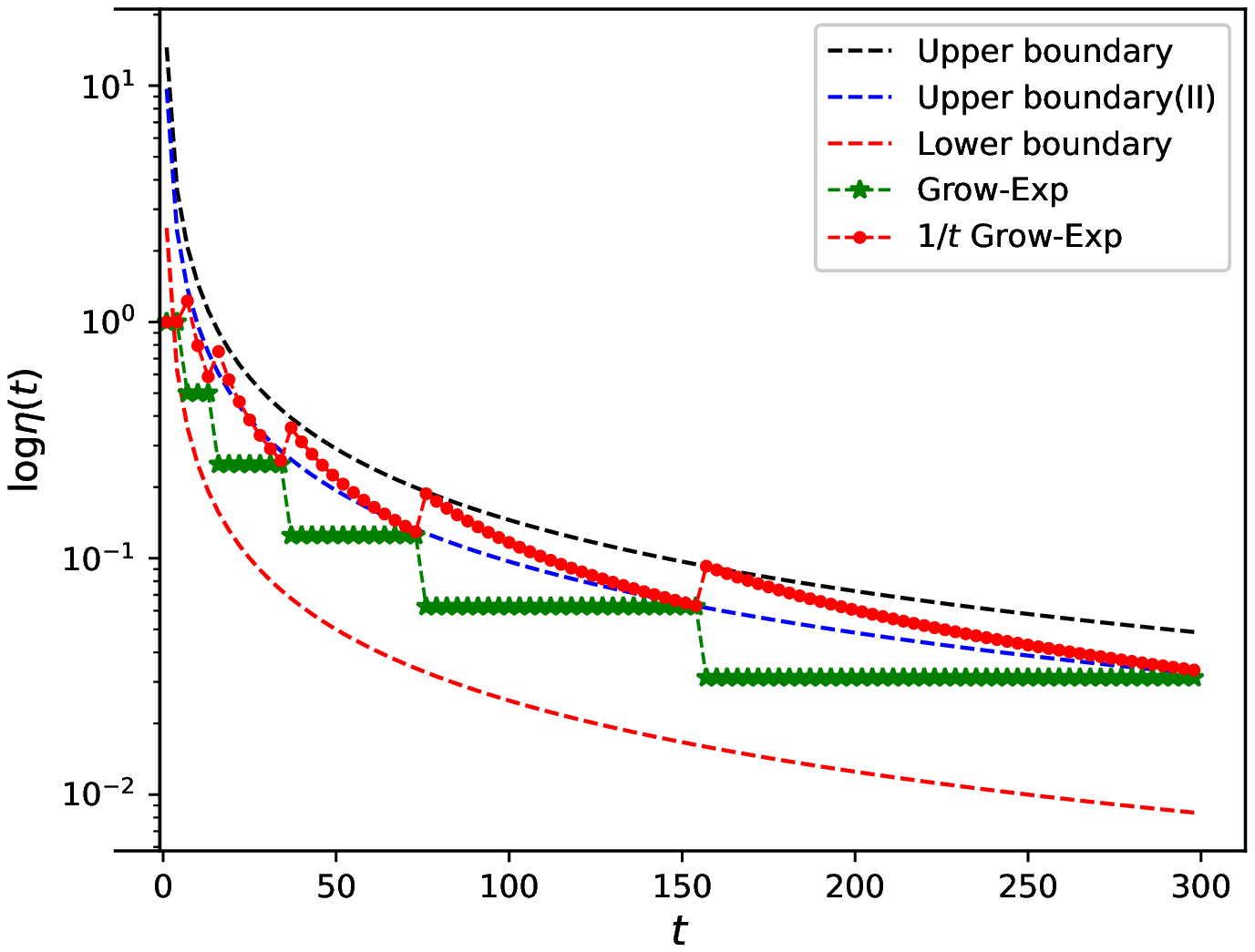}
}
\caption{Different kinds of $1/t$-band step sizes}
\label{fig:1}      
\end{figure}
	
	More general, the step size varies between the minimum $\eta_{\min} =  \lbrace{\eta_{\min}^i\rbrace}_{i\in\N}$ and maximum $\eta_{\max} = \lbrace{\eta_{\max}^i\rbrace}_{i \in \N}$, and locally has the form that
	\begin{equation}
	\eta(t) = \frac{\hat{A}_i}{\hat{B}_it + 1} \in [\eta_{\min}^i, \eta_{\max}^i],  t \in [t_i, t_{i+1}].
	\end{equation}
Especially, we consider $\eta_{\max}^i > \eta_{\min}^{i-1} $, which is called $\bm{1/t}$  {\bf up-down policy}. For $1/t$ Fix-period band and $1/t$ Grow-period band, the baseline  of the step size is  $\eta_{\min} = \eta_0/t$. Based on the known exponentially decaying step size with a growing period (called {\bf Grow-Exp})
\begin{equation}\label{exp_decay}
\eta(t) = \eta_i = \eta_0/2^{i}, t \in [t_i, t_{i+1}],\,\,  T_i = t_{i+1} - t_i= T_02^i,
\end{equation}
which has been used by \citet{hazan2014beyond}.		
Let $\eta_{\min}^i = \eta_i $ in (\ref{exp_decay}) and we define $\eta_{\max}^{i} = \theta \eta_{\min}^{i-1}$ where the up-down ratio $\theta > 1$ (called  $\bm{1/t }$ {\bf Grow-Exp}).  If $\theta$ is too large, a sudden increase in step size might lead to a very negative effect. Therefore, we restrict the ratio $\theta \in (1, 1.5]$. The Grow-Exp step size, $1/t$ Grow-Exp step size and their boundaries are plotted in Figure \ref{stepsize:2} where $T_0=5$ and $\theta=1.5$.  Regardless of Grow-Exp  or $1/t$ Grow-Exp,  we can easily  find that they all belong to $1/t$-band.

\begin{rem}
From Figure \ref{stepsize:1}, we see that the area enclosed by $1/t$ Fixed-period band and $x$-axis  is larger than that of its lower boundary. According to Lemma \ref{sec:lem2}, based on $1/t$ Fixed-period band, we can achieves a lower error bound for $\Gamma_{T}^1$ than that of the boundary $\eta(t)=\eta_0/t$. Thus $1/t$ Fixed-period band could be faster than $1/t$-stepsize ($\eta(t) = \eta_0/t$) at the initial iterations when $\Gamma_{T}^1$ is dominated the error bound of $\E[\left\|x_{T+1}-x^{\ast}\right\|^2]$. We have the similar conclusions for $1/t$ Grow-period band and $1/t$ Grow-Exp.

\end{rem}

 Next some numerical experiments are performed to demonstrate the efficiency of the proposed non-monotonic step sizes. All experiments are implemented in python 3.7.0 on a single node of LSSC-IV\footnote{\url{http://lsec.cc.ac.cn/chinese/lsec/LSSC-IVintroduction.pdf.}}, which is a high-performance computing
cluster maintained at the State Key Laboratory of Scientific and Engineering Computing, Chinese Academy of Sciences. The operating system of LSSC-IV is Red Hat Enterprise Linux Server 7.3. 

\subsection{Parameters Tuning}

In this subsection, we discuss how to choose the parameters when designing the step sizes.

The initial step size $\eta_0$ is chosen from $\left\lbrace 0.1,0.5, 1, 5, 10, 15\right\rbrace$ for the Epoch-SGD algorithm on all step size schedules.  Generally speaking, for $1/t$-band, we do not know exactly the coefficients $m$ and $M$ for the lower and upper boundaries.  In the experiments, the coefficient $m$ is tuned properly using a similar approach as the initial step size $\eta_0$. 
Instead of finding the coefficient $M$ of the upper bound, we tune the bandwidth $s = M/m \in \bigbl 2,3,4,5\bigbr$ for $1/t$ Fix-period band and $1/t$ Grow-period band. The distance of the adjacent nodes $t_i (i\in \N^{+})$ depends on the budget of the outer loop $N$. In our experiments we set $t_{i+1} - t_i = 30$, $t_1=30$ for $1/t$ Fix-period band and $t_{i+1}=2t_i$, $t_1=30$ for $1/t$ Grow-period. From Figure \ref{stepsize:1}, we can see that 
$1/t$ Fix-period , $1/t$ Grow-period  and $1/t$-stepsize coincide in the first cycle and $1/t$ Fix-period also coincides with $1/t$ Grow-period in the second cycle.

The Grow-Exp step size drops by half and the period of per cycle is doubled. The initial period $T_0$ is chosen from $\bigbl 1, 2, 3,  5, 10, 20 \bigbr$. For $1/t$ Grow-Exp, we tune the up-down ratio $\theta \in \bigbl 1.1,1.2, 1.3,1.4,1.5\bigbr$ and  the length of $T_0$ is the same as Grow-Exp.

\subsection{Regularized Logistic Regression}\label{sec:experiment1}
Firstly, we empirically test the above step sizes on the regularized logistic regression problems, which is strongly convex for regularization parameter $\Lambda > 0$
\begin{equation*}
f(x) = \frac{1}{n}\sum_{i=1}^{n} \ln(1+\exp(-b_i\left\langle a_i, x\right\rangle )) + \frac{\Lambda}{2}\left\|x \right\|^2, 
\end{equation*}
where $\left\lbrace a_i, b_i\right\rbrace_{i=1}^n $ is a training sample set with $a_i\in \R^d$ and $b_i\in \left\lbrace -1, +1\right\rbrace $. We use the two binary classification datasets  w8a ($n=49749,d=300$) and rcv1.binary  ($n=20242, d=47236$) from LIBSVM\footnote{\url{https://www.csie.ntu.edu.tw/ cjlin/libsvmtools/datasets/}},  where the 0.75 partition of
the data is used for training and the remaining is for testing. The common parameters $ \Lambda =10^{-4}$, batch size $b=128$, the outer loop $N=120$ and the inner loop $m^{'}=n/128$. 

We plot the average results of 5 runs on w8a in Figure \ref{fig:2}.  For the $x$-axis we always use the number of epochs calculated. The y-axis are the value of the loss function on training dataset (left) and the accuracy (the percent of correctly classified datasets) on testing dataset (right). For $1/t$-stepsize, the best initial step size is achieved at $\eta_0=5$ and we apply the same initial step size for the other step sizes. Other important parameters are set as: $s =3$, $T_0 =2$ and $\theta = 1.2$.  From Figure \ref{fig:2}, we can see that the exponentially decaying step size (Grow-Exp) performs better than $1/t$-stepsize on both training loss and accuracy. Our proposed $1/t$ Fix-period and $1/t$ Grow-period both achieve good performance than $1/t$-stepsize. In addition, $1/t$ Grow-Exp gets higher accuracy than Grow-Exp.  
 \begin{figure}[h]
\centering	
\includegraphics[width=0.49\textwidth,height=2.2in]{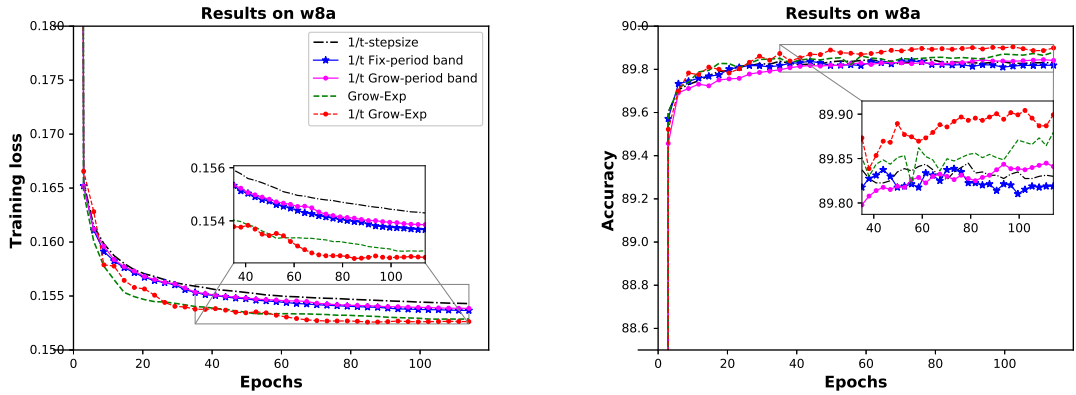}
\hfill 
\includegraphics[width=0.49\textwidth,height=2.2in]{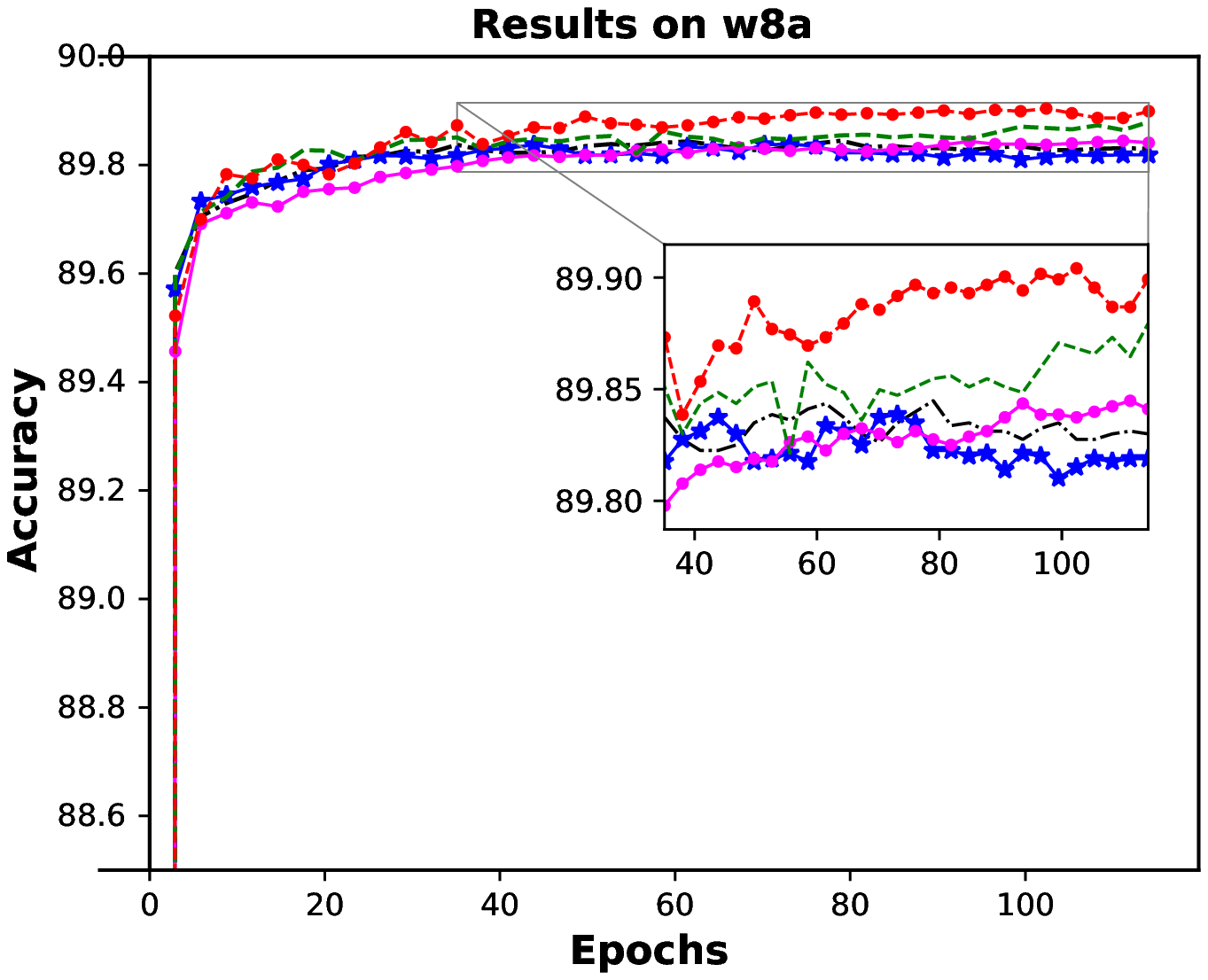}
\caption{Results for regularized logistic regression}
\label{fig:2}  
 \end{figure}

In Figure \ref{fig: rcv1}, we report the average results of 5 runs on rcv1.binary. The best-tuned initial step size $\eta_0$ is 10 for $1/t$-stepsize and we use the same initial step size for other step size schedules. The value of $\theta$ is $1.3$ for $1/t$ Grow-Exp and other parameters are the same as those in w8a.  The similar performance to Figure \ref{fig:2} can be achieved.
 \begin{figure}[h]
\centering	
\includegraphics[width=0.49\textwidth,height=2.2in]{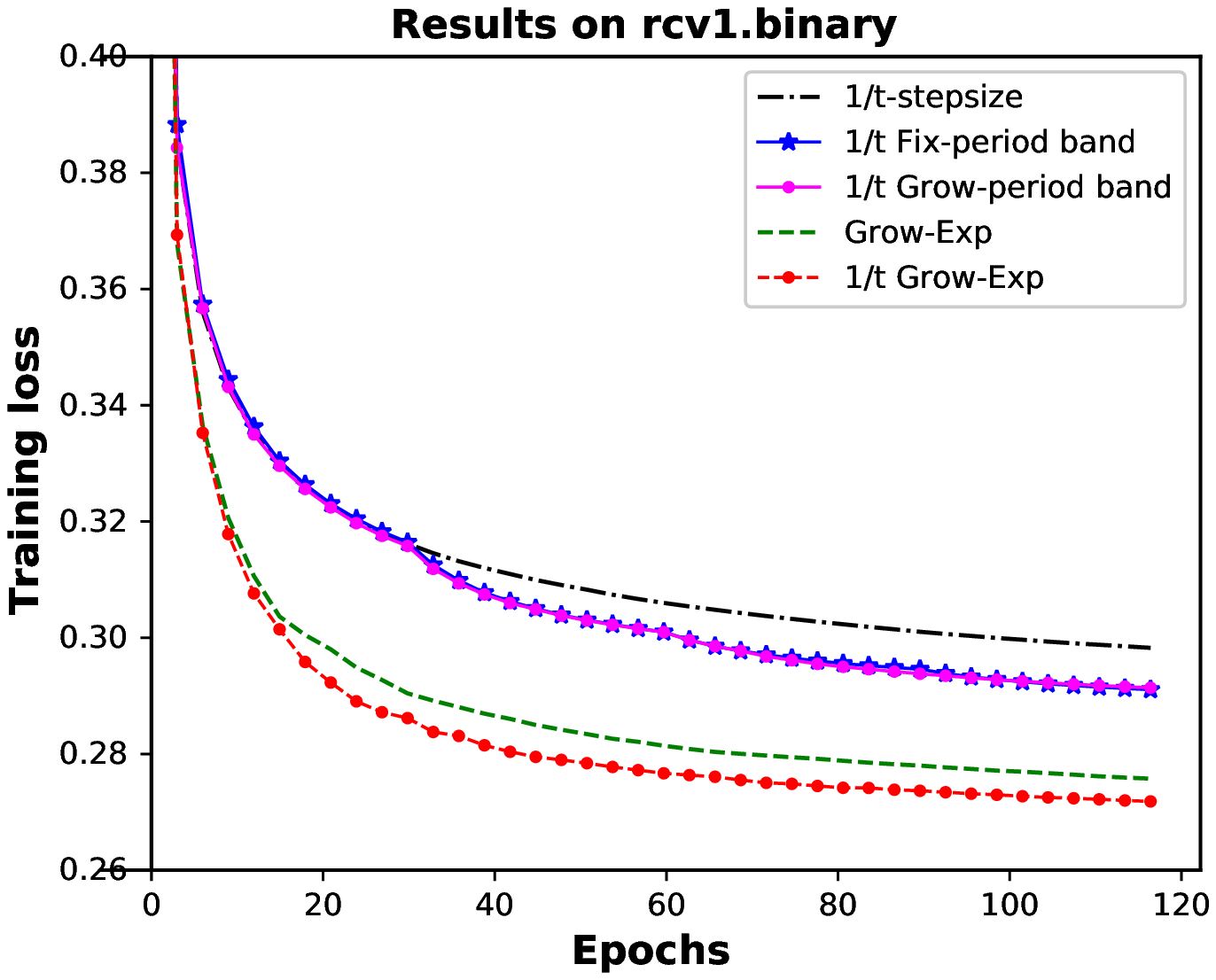}
\hfill 
\includegraphics[width=0.49\textwidth,height=2.2in]{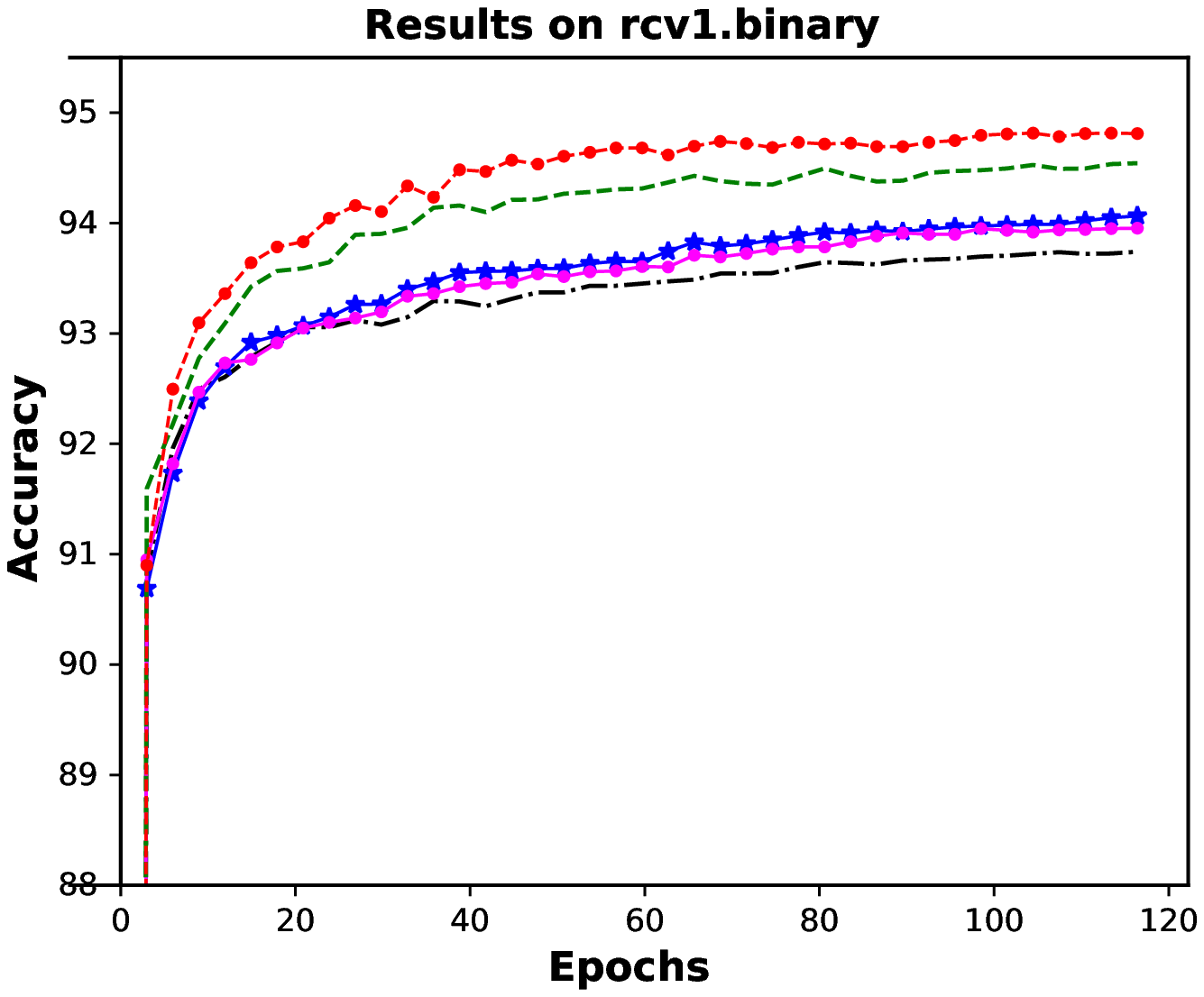}
\caption{Results for regularized logistic regression}
\label{fig: rcv1}  
 \end{figure}
 
From Figures \ref{fig:2} and \ref{fig: rcv1}, the four $1/t$-band  step sizes including $1/t$ Fix-period band, $1/t$ Grow-period band, Grow-Exp and $1/t$ Grow-Exp all perform better than $1/t$-stepsize in strongly convex setting. In particular, the type of the Grow-Exp step size significantly improves the performance of Epoch-SGD over $1/t$-stepsize. This also implies that the relatively large step size at the initial iterations possibly makes the algorithm drop rapidly. It is also observed that the proposed $1/t$ Grow-Exp step size, based on the $1/t$ up-down policy, yields better performance compared to the Grow-Exp step size.

\subsection{Deep Neural Network and Residual Neural Network}\label{sec:experiment2}
In this subsection, we carry out the experiments on some standard datasets, e.g., MNIST and CIFAR-100. 

First of all, we test on a fully-connected 3-layer (784-500-300-10) neural network to train {\bf MNIST}\footnote{\url{http://deeplearning.net/data/mnist/}}, consisting of a training set of 60000 images with 28x28 pixels and a testing set of 10000 images in 10 classes. The batch size $b=128$, the outer loop $N=120$ and the inner loop $m^{'}=n/128$. For the $1/t$-stepsize, the best $\eta_0$ is achieved at $\eta_0 = 0.5$ based on its accuracy.  For $1/t$ Fix-period band and $1/t$ Grow-period band, $\eta_0 $ is the same as that of the $1/t$-stepsize.  We choose $s =3$, that is $\eta(t) \in [\eta_0/t, 3\eta_0/t]$.  For Grow-Exp, the parameters are set as $\eta_0 = 0.5, T_0 = 10$. For $1/t$ Grow-Exp, we set $\theta = 1.3$ and other parameters are the same as Grow-Exp. The average results of 5 runs are given in Figure \ref{fig:4}.  It is easy to see that the Grow-Exp type step size achieves better performance compared to $1/t$-stepsize, $1/t$ Fix-period band and $1/t$ Grow-period band. Besides, our proposed $1/t$ Grow-Exp achieves lower training loss than Grow-Exp.
 \begin{figure}[h]
\centering	
\includegraphics[width=0.49\textwidth,height=2.2in]{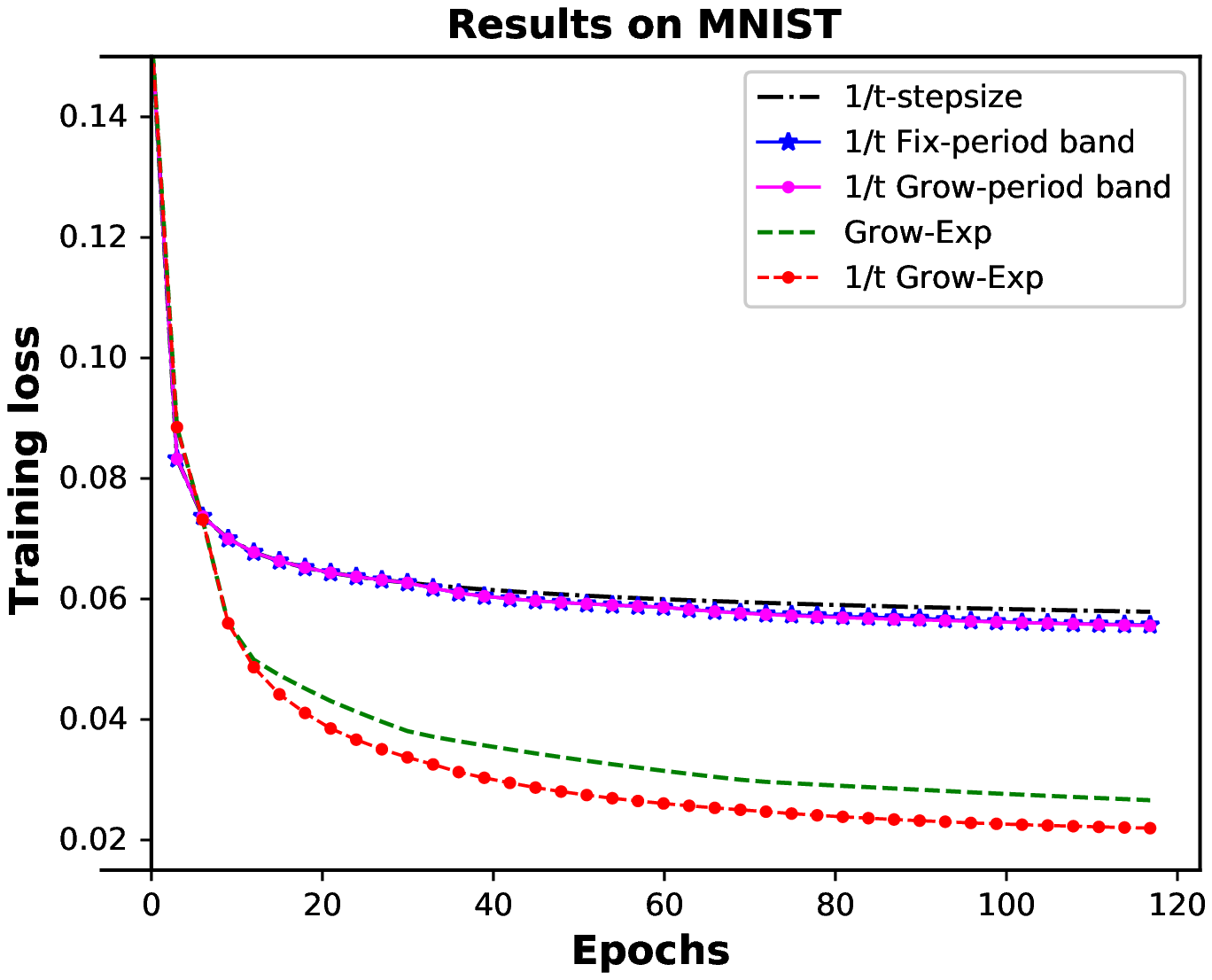}
\hfill 
\includegraphics[width=0.49\textwidth,height=2.2in]{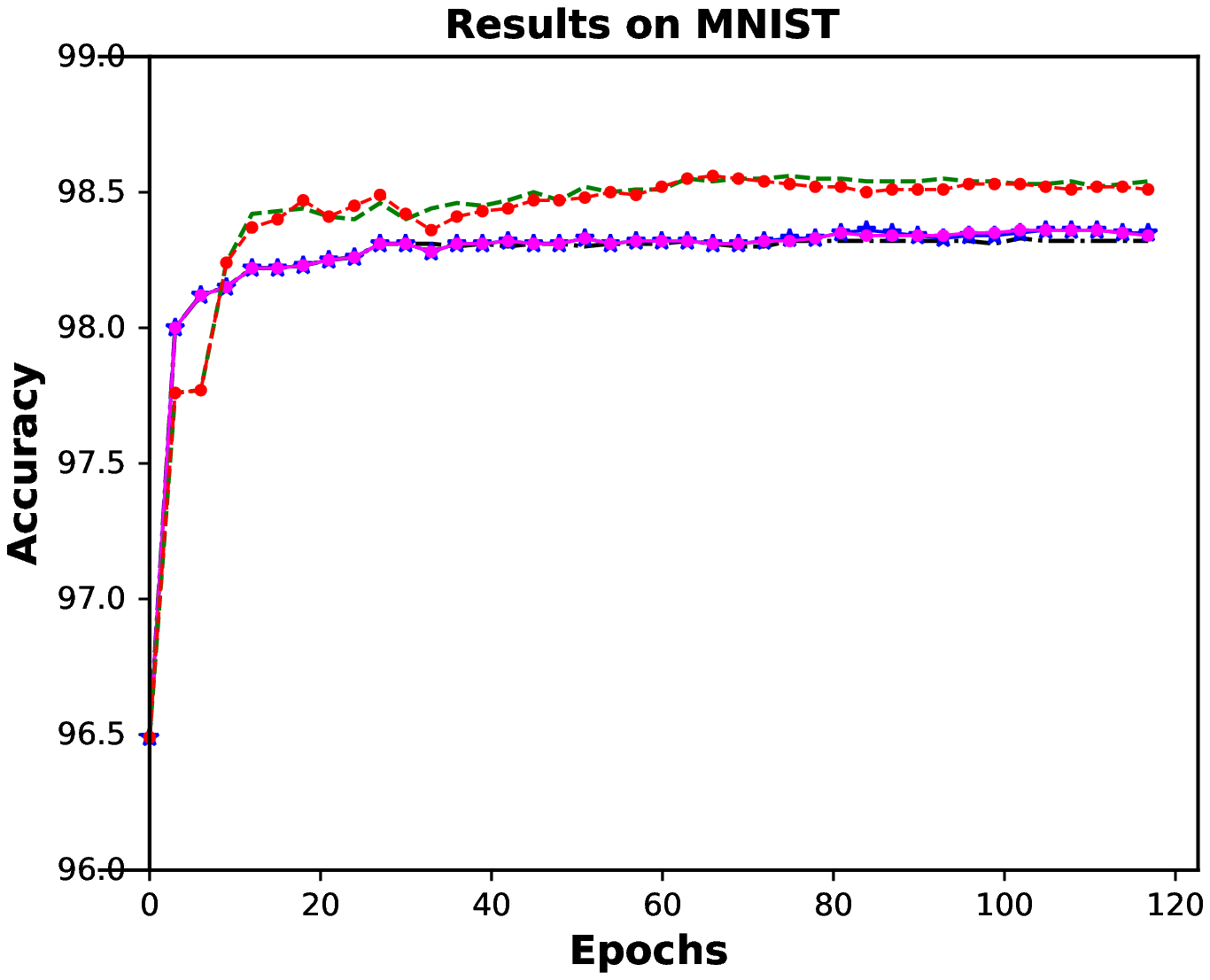}
\caption{Results  on  deep neural networks (DNNs)}
\label{fig:4}  
 \end{figure}

Next, we implement the above five step sizes on {\bf ResNet-18} \citep{he2016deep} with  {\bf CIFAR-100}\footnote{\url{https://www.cs.toronto.edu/~kriz/cifar.html}}. The CIFAR-100 dataset consists of 60000 32x32 color images in 100 classes, and 50000 images for training and remaining 10000 images for testing.
 For $1/t$-stepsize, we set $\eta(t) = \eta_0/ (1+t/10)$, where $\eta_0  \in \left\lbrace 0.1, 0.5, 1, 5, 10, 15\right\rbrace $. The best performance of $1/t$-stepsize is achieved at $\eta_0= 1$. In this case the bandwidth $s=3$. Other important parameters are the same as the experiment in DNNs. For Grow-Exp, $\eta_0 = 0.5$ and $T_0 = 10$. For $1/t$ Grow-Exp, $\eta_0 = 0.5$, $T_0=10$ and $\theta=1.3$. 

We repeat the training process 5 times and the average results (the left is testing loss function, the right is the accuracy on testing dataset) are presented in Figure \ref{fig:5}. In this case, we see that the sudden increase in $1/t$ Fix-period band and $1/t$ Grow-period band may lead to a short-term negative effect but overall outperform $1/t$-stepsize at long term training. Especially, $1/t$ Grow-period band performs better than $1/t$ Fix-period band. The frequently going up and down makes $1/t$ Fix-period band less stable than $1/t$ Grow-period band. This may be the main reason for this phenomenon. 
 \begin{figure}[h]
\centering	
\includegraphics[width=0.49\textwidth,height=2.2in]{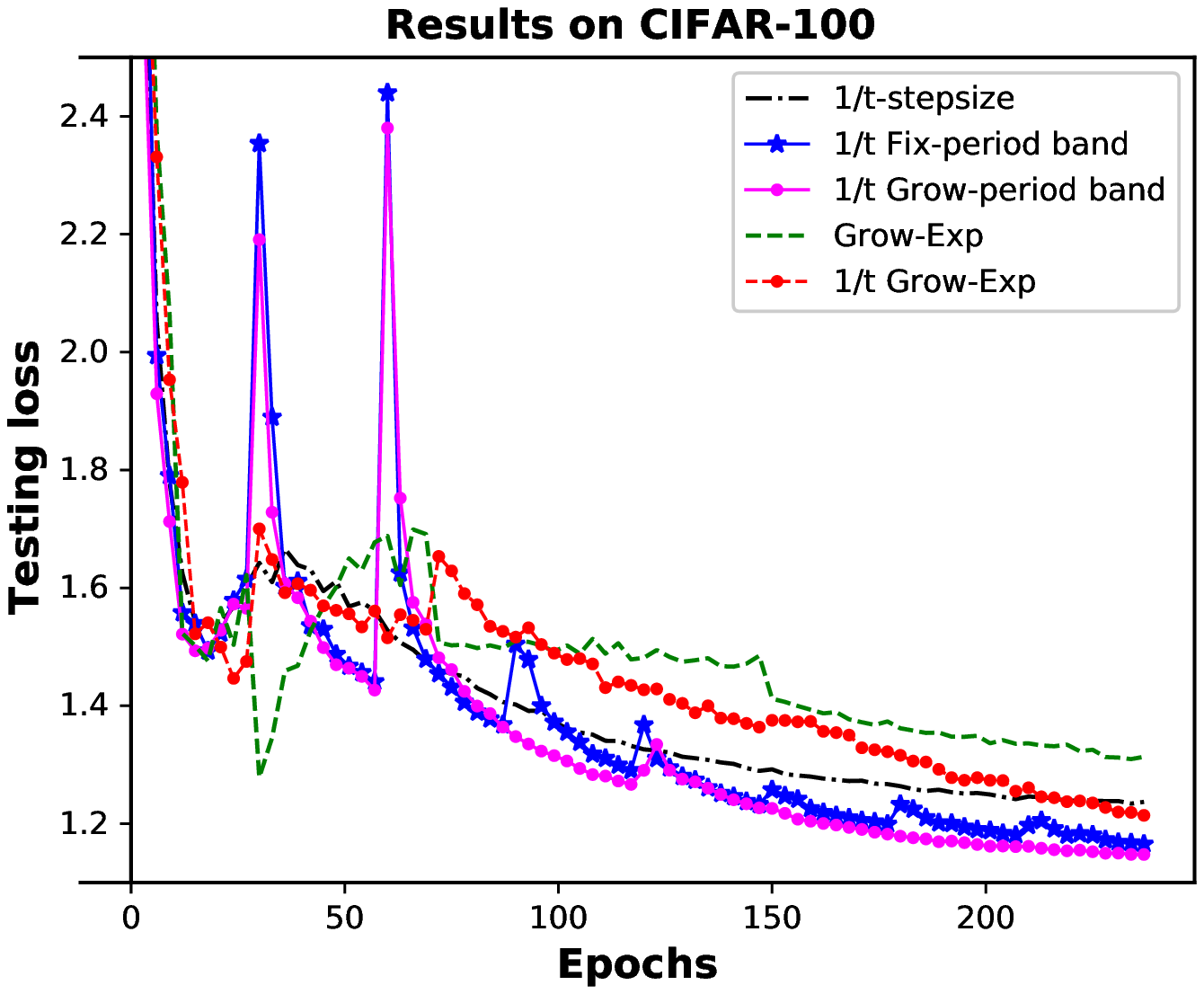}
\hfill 
\includegraphics[width=0.49\textwidth,height=2.2in]{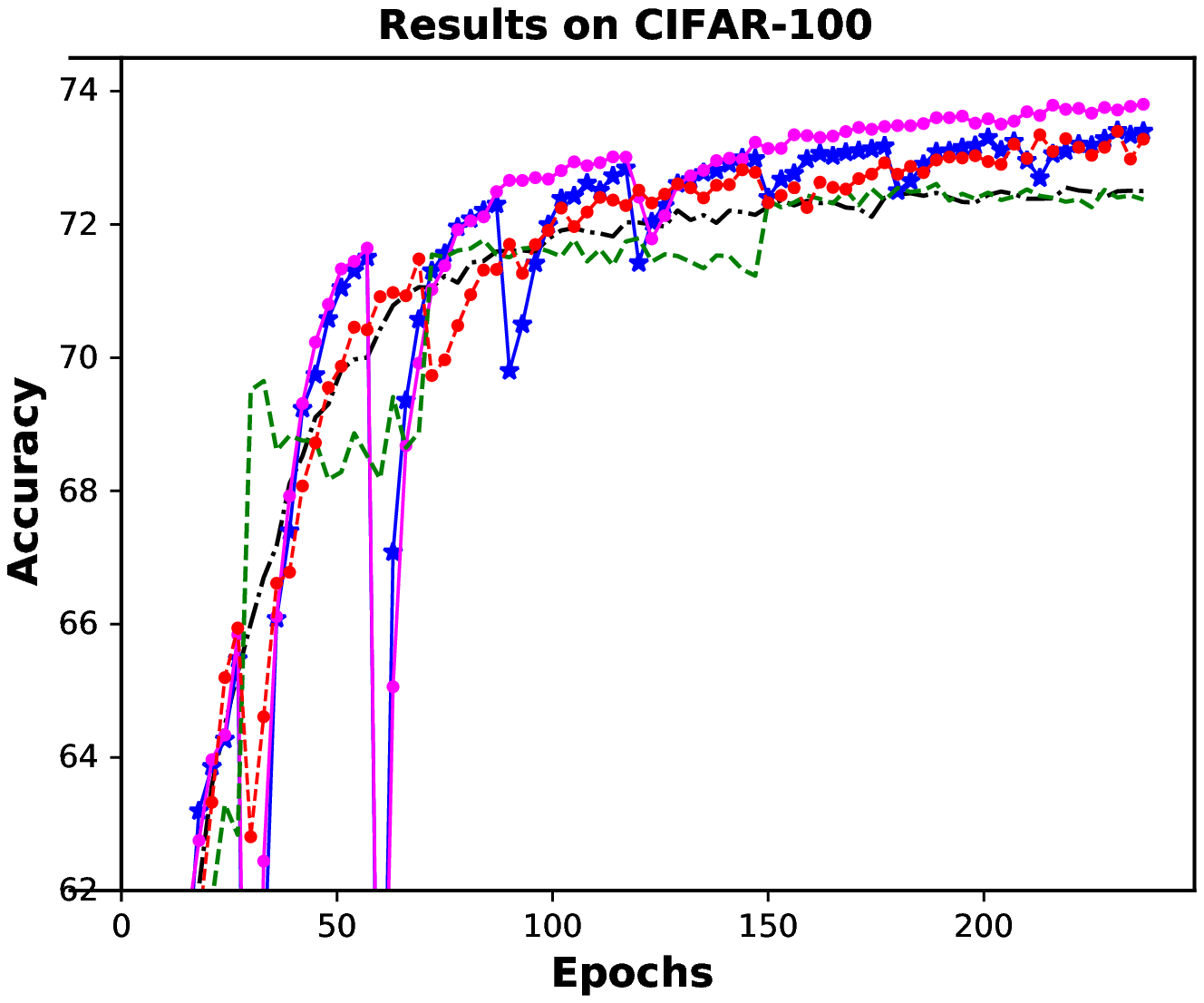}
\caption{Results  on  ResNet-18}
\label{fig:5}  
 \end{figure}

Another observation from Figure \ref{fig:5} is that the Grow-Exp type step size does not work well as Section \ref{sec:experiment1}.  This may be because that a growing number of epochs in Grow-Exp might reduce its generalization at the final stage per cycle.  Nevertheless, we can find that $1/t$ Grow-Exp yields better performance than Grow-Exp. Indeed, the $1/t$-stepsize scheme may not be the best baseline for solving nonconvex problems.  We take it as an example here and empirically demonstrate that the step size based on bandwidth is potential and often helps in practice.


 
\subsection{Additional Experiments on Other Algorithms and Step Sizes}

For further investigation, more experiments are carried out to compare different step sizes for Epoch-SGD and other default algorithms in deep learning including  SGD with momentum (called Momentum for short), averaged SGD (called ASGD) \citep{polyak1992} and Adam \citep{Adam}.  We use two popular datasets: CIFAR-10\footnote{\url{http://www.cs.toronto.edu/~kriz/cifar.html}}  and CIFAR-100 for image classifications. The CNN architectures VGG-16 \citep{Simonyan2015}  and ResNet-18 \citep{he2016deep} are adopted for training CIFAR-10 and CIFAR-100, respectively. 

In addition to the step sizes tested in the above subsections,  we implement the popular exponentially decaying step size with a fixed period $T_0$ (called {\bf Fix-Exp}) which has been discussed in Section \ref{sec:3:1}: 
\begin{equation}\label{exp_decay_const}
\eta(t) = \eta_i = \eta_0/10^i, t \in [T_i, T_{i+1}), \,\, T_{i+1} - T_i = T_0, i \in \N.
\end{equation}
Let $\eta_{\min}^i = \eta_i $ for $i \in \N^{+}$ and we define $\eta_{\max}^{i} = \theta \eta_{\min}^{i-1}$ where $\theta \in (1, 1.5]$. Based on (\ref{exp_decay_const}), we propose the following step size (called $\bm{1/t}$ {\bf Fix-Exp}):
\begin{equation}
\eta(t) = \frac{\hat{A}_i}{\hat{B}_i t + 1} \in [\eta_{\min}^i, \eta_{\max}^{i}], \, \, t \in [T_i, T_{i+1}) \,,\, T_{i+1} - T_i = T_0.
\end{equation}
This is similar to $1/t$ Grow-Exp but the number of epochs per cycle is the same. Besides, we also implement the two cyclical step sizes: triangular policy \citep{smith2017cyclical} and cosine annealing \citep{loshchilov2016sgdr}.

Firstly, we test on  VGG-16  for training CIFAR-10. The baseline initial step size is set as $\eta_0=1$ for SGD and ASGD, $\eta_0=0.1$ for Momentum and $\eta_0=0.001$ for Adam. For Momentum, $\beta = 0.9$. In Adam, $(\beta_1,\beta_2)=(0.9,0.99)$ is used.  The best-tuned  value of weight decay is $10^{-4}$ for SGD and ASGD, $5\times 10^{-4}$ for Momentum and $10^{-5}$ for Adam.
The common parameters $N = 120$ and $b=128$ for all algorithms.
We perform the above algorithms with Fix-Exp ($T_0=30$) and $1/t$ Fix-Exp ($T_0=30, \theta = 1.3$). The average results of five runs are presented in Figure \ref{fig:6}.  We can find that $1/t$ Fix-Exp overall shows better performance than Fix-Exp on SGD, Momentum and ASGD, respectively. 
However, the results of  Adam based on Fix-Exp and $1/t$ Fix-Exp almost coincide which implies that the up-down policy may not work well for Adam. 
\begin{figure}[h]
\centering	
\includegraphics[width=0.49\textwidth,height=2.2in]{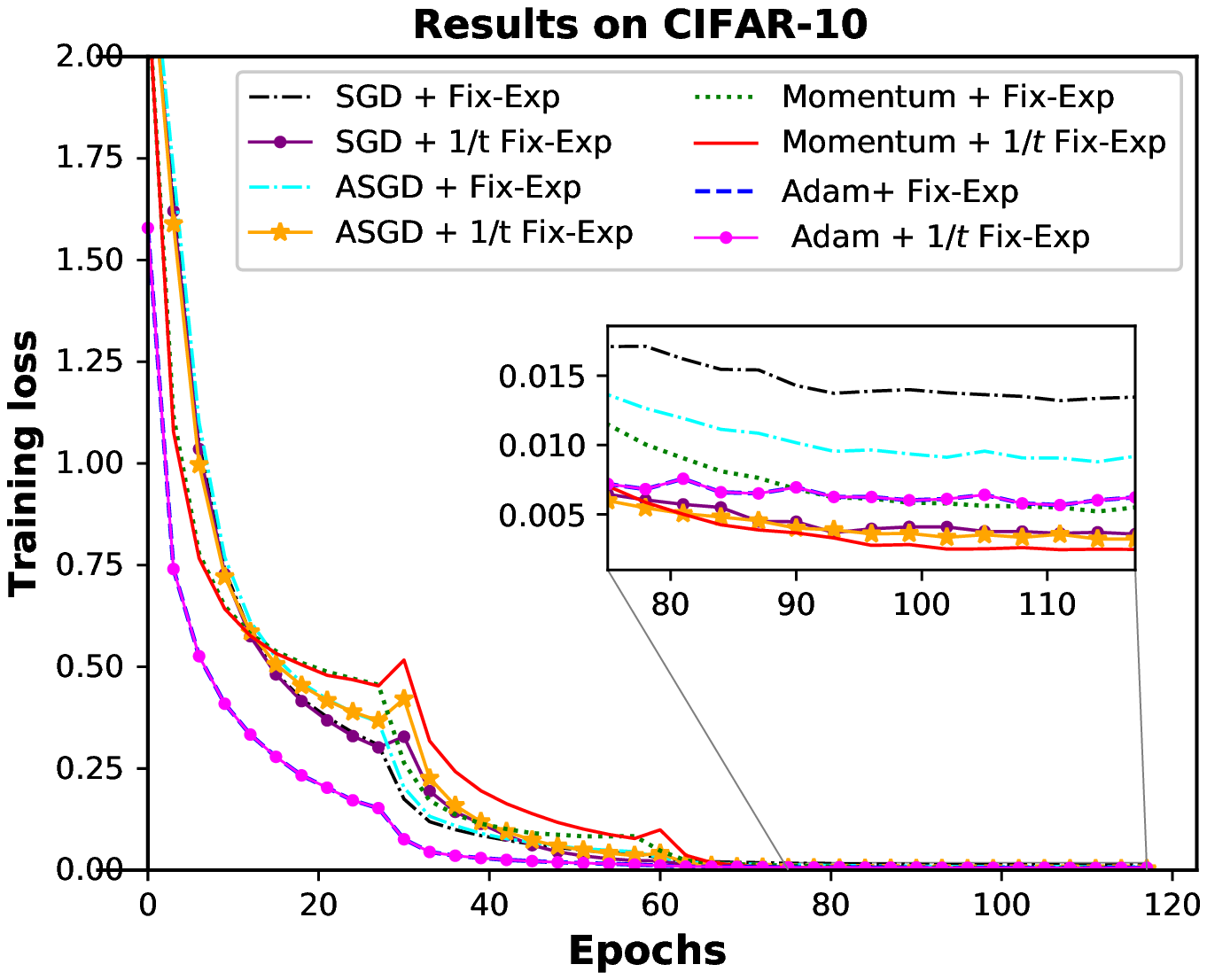}
\hfill
\includegraphics[width=0.49\textwidth,height=2.2in]{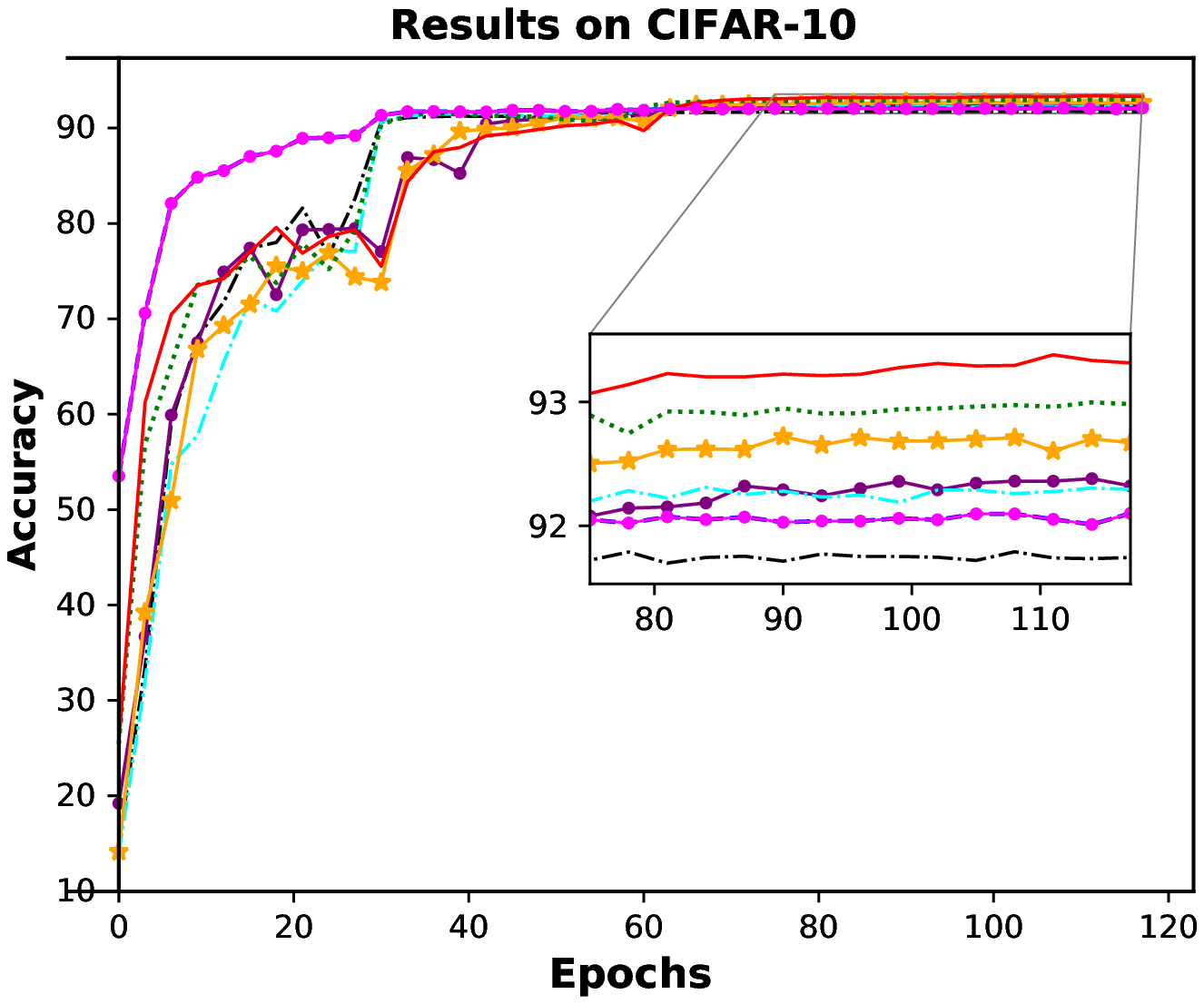}
\caption{Results  on  VGG-16 for CIFAR-10}
\label{fig:6}  
 \end{figure}
 
Besides, we test Momentum with the following step sizes: ({\bf 1}) $1/t$-stepsize ($\eta(t) = \eta_0/(1+t/5)$); ({\bf 2}) $1/t$ Fix-period band ($t_{i+1}-t_i = 30, s=3$); ({\bf 3}) Fix-Exp ($T_0 = 30$); ({\bf 4}) $1/t$ Fix-Exp ($T_0=30,\theta = 1.3$); ({\bf 5}) triangular policy based on (\ref{exp_decay_const}), called ``Triangular" (rise and fall ratio is 1.5); ({\bf 6}) cosine annealing, called ``Cosine" (we use the last iterations as the initial point of restart cycle). All the step sizes are best-tuned with $\eta_0 = 0.1$ and the period of each cycle is $30$ for triangular policy and cosine annealing. The average results of 5 runs are shown in Figure \ref{fig:7}. We observe that $1/t$ Fix-Exp shows its advantages over $1/t$-stepsize, $1/t$ Fix-period band, Fix-Exp and triangular policy after 80 epochs and the final results are comparable to cosine annealing.

 \begin{figure}[h]
\centering	
\includegraphics[width=0.49\textwidth,height=2.2in]{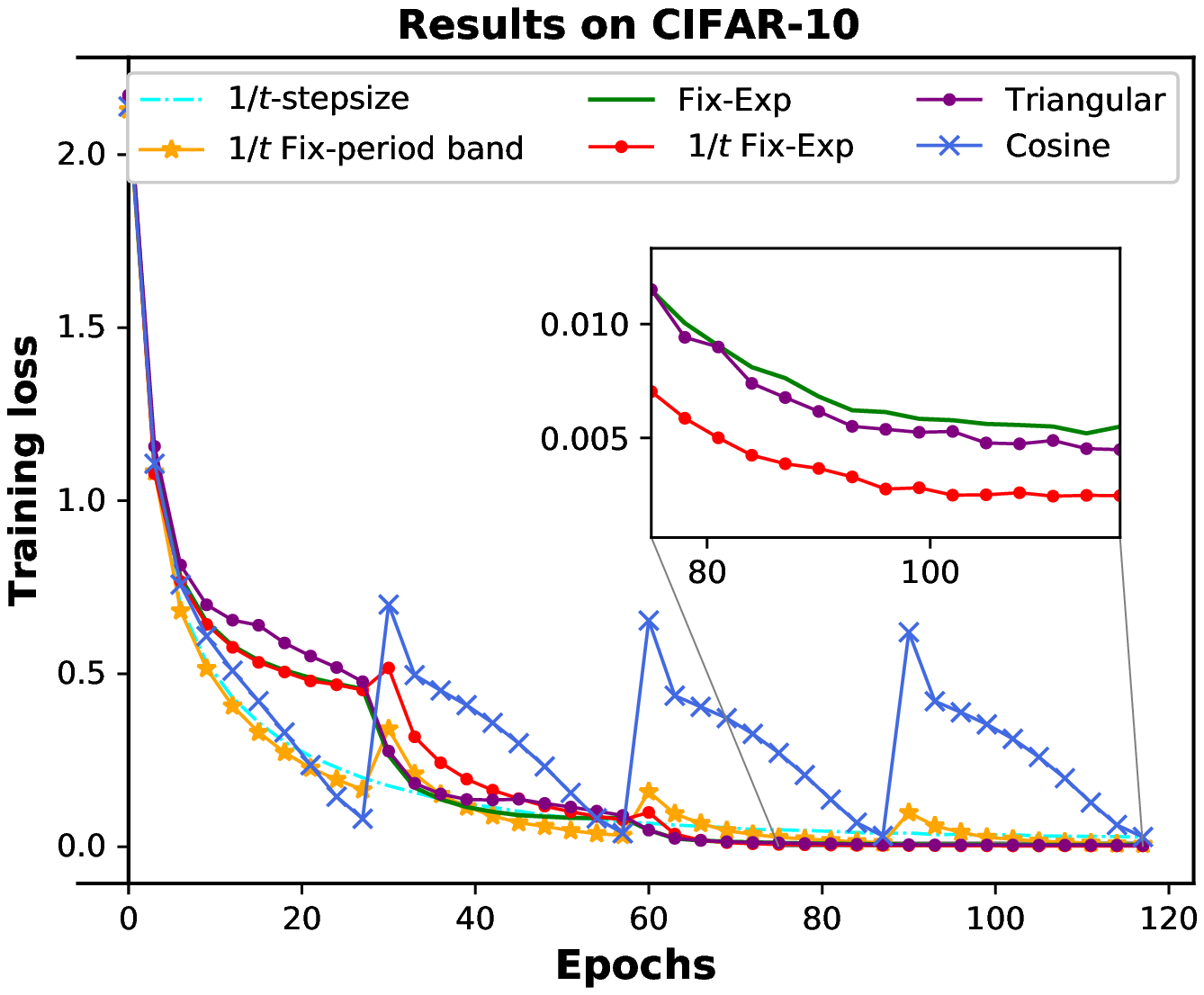}
\hfill
\includegraphics[width=0.49\textwidth,height=2.2in]{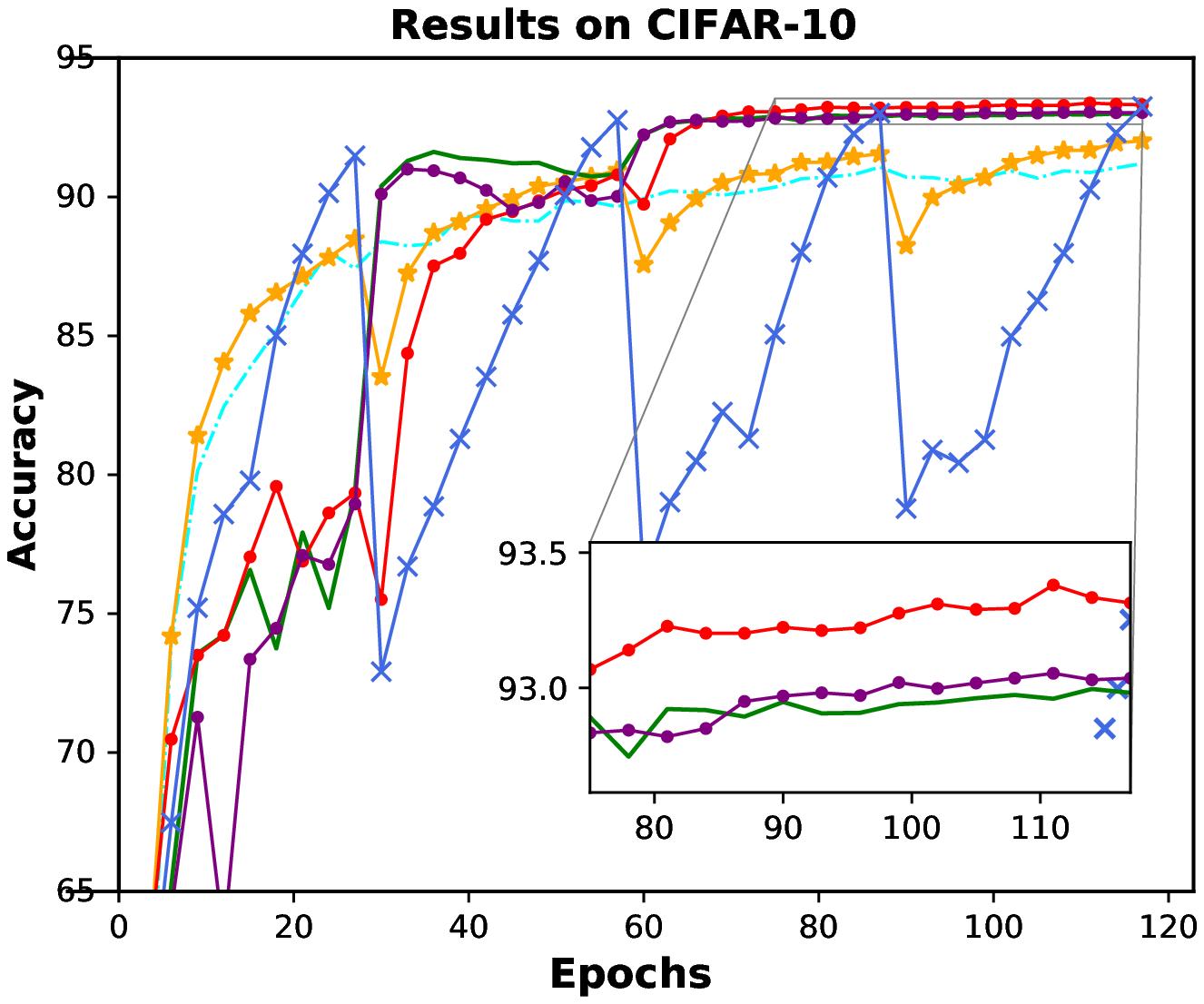}
\caption{Results  of different step sizes for  CIFAR-10}
\label{fig:7}  
 \end{figure}

Next, we implement the above algorithms with Fix-Exp and $1/t$ Fix-Exp on ResNet-18 for training CIFAR-100. The average results of five runs are reported in Figure \ref{fig:8}.  The budget  of the outer iteration $N = 240$ and the period of each cycle $T_0=60$. The other parameters are chosen the same as the experiments on CIFAR-10. Similarly, we can conclude that the up-down policy in $1/t$ Fix-Exp indeed lead to improvements after the second cycles over Fix-Exp on SGD, ASGD and Momentum, respectively.  It is also observed that the up-down policy does not work for Adam but does not make Adam worse.

 \begin{figure}[h]
\centering	
\includegraphics[width=0.49\textwidth,height=2.2in]{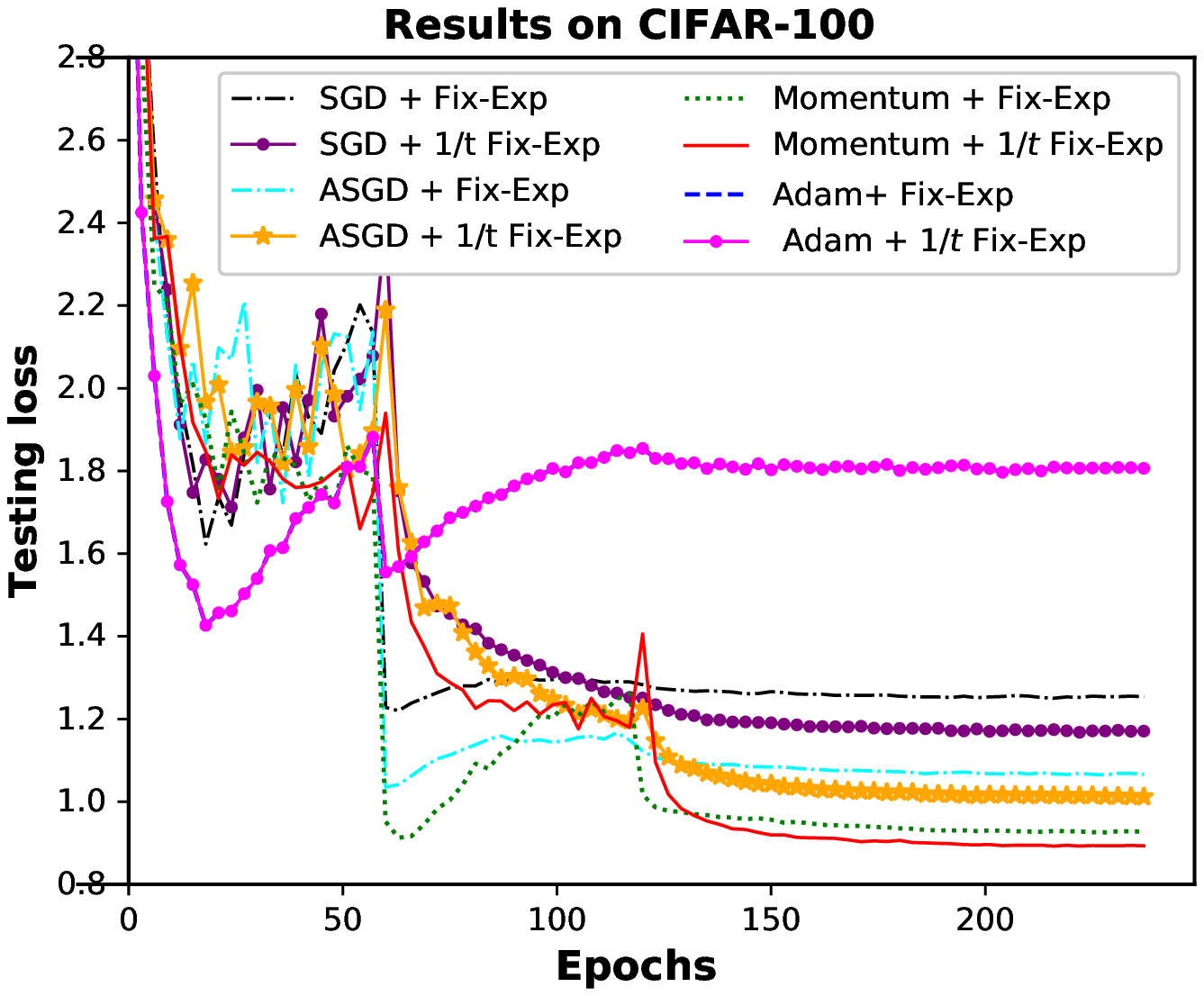}
\hfill
\includegraphics[width=0.49\textwidth,height=2.2in]{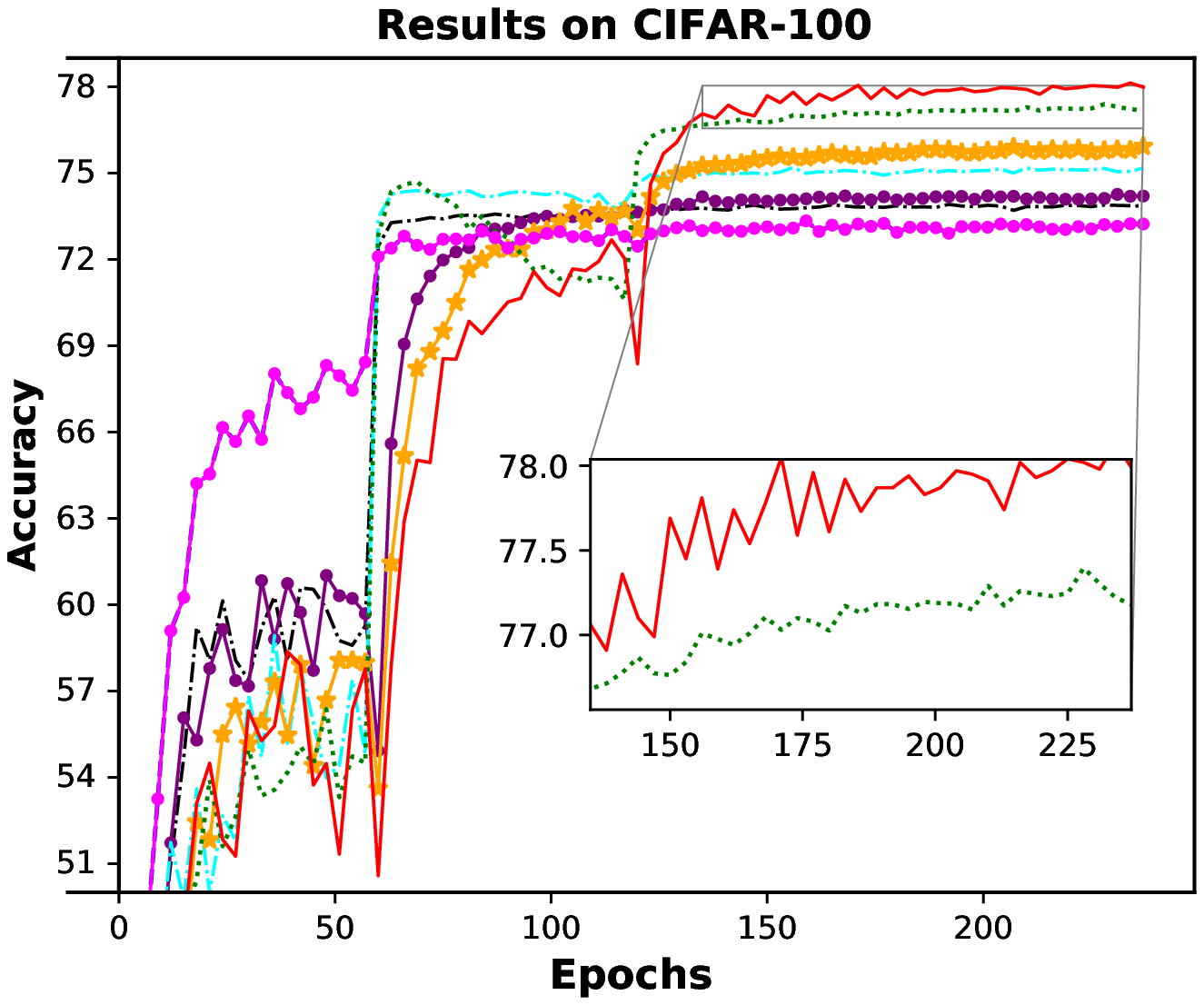}
\caption{Results  on  ResNet-18 for CIFAR-100}
\label{fig:8}  
 \end{figure}
 
In Figure \ref{fig:9}, we report the average results of five runs on the above step sizes for Momentum.  The period for $1/t$ Fix-period band  is $t_{i+1}-t_i = 60$. For Fix-Exp,$1/t$ Fix-Exp, triangular policy (the ratio of rise and fall  is 2) and cosine annealing, the period per cycle $T_0=60$ and other parameters are the same as those of CIFAR-10. As the figures shows, $1/t$ Fix-Exp is enable to reach lower testing loss and higher accuracy than the other step sizes after about 150 epochs.

 \begin{figure}[h]
\centering	
\includegraphics[width=0.49\textwidth,height=2.2in]{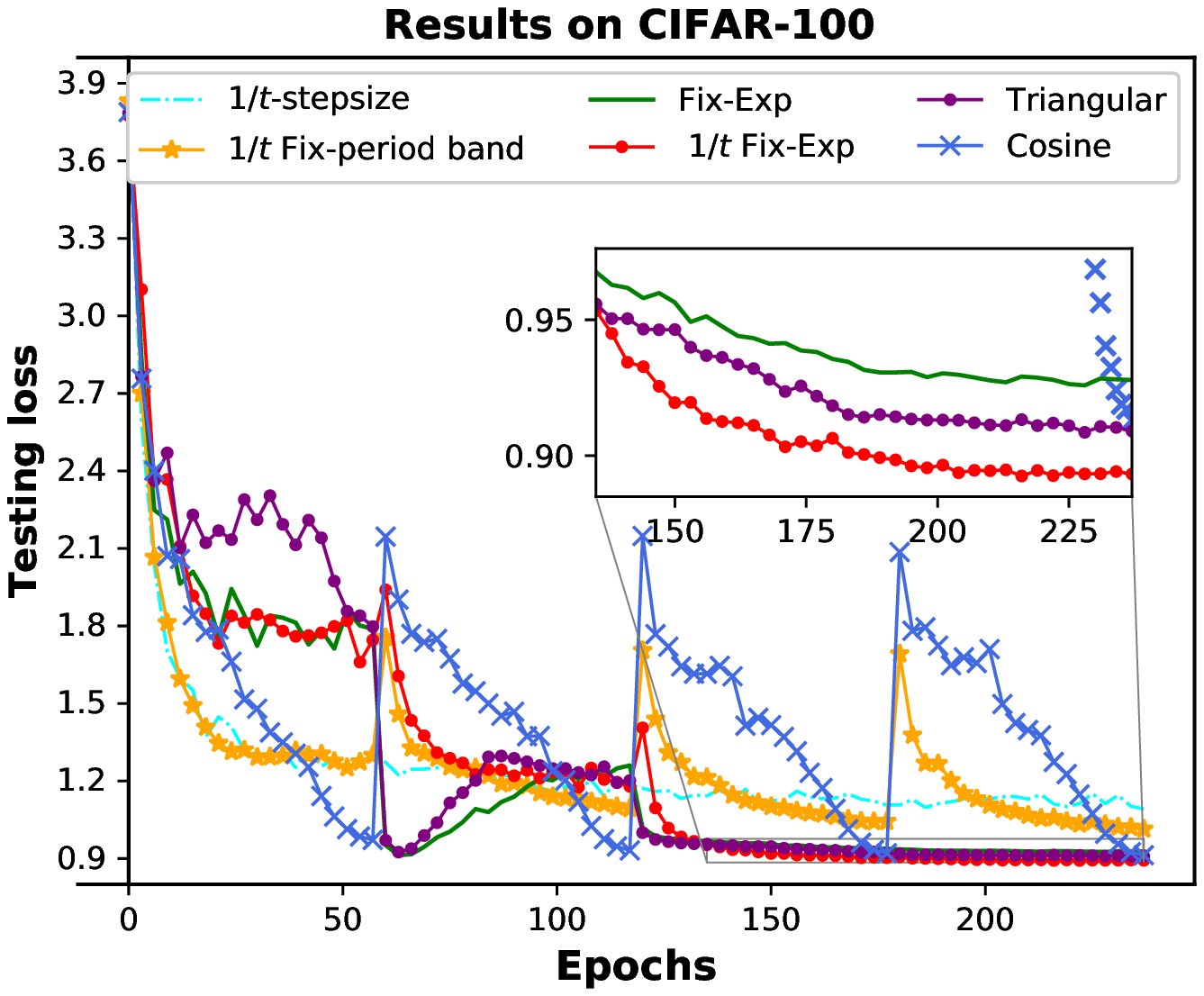}
\hfill
\includegraphics[width=0.49\textwidth,height=2.2in]{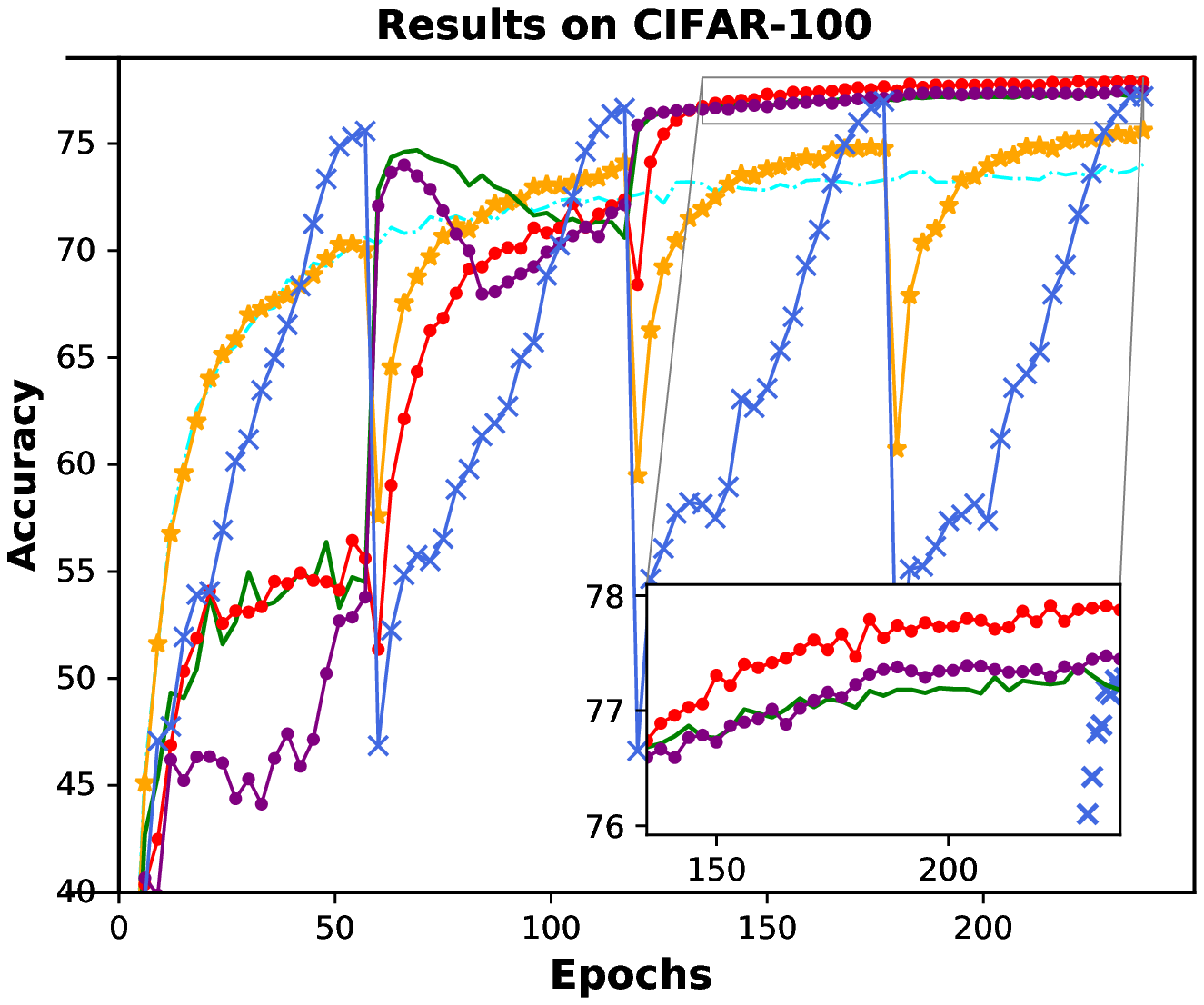}
\caption{Results  of  different step sizes on CIFAR-100}
\label{fig:9}  
 \end{figure}

\section{Conclusion}\label{sec:7}
We propose a bandwidth-based framework for SGD which allows the step size to vary in a banded region and be non-monotonic. We have investigated the conditions where the SGD method achieves an $\mathcal{O}(1/T)$ convergence rate and extended its boundaries at the initial iterations which would be helpful in practical applications. Moreover, we have discussed three different situations which cover most general cases and given explicit error bounds. In many cases, such as $\eta(t) = \eta_0/(t\ln(t))$ and  $\eta_0/\sqrt t$, we have achieved the better upper bounds than those of theorem 10 in \citet{nguyen2018new}. The bandwidth-based step size with different order of the lower and upper bounds often gets worse convergence rates compared to its boundaries. The
convergence rate for some
existing step sizes such as exponentially decaying step size \citep{hazan2014beyond}, cyclical policy \citep{smith2017cyclical} and cosine annealing \citep{sinewave} can be revealed by our analysis if their boundaries 
satisfy the conditions discussed in this paper.

The bandwidth-based framework has more freedom when designing the step size. We have proposed four non-monotonic step sizes based on $1/t$-stepsize and exponentially decaying step size. The numerical results empirically demonstrate their efficiency and potentiality for solving both convex and nonconvex problems, especially for nonconvex problems (e.g., deep neural networks and convolutional neural networks). Besides, we found that the bandwidth-based step size also works for averaged SGD and momentum. 
It is worthwhile to explore SGD and its variants (e.g., momentum) with bandwidth-based step size on nonconvex optimization in the future. We believe that the bandwidth scheme can inspire possibilities for designing more effective step sizes for nonconvex optimization.

The proposed schedule leads to a new prospect based on step size which might be helpful to avoid the saddle points. As we can see, a great of effort has been made to avoid saddle points  by incorporating the noise into search direction per iteration \citep{ge2015escaping,jin2017escape,du2017gradient}.  Whether incorporating the noise or intermediate increasing to step size would help to avoid the saddle points and bad local minimizers will be a very interesting subject for future research.

\acks{The authors would like to thank three anonymous referees and the editor for handling this paper. We are grateful for their comments and suggestions which led to important improvements. The authors thank Dr. Shuxiong Wang for polishing the manuscript.
This research is partially supported by the National Natural Science Foundation of China 11331012 and 11688101. }


\appendix
\section*{Appendix A.} 
\label{app:theorem}

\begin{proof}[{\bf Lemma \ref{sec:lem1}}]
	Due to the $\mu$-strongly convex property of the objective function $f(x)$ for $x \in \R^d$ and $\nabla f(x^{\ast})=0$, let $x=x$ and $\hat{x} = x^{\ast}$ in (\ref{stronglyconvex}), we have 
	\begin{equation}\label{Str_ineq}
	\begin{split}
	f(x) & \geq  f(x^{\ast})  + \left\langle \nabla f(x^{\ast}), x-x^{\ast}\right\rangle  + \frac{\mu}{2}\left\|x -x^{\ast} \right\|^2 \\
	& \geq f(x^{\ast}) +  \frac{\mu}{2}\left\|x -x^{\ast} \right\|^2. \\
	\end{split}
	\end{equation}	
	Besides, letting $x=x^{\ast}$ and $\hat{x} = x$ in (\ref{stronglyconvex}) gives
	\begin{equation*}
	f(x^{\ast}) \geq f(x) + \left\langle \nabla f(x), x^{\ast} - x\right\rangle  + \frac{\mu}{2}\left\|x -x^{\ast} \right\|^2.
	\end{equation*}
	Re-arranging the above inequality,  we have
	\begin{equation}\label{inequ_gradx}
	\left\langle \nabla f(x), x - x^{\ast} \right\rangle \geq f(x) - f(x^{\ast}) + \frac{\mu}{2}\left\|x -x^{\ast} \right\|^2.
	\end{equation}
	Then applying (\ref{Str_ineq}) into (\ref{inequ_gradx}), we obtain that
	\begin{equation}\label{inequ_gradx2}
	\left\langle \nabla f(x), x - x^{\ast}\right\rangle  \geq \mu \left\|x -x^{\ast} \right\|^2.
	\end{equation}
 Introducing a constant $\tau \in [1,2)$ and multiplying (\ref{inequ_gradx}) and (\ref{inequ_gradx2}) by $(2-\tau)$ and $(\tau-1)$, respectively, we have 
	\begin{equation}\label{inequ_gradxf}
	\begin{split}
	\left\langle \nabla f(x), x - x^{\ast}\right\rangle  & = (2-\tau) \left\langle \nabla f(x), x - x^{\ast}\right\rangle  + (\tau-1)\left\langle \nabla f(x), x - x^{\ast}\right\rangle  \\
	& \geq (2-\tau)(f(x) - f(x^{\ast}))  + \frac{\tau\mu }{2}\left\|x -x^{\ast} \right\|^2, 
	\end{split}
	\end{equation}
	as required.	
\end{proof}
\begin{proof}[{\bf Lemma \ref{sec:lem2}}]
		Considering the mini-batch version of the SGD algorithm, we have 
		\begin{equation}{\label{main_ineq}}
		\begin{split}
		\E[ \left\|x_{t+1}-x^{\ast} \right\|^2 | \mathcal{F}_{t}] & =  \E[ \left\|x_{t}- \eta(t) g_t - x^{\ast}\right\|^2 | \mathcal{F}_{t}] \\
		& = \E[\left\|x_{t}- x^{\ast}\right\|^2 | \mathcal{F}_{t}] - \E[2\eta(t) \left\langle g_t, x_{t} - x^{\ast} \right\rangle | \mathcal{F}_{t} ] + \eta(t)^2 \E[\left\|g_t \right\|^2 |  \mathcal{F}_{t}] \\
		& = \left\|x_{t}- x^{\ast}\right\|^2 - 2\eta(t) \left\langle  \nabla f(x_{t}), x_{t}- x^{\ast}\right\rangle + \eta(t)^2 \E[\left\|g_t \right\|^2 |  \mathcal{F}_{t}],
		\end{split}
		\end{equation}
		where the last equality uses the fact that the stochastic gradient $g_t = \frac{1}{b}\sum_{l\in \Omega_t}\nabla f(x_t, \xi_l)$ ($|\Omega_t|=b$) is an unbiased estimation of $\nabla f(x_t)$ at $x_t$.	We assume that Assumption \ref{expected_smooth} holds which induces that there exists a constant $L_f > 0$ such that 
	\begin{equation}\label{sec:lem2:ineq1}
	\begin{split}
	\E[ \left\|g_t \right\|^2 \mid \mathcal{F}_t] & = \E \left[ \left\|\, \frac{1}{b}\sum_{l\in \Omega_t} \nabla f(x_t; \xi_l)\,\right\|^2 \mid \mathcal{F}_t \right] \\ & \leq  \frac{b}{b^2}\sum_{l\in \Omega_t} \E[ \left\|\,\nabla f(x_t; \xi_l)\,\,\right\|^2  \mid \mathcal{F}_t] 
\leq 4L_f(f(x_t) - f^{\ast}) + 2\sigma^2,
	\end{split}
	\end{equation} 
	where $\E[\left\|\nabla f(x^{\ast};\xi) \right\|^2 ] =\sigma^2$.
	Since $f$ is $\mu$-strongly convex, by Lemma \ref{sec:lem1}, the inequality (\ref{inequ_gradxf}) holds.
	 Let $x=x_t$ in (\ref{inequ_gradxf}), together with  (\ref{sec:lem2:ineq1}), then  (\ref{main_ineq}) can be evaluated by
\begin{equation}\label{lem2_inequ1}
\begin{split}
	\E[  \left\|x_{t+1}-x^{\ast} \right\|^2 \mid \mathcal{F}_t] 
  \leq  & \left( 1- \tau\mu\eta(t)\right)\left\|x_t -x^{\ast} \right\|^2\\
 &  + 2\eta(t)^2 \sigma^2 + (4L_f\eta(t)^2 - 2(2-\tau)\eta(t))[f(x_t) - f(x^{\ast})]. 
 \end{split}
\end{equation}
	
Let $n_0 := \sup \left\lbrace t \in \N^{+}: \eta(t) >  \frac{(2-\tau)}{2L_f} \right\rbrace $. For $t > n_0$, we have $4L_f\eta(t)^2 - 2(2-\tau)\eta(t) \leq 0$. Then the inequality (\ref{lem2_inequ1}) can be 
\begin{equation}\label{sec:thm1_ineq1}
\E[ \left\|x_{t+1}-x^{\ast} \right\|^2 \mid \mathcal{F}_t] \leq (1- \tau\mu\eta(t))\left\|x_t -x^{\ast} \right\|^2  + 2\eta(t)^2 \sigma^2.
\end{equation}
Let $\chi_{n_0} = \mathop{\max}\limits_{1 \leq t \leq n_0}\left\lbrace 4L_f\eta(t)^2 - 2(2-\tau)\eta(t)\right\rbrace  $ and $f_{n_0} = \mathop{\max}\limits_{1 \leq t \leq n_0}\left\lbrace f(x_t) - f(x^{\ast})\right\rbrace $. Because $n_0$ is supposed to be a constant which is independent of $T$, the sequence $\left\lbrace f(x_t) - f(x^{\ast})\right\rbrace_{t=1}^{n_0} $ is bounded by a constant $f_{n_0}$. For $ 1 \leq t \leq n_0$, we have
\begin{equation}\label{sec:thm1_ineq}
\E[  \left\|x_{t+1}-x^{\ast} \right\|^2 \mid \mathcal{F}_t]  \leq (1- \tau\mu\eta(t))\left\|x_t -x^{\ast} \right\|^2  + 2\eta(t)^2 \sigma^2 + \chi_{n_0}f_{n_0}.
\end{equation}
For $t > n_0$, taking expectations again and applying the recursion of (\ref{sec:thm1_ineq1}) and (\ref{sec:thm1_ineq}) from $1$ to $t$, we have
\begin{align}\label{main_ineq2}
& \E[  \left\|x_{t+1}-x^{\ast} \right\|^2]  \notag \\ \leq \, &  \prod_{l=1}^{t}(1-\tau\mu\eta(l))\left\|x_1- x^{\ast} \right\|^2 + 2\sigma^2\sum_{l=1}^{t}\eta(l)^2\prod_{u > l}^{t}(1 - \tau\mu\eta(u))  + \chi_{n_0}f_{n_0}\sum_{l=1}^{n_0}\prod_{u > l}^{t}(1 - \tau\mu\eta(u)) \notag \\
 \leq \, &  \exp\left( -\tau\mu\sum_{l=1}^{t}\eta(l)\right) \Delta_{n_0}^0  + 2\sigma^2\sum_{l=1}^{t}\eta(l)^2\exp\left( -\tau\mu\sum_{u>l}^{t}\eta(u)\right),
\end{align}
where $\Delta_{n_0}^0  = \left\|x_1- x^{\ast} \right\|^2 + \frac{n_0 \chi_{n_0}f_{n_0}}{\exp\left( -\tau\mu \sum_{l=1}^{n_0}\eta(l)\right) }$.
The last inequality of (\ref{main_ineq2}) uses the fact that $ 1 + x \leq \exp(x)$ for all $x \in \R$. Note that the coefficient $ 1- \tau\mu\eta(l)$ of $\E[\left\|x_{l} -x^{\ast} \right\|^2 ]$ may be negative for the previous finite terms $1 \leq l  \leq t$, so the recursive process starting from $t=1$ is not appropriate. However, due to that $ \exp(-\tau\mu \eta(l))$ is always positive, we might as well relax the upper bound of $\E[\left\|x_{t+1}-x^{\ast} \right\|^2]$ as (\ref{main_ineq2}).
	
\end{proof}
\section*{Appendix B.}
\begin{proof}[{\bf Theorem \ref{sec:thm1}}]		
	In this case, the sequence of step size $\eta(t)$ satisfies that $$  0 <\frac{m}{t} \leq \eta(t) \leq \frac{M}{t}, \,\text{for}\, 1 \leq t \leq T.$$  	
	It is known that 
\begin{subequations}	
	\begin{equation}\label{inequ_sum1/l}
	\ln(t+1) \leq \sum_{l=1}^{t} \frac{1}{l} \leq \ln(t)+1
	\end{equation}
	and 
	\begin{equation}\label{inequ_sum21/l}
	\int_{u = l}^{t+1} \frac{du}{u} \leq \sum_{u=l}^{t} \frac{1}{u} \leq \int_{u=l-1}^{t}\frac{du}{u}, \,	\text{for any\,}\, l>1.
	\end{equation}
\end{subequations}	
	Then we have 
\begin{subequations}	
	\begin{equation}\label{inequ_sum1}
	\sum_{l=1}^{t}\eta(l) \geq \sum_{l=1}^{t} \frac{m}{l} \geq m\ln(t+1)
	\end{equation} 
	and 
	\begin{equation}\label{inequ_sum2}
	\sum_{u>l}^{t}\eta(u) \geq \sum_{u>l}^{t}\frac{m}{u} = \sum_{u=1}^{t}\frac{m}{u} - \sum_{u=1}^{l}\frac{m}{u} \geq m(\ln(t+1) - \ln(l)-1).
	\end{equation}
\end{subequations}	
	Let $n_0 := \sup \left\lbrace t \in \N^{+}: \eta(t) >  \frac{2-\tau}{2L_f} \right\rbrace $. In this case, when $ t \geq 2ML_f/(2-\tau)$, we have 
	\begin{equation}
	\eta(t) \leq \frac{M}{t} \leq \frac{2-\tau}{2L_f}.
	\end{equation}
	Thus, $n_0 \leq 2ML_f/(2-\tau)$ which is independent of  $T$.

From Lemma \ref{sec:lem2}, we know that for $T >  n_0$, $\E[ \left\|x_{T+1}-x^{\ast} \right\|^2] $ can be estimated as 
	\begin{equation}\label{sec:thm1_ineq2}
	\E[ \left\|x_{T+1}-x^{\ast} \right\|^2] \leq \Gamma_T^2 + \Gamma_{T}^2,
	\end{equation} 
where\begin{equation*}
\Gamma_{T}^1   := \exp\left(-\tau\mu\sum_{l=1}^{T}\eta(l)\right) \Delta_{n_0}^0,  \,\,\,\, 
\Gamma_{T}^2  :=2 \sigma^2\sum_{l=1}^{T}\eta(l)^2\exp\left(-\tau\mu\sum_{u>l}^{T}\eta(u)\right).
\end{equation*} 

	Applying (\ref{inequ_sum1}) into $\Gamma_T^1$, we can achieve that
	\begin{align}
	\Gamma_{T}^1 & \leq \exp\left(-\tau\mu m \ln(T+1)\right) \Delta_{n_0}^0 = \frac{\Delta_{n_0}^0 }{(T+1)^{\tau\mu m}} .
	\end{align}
	
	Now, we proceed to obtain the upper bound for $\Gamma_T^2$. Using the upper bound of $\eta(t)$ and (\ref{inequ_sum2}) gives
	\begin{equation*}
	\begin{split}
	\Gamma_{T}^2 & = 2\sigma^2\sum_{l=1}^{T}\eta(l)^2\exp\left(-\tau\mu\sum_{u>l}^{T}\eta(u)\right) \\
	& \leq  2\sigma^2\sum_{l=1}^{T}\eta(l)^2\exp(-\tau\mu m (\ln(T+1) -\ln(l) - 1))\\
	& \leq  \frac{2\sigma^2 M^2\exp(\tau\mu m)}{(T+1)^{\tau\mu m}} \sum_{l=1}^{T} \frac{1}{l^2}\cdot\exp(\tau\mu m \ln(l))  \leq \frac{2\sigma^2 M^2\exp(\tau\mu m)}{(T+1)^{\tau\mu m}} \sum_{l=1}^{T} \frac{l^{\tau\mu m}}{l^2}. \\
	\end{split}
	\end{equation*}  	
	If $m= \frac{1}{\tau\mu}$, then 
	\begin{equation*}
	\Gamma_{T}^2 \leq 2\sigma^2 M^2\exp(1)\cdot\frac{\ln(T) +1}{T+1}.
	\end{equation*}
	However when $m \neq \frac{1}{\tau\mu}$, whether $ \tau\mu m $ is greater than 2 or the other case, we have
	\begin{equation}\label{sum_upper}
	\sum_{l=1}^{T} \frac{l^{\tau\mu m}}{l^2} = \sum_{l=1}^{T} l^{(\tau\mu m -2)} \leq \int_{l=1}^{T+1} l^{(\tau\mu m -2)}dl +1,
	\end{equation}
	then
	\begin{equation*}
	\Gamma_{T}^2 \leq \frac{2\sigma^2 M^2\exp(\tau\mu m)}{(\tau\mu m -1)} \cdot \frac{(T+1)^{\tau\mu m -1} + \tau\mu m -2 }{(T+1)^{\tau\mu m}}.
	\end{equation*}
	
Substituting the upper bounds of $\Gamma_1^T$ and $\Gamma_2^T$ into (\ref{sec:thm1_ineq2}), we get the desired result.
	
\end{proof}

\begin{proof}[{\bf Theorem \ref{sec:thm2_average}}]	
 Let $n_1 := \sup \left\lbrace t \in \N^{+}: \eta(t) >  \frac{2-\tau}{4L_f} \right\rbrace $. In this case, $ \frac{m}{t} \leq \eta(t) \leq \frac{M}{t}$ which implies that $\delta_1(t) = \delta_2(t) = 1/t$.  When $t \geq 4ML_f/(2-\tau)$, we have $\eta(t) \leq (2-\tau)/(4L_f)$. Thus we know $n_1 \leq 4ML_f/(2-\tau)$  which is independent of $T$.  Let $\chi_{n_1} = \mathop{\max}\limits_{1 \leq t \leq n_1}\left\lbrace 4L_f\eta(t)^2 - 2(2-\tau)\eta(t)\right\rbrace  $ and $f_{n_1} = \mathop{\max}\limits_{1 \leq t \leq n_1}\left\lbrace f(x_t) - f(x^{\ast})\right\rbrace $. Because $n_1$ is a constant, the sequence $\left\lbrace f(x_t) - f(x^{\ast})\right\rbrace_{t=1}^{n_1} $ can be bounded by $f_{n_1}$ which is a constant. For $t > n_1$, $4L\eta(t)^2 - 2(2-\tau)\eta(t) \leq - (2-\tau)\eta(t) $, then the inequality (\ref{lem2_inequ1}) in Lemma \ref{sec:lem2} will be 
	\begin{align}
	&  \E[  \left\|x_{t+1}-x^{\ast} \right\|^2 \mid \mathcal{F}_t ]   \notag\\
	\leq & \, (1- \tau\mu\eta(t))\left\|x_t -x^{\ast} \right\|^2 + 2\eta(t)^2 \sigma^2 + (4L\eta(t)^2 - 2(2-\tau)\eta(t))[f(x_t) - f(x^{\ast})] \notag\\
	\leq  	& \, (1- \tau\mu\eta(t))\left\|x_t -x^{\ast} \right\|^2 + 2\eta(t)^2 \sigma^2  - (2-\tau)\eta(t)[f(x_t) - f(x^{\ast})].\label{sec:lem3_ineq1}
	\end{align}
	Shifting $[f(x_t) - f(x^{\ast})]$ to the left side and $\E[  \left\|x_{t+1}-x^{\ast} \right\|^2 \mid \mathcal{F}_t ] $ to the right side, we obtain
	\begin{equation*}	
	(2-\tau)\eta(t)[f(x_t) - f(x^{\ast})]  \leq (1- \tau\mu \eta(t))\left\|x_t -x^{\ast} \right\|^2 - \E[  \left\|x_{t+1}-x^{\ast} \right\|^2 \mid  \mathcal{F}_t]  + 2\eta(t)^2 \sigma^2.
	\end{equation*}
Applying the lower bound of $\eta(t) $ into the left side and then dividing the above inequality by $(2-\tau)m \delta_1(t)\delta_1(t+t_0)$ ($t_0 \in \N$) gives
	\begin{equation*}
	\begin{split}
	\frac{f(x_t) - f(x^{\ast})}{\delta_1(t + t_0)}
	& \leq  \frac{1}{(2-\tau)m} \left\lbrace \left(\frac{1}{\delta_1(t) \delta_1(t+t_0)} - \frac{\tau\mu m}{\delta_1(t+t_0)}\right) \left\| x_t - x^{\ast}  \right\|^2 - \frac{\E[\left\| x_{t+1} - x^{\ast} \right\|^2]}{\delta_1(t)\delta_1(t+t_0)} \right\rbrace \\  &\,\,\, +  \frac{2\eta(t)^2 \sigma^2}{(2-\tau)m\delta_1(t)\delta_1(t+t_0)}.
	\end{split}
	\end{equation*}
Summing the above inequality for $t$ from $n_1$ to $ T$,  we get that	
\begin{align}
&  \E\left[ f\left( \frac{\sum_{t=1}^{T}\frac{1}{\delta_1(t + t_0)}x_t}{\sum_{t=1}^{T}\frac{1}{\delta_1(t+t_0)}}\right)  - f(x^{\ast})\right]  \notag\\
\leq \,\, 	& \frac{1}{\sum_{t=1}^{T} \frac{1}{\delta_1(t+t_0)}}\left( \sum_{t = 1}^{n_1}\E\left[ \frac{f(x_t) - f(x^{\ast})}{\delta_1(t + t_0)}\right] +  \sum_{t = n_1+1}^{T}\E\left[ \frac{f(x_t) - f(x^{\ast})}{\delta_1(t + t_0)}\right] \right)  \notag\\
\leq \,\,  & \frac{1}{\sum_{t=1}^{T}\frac{(2-\tau)m}{\delta_1(t+t_0)}}\sum_{t=n_1+1}^{T}\left\lbrace \left( \frac{1}{\delta_1(t) \delta_1(t+t_0)} - \frac{\tau\mu m}{\delta_1(t+t_0)}\right)  \E[\left\| x_t - x^{\ast}  \right\|^2] - \frac{\E[\left\| x_{t+1} - x^{\ast} \right\|^2]}{\delta_1(t)\delta_1(t+t_0)} \right\rbrace  \notag \\
&  + \frac{1}{\sum_{t=1}^{T}\frac{1}{\delta_1(t+t_0)}}\sum_{t=1}^{n_1}\frac{f_{n_1}}{\delta_1(t + t_0)}  + \frac{1}{\sum_{t=1}^{T}\frac{(2-\tau)m}{\delta_1(t+t_0)}}\sum_{t=n_1+1}^{T}\frac{2\eta(t)^2 \sigma^2}{\delta_1(t)\delta_1(t+t_0)}, \label{func_inequ}
\end{align}
where the first inequality follows from the well-known \emph{Jensen inequality} if  $f$ is convex.
If $\tau\mu m$ satisfies the following condition:
\begin{equation}\label{taulambda}
\tau\mu m \geq \frac{1}{\delta_1(t+1)} - \frac{\delta_1(t+t_0 +1)}{\delta_1(t)\delta_1(t+t_0)}   \,(\forall \,t >  n_1),
\end{equation}
by simple calculations, we can show that the coefficient of $\E[\left\|x_t -x^{\ast} \right\|^2 ]$ ($t > n_1$) is non-positive. Taking the form $\delta_1(t) = 1/t$, if $\tau\mu m \geq 1$, the condition (\ref{taulambda}) will hold. Then let $\hat{x}_T = \frac{\sum_{t=1}^{T}(t + t_0)x_t}{S_1}$ and $S_1 = \sum_{t=1}^{T}(t+t_0)$, applying the inequality (\ref{func_inequ}), we get
\begin{align}
\E\left[ f\left( \hat{x}_T \right)  - f(x^{\ast})\right] & \leq \frac{(n_1+t_0+1)}{(2-\tau)mS_1}\left(n_1+1 - \tau\mu m\right)  \E[ \left\| x_{n_1+1} - x^{\ast}  \right\|^2]   + \frac{(1+t_0)(n_1+t_0)f_{n_1}}{2S_1}  \notag \\
&\,\,\,+ \frac{2 \sigma^2M^2}{(2-\tau)mS_1} \sum_{t=n_1+1}^{T} \frac{t(t+t_0)}{t^2}.\label{sec:thm2:core}
\end{align}
By Lemma \ref{sec:lem2}, for $ 1 \leq  t \leq n_1$,  we have that 
\begin{equation}\label{sec:thm2_ineq}
\E[  \left\|x_{t+1}-x^{\ast} \right\|^2 \mid \mathcal{F}_t]  \leq (1- \tau\mu\eta(t))\left\|x_t -x^{\ast} \right\|^2 + 2\eta(t)^2 \sigma^2 + \chi_{n_1}f_{n_1}.
\end{equation}
Applying the recursion of (\ref{sec:thm2_ineq}) for $t$ from $1$ to $n_1$ and taking expectation again gives
\begin{align*}
& \E[  \left\|x_{n_1+1}-x^{\ast} \right\|^2] \notag \\
\leq \,\, &  \exp\left(-\tau\mu\sum_{t=1}^{n_1} \eta(t)\right) \left\|x_1 -x^{\ast} \right\|^2 + 2\sigma^2\sum_{l=1}^{n_1}\eta(l)^2\exp\left( -\tau\mu\sum_{u > l}^{n_1}\eta(u)\right)   \notag \\
&  + \chi_{n_1}f_{n_1}\sum_{l=1}^{n_1}\exp\left(  - \tau\mu \sum_{u > l}^{n_1}\eta(u)\right)   \notag \\
\leq \,\, &  \exp\left(-\tau\mu m \ln(n_1+1)\right)\left\|x_1 -x^{\ast} \right\|^2 + 2\sigma^2M^2\sum_{l=1}^{n_1}\frac{1}{l^2} + n_1\chi_{n_1}f_{n_1}  \notag \\
\leq \,\, & \frac{\left\|x_1 -x^{\ast} \right\|^2}{(n_1+1)^{\tau\mu m}} + 4\sigma^2M^2 + n_1\chi_{n_1}f_{n_1}.
\end{align*}                                                             
 Incorporating the above bound of $\E[  \left\|x_{n_1+1}-x^{\ast} \right\|^2] $   into (\ref{sec:thm2:core}), we can obtain that 
\begin{align*}
\E\left[ f\left( \hat{x}_T \right)  - f(x^{\ast})\right] & \leq \frac{(n_1+t_0+1)\left(n_1+1 - \tau\mu m\right)}{(2-\tau)m S_1}\left[ \frac{\left\|x_1 -x^{\ast} \right\|^2}{(n_1+1)^{\tau\mu m}}  + 4\sigma^2M^2 + n_1\chi_{n_1}f_{n_1}\right]  \notag \\
& \,\,\, + \frac{(1+t_0)(n_1+t_0)f_{n_1}}{2S_1}  + \frac{2 \sigma^2M^2}{(2-\tau)mS_1}(T-n_1+t_0\ln(T/n_1)) \\
&  = \frac{1}{(2-\tau)mS_1}\left[\upsilon_1\Delta_{n_1}^0 + \upsilon_2(1-\frac{\tau}{2})mf_{n_1}+2 \sigma^2M^2(T-n_1+t_0\ln(T/n_1))\right],
\end{align*}
where $\hat{x}_T = \frac{\sum_{t=1}^{T}(t + t_0)x_t}{S_1}$, $S_1 = \frac{T(T+t_0)(t_0+1)}{2}$, $\Delta_{n_1}^0 = \frac{\left\|x_1 -x^{\ast} \right\|^2}{(n_1+1)^{\tau\mu m} } + 4\sigma^2M^2 + n_1\chi_{n_1}f_{n_1}$,   $\upsilon_1=(n_1+t_0+1)\left(n_1+1 - \tau\mu m\right)$  and  $\upsilon_2=(1+t_0)(n_1+t_0)$.  
                                          
\end{proof}

\begin{proof}[{\bf Theorem \ref{sec:thm3}}]	
In this case, we assume that $\eta(t)$ satisfies conditions $(A_1)$ and $(B)$.	Similar to Theorem \ref{sec:thm1}, let  $n_0 := \sup \left\lbrace t \in \N^{+}: \eta(t) >  \frac{2-\tau}{2L_f} \right\rbrace $. We know $n_0  \leq 2ML_f/(2-\tau)$ which is independent of  $T$. Then for $T >  n_0$, the conclusion of Lemma \ref{sec:lem2} is true.

 Let $t^{\ast} = 1$ in $(A_1)$, we have 
	\begin{equation*}
	\sum_{t=1}^{T}\eta(t) \geq C\ln(T+1),
	\end{equation*}
then $\Gamma_T^1$ defined by (\ref{gamma1}) can be evaluated as follows
	\begin{equation}\label{Gamma_C1}
	\Gamma_{T}^1  = \exp\left( -\tau\mu\sum_{l=1}^{T}\eta(l)\right)\Delta_{n_0}^0  \leq \frac{1}{(T+1)^{(\tau\mu C)}} \Delta_{n_0}^0 .
	\end{equation}
	Recalling the definition of $\Gamma_T^2$ in (\ref{gamma2}), we have
	\begin{equation*}
	\begin{split}
	\Gamma_{T}^2 & = 2\sigma^2 \sum_{t=1}^{T}\eta(t)^2\exp\left( -\tau\mu\sum_{u>t}^{T}\eta(u)\right)  \leq 2\sigma^2 M^2 \sum_{t=1}^{T}\frac{1}{t^2}\cdot\exp\left( -\tau\mu\sum_{u>t}^{T}\eta(u)\right)  \\
	& \leq 2\sigma^2M^2 \sum_{t=1}^{T}\frac{1}{t^2}\cdot \exp\left( -\tau\mu C\ln\left( \frac{T+1}{t+1}\right) \right)  =  2\sigma^2M^2 \sum_{t=1}^{T}\frac{(t+1)^2}{t^2} \cdot \frac{(t+1)^{(\tau\mu C-2)}}{(T+1)^{(\tau\mu C)}}\\
	& \leq 8\sigma^2M^2 \frac{\sum_{t=1}^{T} (t+1)^{(\tau\mu C-2)}}{(T+1)^{(\tau\mu C)}},
	\end{split}
	\end{equation*}
	where the first inequality uses condition $(B)$, the second inequality follows from condition $(A_1)$ for $t+1=t^{\ast}$, and the third inequality is derived from $(t+1)^2/t^2 \leq 4$ for all $t\geq 1$. 
	
	No matter whether $\tau\mu C > 2$ or not, we have $\sum_{t=1}^{T} t^{(\tau\mu C-2)} \leq \int_{t=1}^{T+1} t^{(\tau\mu C-2)} dt +1 $. 	
	When $ C > \frac{1}{\tau\mu}$, then $\Gamma_2^T$ can be estimated by
	\begin{equation}\label{Gamma_C2}
	\begin{split}
	\Gamma_{T}^2  \leq \frac{8\sigma^2M^2}{(\tau\mu C-1)}\cdot\frac{(T+2)^{(\tau\mu C-1)} + \tau\mu C-2}{(T+1)^{(\tau\mu C)}} \leq \frac{8\sigma^2M^2 \exp(1)}{(\tau\mu C-1)}\cdot \frac{1}{T+1} +  \frac{8\sigma^2M^2}{(T+1)^{(\tau\mu C)}} .
	\end{split}
	\end{equation}
	Combining (\ref{Gamma_C1}) and (\ref{Gamma_C2}) together, we have
	\begin{equation*}
	\begin{split}
	\E[\left\|x_{T+1} - x^{\ast} \right\|^2 ] & = \Gamma_T^1 + \Gamma_T^2 \\
	& \leq \frac{ \Delta_{n_0}^0 }{(T+1)^{(\tau\mu C)}}  +\frac{8\sigma^2M^2 \exp(1)}{(\tau\mu C-1)}\cdot \frac{1}{T+1} +  \frac{8\sigma^2M^2}{(T+1)^{(\tau\mu C)}} \\
	& \leq \frac{ \Delta_{n_0}^0  + 8\sigma^2M^2}{(T+1)^{(\tau\mu C)}} +  \frac{8\sigma^2M^2 \exp(1)}{(\tau\mu C-1)}\cdot\frac{1}{T+1}.
	\end{split}
	\end{equation*}
\end{proof}

\begin{proof}[{\bf Theorem \ref{sec:opt_thm}}]
We assume that the step size $\eta(t)$ satisfies $(A)$ and $(B_1)$. Let $n_0  := \sup \left\lbrace t \in \N^{+}: \eta(t) >  \frac{2-\tau}{2L_f} \right\rbrace $. By condition $(B_1)$, for $ t \geq (2ML_f/(2-\tau))^{1/r}$, we have
	\begin{equation}
	\eta(t) \leq \frac{M}{t^r} \leq \frac{2-\tau}{2L_f},
	\end{equation} 
	which means that $n_0 $ is a constant that is independent of $T$. Thus Lemma \ref{sec:lem2} holds.
Since the step size $\eta(t)$ satisfies $(A)$, it follows that
	\begin{equation*}
	\begin{split}
	\sum_{t=1}^{T} \eta(t) \geq \sum_{t=1}^{T}  \frac{m }{t}& \geq m\ln(T+1). 
	\end{split}
	\end{equation*}
	Recalling $\Gamma_T^1$ defined by (\ref{gamma1}),	we have
	\begin{equation}\label{sec:thm4_inequ1}
	\Gamma_{T}^1  = \exp\left( -\tau\mu\sum_{l=1}^{T}\eta(l)\right)\Delta_{n_0}^0   \leq \frac{1 }{(T+1)^{(\tau\mu m)}}\Delta_{n_0}^0 .
	\end{equation}
	In order to achieve an optimal rate $\mathcal{O}(1/T)$, we require that $\tau\mu m \geq 1$, that is $ m \geq \frac{1}{\tau\mu}$.

	Recalling the definition of $\Gamma_T^2$ as 	(\ref{gamma2}), we have 
	\begin{equation*}
	\begin{split}
	\Gamma_{T}^2 & = 2\sigma^2 \sum_{t=1}^{T}\eta(t)^2\exp\left( -\tau\mu\sum_{u>t}^{T}\eta(u)\right)   \\	 
	& = 2\sigma^2\sum_{t=1}^{C_1T^p}\eta(t)^2\exp\left( -\tau\mu\sum_{u>t}^{T}\eta(u)\right)  + 2\sigma^2 \sum_{t=C_1T^p+1}^{T}\eta(t)^2\exp\left( -\tau\mu\sum_{u>t}^{T}\eta(u)\right).
	\end{split}
	\end{equation*}
	Let $
	\Theta_1 := \sum_{t=1}^{C_1T^p}\eta(t)^2\exp\left(-\tau\mu\sum_{u>t}^{T}\eta(u)\right)$ and	$\Theta_2  :=\sum_{t=C_1T^p+1}^{T}\eta(t)^2\exp\left(-\tau\mu\sum_{u>t}^{T}\eta(u)\right)$,
then $\Gamma_{T}^2$ can be rewritten as 
	\begin{equation}\label{Gamma_2T}
	\Gamma_{T}^2 = 2\sigma^2(\Theta_1 + \Theta_2).
	\end{equation}	
	Next, we will estimate the upper bounds of $\Theta_1$ and $\Theta_2$, separately. 
	
	Let us proceed $\Theta_1$ firstly. Since the conditions $(A)$ and $(B_1)$ are satisfied, it gives that 	
\begin{align}\label{thm4:inequ1}
	\Theta_1  & = \sum_{t=1}^{C_1T^p}\eta(t)^2\exp\left( -\tau\mu\sum_{u>t}^{T}\eta(u)\right)  \leq \sum_{t=1}^{C_1T^p} \left( \frac{M_1}{t^r}\right)^2 \exp\left( -\tau\mu\sum_{u>t}^{T}\eta(u)\right)  \notag \\
	& \leq  \frac{\sum_{t=1}^{C_1T^p} \left( \frac{M_1}{t^r}\right)^2 \exp\left( -\tau\mu\sum_{u>t}^{(C_1T^p)}\eta(u)\right) }{\exp\left( \tau\mu\sum_{u>C_1T^p}^{T}\eta(u)\right) }  \leq \frac{M_1^2 \exp\left( \tau\mu m\right) \sum_{t=1}^{C_1T^p} \frac{ t^{(\tau\mu m)}}{t^{2r}}}{(C_1T^{p}+1)^{(\tau\mu m)}\exp\left( \tau\mu\sum_{u>C_1T^p}^{T}\eta(u)\right) } \notag \\
	& \leq \frac{M_1^2 \exp(\tau\mu m)\sum_{t=1}^{C_1T^p}  t^{(\tau\mu m -2r)}}{(C_1T^{p}+1)^{(\tau\mu m)}\exp\left( \tau\mu\sum_{u>C_1T^p}^{T}\eta(u)\right)}. 
\end{align}
	We know $r\in (0,1)$ and $m > \frac{1}{\tau \mu}$, so $ \tau\mu m -2r +1 > 0$. The above inequality will be 
\begin{align}\label{thm4:inequ2}	
	\Theta_1 	& \leq \frac{M_1^2\exp(\tau\mu m)}{(\tau\mu m -2r+1)}\cdot\frac{\left[ (C_1T^p+1)^{(\tau\mu m -2r+1)} + \tau\mu m -2r\right] }{(C_1T^{p}+1)^{(\tau\mu m)}\exp\left( \tau\mu\sum_{u>C_1T^p}^{T}\eta(u)\right) } \notag \\
	& \leq \frac{M_1^2\exp(\tau\mu m)}{\exp\left( \tau\mu\sum_{u>C_1T^p}^{T}\eta(u)\right) } \left[ \frac{(C_1T^p+1)^{(-2r+1)}}{\tau\mu m -2r+1} + \frac{1}{ (C_1T^{p}+1)^{(\tau\mu m)}}\right]. 
	\end{align}
	Because $\eta(t) \geq \frac{m}{t}$ for all $ 1 \leq t \leq T$, this implies that
	\begin{align}\label{thm4:inequ3}
	\exp\left( \tau\mu\sum_{u>C_1T^p}^{T}\eta(u)\right)   \geq  \exp\left( \tau\mu m\sum_{u>C_1T^p}^{T}\frac{1}{u}\right) & \geq \exp\left( \tau\mu m\int_{u=CT^p + 1}^{T+1}\frac{du}{u}\right) \notag \\ &  = \frac{(T+1)^{(\tau\mu m)}}{(C_1T^p+1)^{(\tau\mu m)}}.
	\end{align}
	Substituting (\ref{thm4:inequ3}) into (\ref{thm4:inequ2}), we have
	\begin{equation}\label{Theta_1}
	\begin{split}
	\Theta_1 &\leq  \frac{M_1^2\exp(\tau\mu m)}{(\tau\mu m -2r+1)}\cdot
	\frac{(C_1T^p+1)^{(1-2r+\tau\mu m)}}{(T+1)^{(\tau\mu m)}} + \frac{M_1^2\exp(\tau\mu m)}{(T+1)^{(\tau\mu m)}}\\
	& \leq  \frac{M_1^2\exp(\tau\mu m)}{(\tau\mu m -2r+1)}\cdot \frac{(C_1+1)^{(1-2r+\tau\mu m)}}{T^{(1-p)\tau\mu m+p(2r-1)}} + \frac{M_1^2\exp(\tau\mu m)}{(T+1)^{(\tau\mu m)}}.\\ 	 	
	\end{split}
	\end{equation}
Then we turn to bound $\Theta_2$ as follows
	\begin{equation}\label{Theta_2}
	\begin{split}
	\Theta_2 &  =\sum_{t=C_1T^p+1}^{T}\eta(t)^2\exp\left( -\tau\mu\sum_{u>t}^{T}\eta(u)\right)   \leq M_2^2\sum_{t=C_1T^p+1}^{T}\frac{1}{t^2}\cdot\exp\left( -\tau\mu m\sum_{u>t}^{T}\frac{1}{u}\right)  \\
	& \leq \frac{M_2^2 \exp(\tau\mu m)}{(T+1)^{(\tau\mu m)}}\sum_{t=C_1T^p+1}^{T}t^{(\tau\mu m-2)} \leq \frac{ M_2^2 \exp(\tau\mu m)}{(\tau\mu m -1)}\cdot\frac{(T+1)^{\tau\mu m-1} + \tau\mu m-2}{(T+1)^{(\tau\mu m)}} \\
	&
	\leq \frac{ M_2^2 \exp(\tau\mu m)}{(\tau\mu m -1)}\cdot\frac{1}{T+1} + \frac{M_2^2\exp(\tau\mu m)}{(T+1)^{(\tau\mu m)}}, \\
	\end{split}
	\end{equation}
	where the second inequality is due to (\ref{inequ_sum2}) and the third inequality follows from (\ref{sum_upper}). Therefore,  incorporating (\ref{Theta_1}) and (\ref{Theta_2}) into $\Gamma_T^2$, we have
	\begin{equation}\label{sec:thm4_inequ2}
	\begin{split}
	\Gamma_T^2  = & \, 2\sigma^2(\Theta_1 + \Theta_2) \\
   \leq	& \,  \frac{2\sigma^2 M_1^2\exp(\tau\mu m)}{(\tau\mu m -2r+1)}
	\frac{(C_1+1)^{(1-2r+\tau\mu m)}}{T^{(1-p)\tau\mu m + p(2r-1)}} \\ & + \frac{2\sigma^2(M_1^2+M_2^2)\exp(\tau\mu m)}{(T+1)^{(\tau\mu m)}}
	+ \frac{ 2M_2^2\sigma^2 \exp(\tau\mu m)}{(\tau\mu m -1)}\frac{1}{T+1}.
	\end{split}
	\end{equation}	
	When $\tau\mu m > 1 $, applying the inequalities (\ref{sec:thm4_inequ1}) and (\ref{sec:thm4_inequ2}) into (\ref{lem2_equ2}) of Lemma \ref{sec:lem2}, we get the desired result.

\end{proof}

\begin{proof}[{\bf Theorem \ref{sec2:thm6}}]	
In this case, we assume that 
\begin{align}
 m_1\leq \eta(t) \leq M_1, & \,\text{for}\,\, t \in [C_1T^{p}] \,\, \text{and }\notag \\
 \frac{m_2}{t} \leq \eta(t) \leq \frac{M_2}{t}, & \, \text{for}\,\, t \in  [T]\backslash [C_1T^{p}], \notag 
\end{align}
where  $p\in(0,1) $. Then we have 
\begin{subequations}
\begin{equation}\label{sec:thm5_inequ_1} 
	m_1 C_1T^{p}  \leq \sum_{t=1}^{C_1T^p} \eta(t) \leq M_1 C_1T^{p}, 
\end{equation}
\begin{equation}\label{sec:thm5_inequ_2}
	m_2 [\ln(T+1) - \ln(C_1T^p) -1]  \leq \sum_{C_1T^p+1}^{T} \eta(t) \leq M_2 [\ln(T) - \ln(C_1T^p)],
\end{equation}
\end{subequations}
	where (\ref{sec:thm5_inequ_2}) follows from the inequalities (\ref{inequ_sum21/l}) and (\ref{inequ_sum2}).
Let $n_0  := \sup \left\lbrace t \in \N^{+}: \eta(t) > \frac{2-\tau}{2L_f} \right\rbrace $. In this case, we assume that $n_0$ is a constant which is independent of $T$. Thus the results of Lemma \ref{sec:lem2} hold.	
	
	Recalling the definition of $\Gamma_1^T$ in (\ref{gamma1}) and applying (\ref{sec:thm5_inequ_1}) and (\ref{sec:thm5_inequ_2}), we have
	\begin{align}\label{sec2:thm6_gamma1}
	\Gamma_1^T   &  = \exp\left( -\tau\mu\sum_{t=1}^{T}\eta(t)\right) \Delta_{n_0}^0  \notag \\
	&  \leq  \exp\left( -\tau\mu\left( m_1 C_1T^{p}+m_2 (\ln(T+1) - \ln(C_1T^p)-1)\right)  \right) \Delta_{n_0}^0  \notag \\
	&  \leq  \frac{\exp(\tau\mu m_2)\Delta_{n_0}^0  }{T^{(\tau\mu m_2 (1-p))}\exp(\tau\mu m_1C_1T^{p})} \leq \frac{\exp(\tau\mu m_2)\Delta_{n_0}^0  }{T^{(\tau\mu m_2 (1-p))}\left( \tau\mu m_1C_1T^{p}+1\right)} \notag \\
	& \leq  \frac{\exp(\tau\mu m_2)}{\tau\mu m_1C_1}\cdot\frac{ \Delta_{n_0}^0 }{T^{(\tau\mu m_2 (1-p)+p)}},             
	\end{align}
	where the last inequality dues to the fact that $\exp(x) \geq 1 + x$ for $x \in \R$.	
	After that, we start to estimate $\Gamma_2^T$ which is divided into two parts as follows.
	\begin{equation*}
	\begin{split}
	\Gamma_2^T &  = 2\sigma^2\sum_{t=1}^{T}\eta(l)^2\exp(-\tau\mu\sum_{u>t}^{T}\eta(u)) \\
	& \leq 2\sigma^2\left[ \sum_{t=1}^{C_1T^p}\eta(l)^2\exp(-\tau\mu\sum_{u>t}^{T}\eta(u)) + \sum_{t=C_1T^p+1}^{T}\eta(l)^2\exp(-\tau\mu\sum_{u>t}^{T}\eta(u)) \right].
	\end{split}
	\end{equation*}
Let 
\begin{align}
\Theta_1 = \sum_{t=1}^{C_1T^p}\eta(l)^2\exp(-\tau\mu\sum_{u>t}^{T}\eta(u)),\,\, & \Theta_2 = \sum_{t=C_1T^p+1}^{T}\eta(l)^2\exp(-\tau\mu\sum_{u>t}^{T}\eta(u)).
\end{align}
 Proceeding as Theorem \ref{sec:opt_thm}, we have
	\begin{equation*}
	\Gamma_2^T \leq 2\sigma^2(\Theta_1 + \Theta_2 ).
	\end{equation*}
	In order to get the upper bound of $\Gamma_2^T$, we will estimate $\Theta_1$ and $\Theta_2$ separately. Let us evaluate $\Theta_1$ firstly.
	\begin{align*}
	\Theta_1 & = \sum_{t=1}^{C_1T^p}\eta(l)^2\exp\left( -\tau\mu\sum_{u>t}^{T}\eta(u)\right)  \leq M_1^2\sum_{t=1}^{C_1T^p}\exp\left( -\tau\mu\sum_{u>t}^{T}\eta(u)\right)  \\ 
	& \leq M_1^2\sum_{t=1}^{C_1T^p}  \frac{\exp(\tau\mu m_1 t)}{\exp\left( \tau\mu m_1 C_1T^p\right) }\exp\left( -\tau\mu\sum_{u>C_1T^p}^{T}\eta(u)\right)  \\
	& \leq M_1^2 \exp\left( -\tau\mu\sum_{u>C_1T^p}^{T}\eta(u)\right) \sum_{t=1}^{C_1T^p}\frac{\exp(\tau\mu m_1 t)}{\exp\left( \tau\mu m_1 C_1T^p\right) } \\
	& \leq  \frac{M_1^2 \exp(\tau\mu m_2) (C_1T^p)^{(\tau\mu m_2)}}{(T+1)^{(\tau\mu m_2)}}\cdot  \frac{\int_{t=1}^{C_1T^p+1}\exp(\tau\mu m_1 t)dt}{\exp(\tau\mu m_1 C_1T^p)} \\
	& \leq \frac{M_1^2\exp(\tau\mu m_2)  (C_1T^p)^{(\tau\mu m_2)}}{(T+1)^{(\tau\mu m_2)}} \cdot\frac{\exp(\tau\mu m_1(C_1T^p+1))-\exp(\tau\mu m_1)}{\tau\mu m_1 \exp(\tau\mu m_1 C_1T^p)} \leq  \frac{M_1^2\exp(\tau\mu m_2) C_1^{(\tau\mu m_2)}}{\tau\mu m_1T^{(\tau\mu m_2)(1-p)}},
	\end{align*}
	where the fourth inequality follows from (\ref{sec:thm5_inequ_2}). Next we bound $\Theta_2$ as follows.
	\begin{align*}\label{sec:thm6_theta2}
\Theta_2 & = \sum_{t=C_1T^p+1}^{T}\eta(l)^2\exp\left( -\tau\mu\sum_{u>t}^{T}\eta(u)\right)  \leq  M_2^2 \sum_{t=C_1T^p+1}^{T} \frac{1}{t^2}\cdot\exp\left( -\tau\mu m_2\sum_{u>t}^{T}\frac{1}{u}\right)  \\ 
	& \leq M_2^2 \sum_{t=C_1T^p+1}^{T} \left( \frac{1}{t}\right) ^2\exp\left( -\tau\mu m_2(\ln(T+1) - \ln(t+1) -1)\right)  \\
	& \leq \,\, \frac{M_2^2\exp(\tau\mu m_2)}{(T+1)^{(\tau\mu m_2)}}\cdot \sum_{t=C_1T^p+1}^{T} \frac{t^{(\tau\mu m_2)}}{t^2} \leq \frac{M_2^2\exp(\tau\mu m_2)}{(T+1)^{(\tau\mu m_2)}}\cdot \int_{t=C_1T^p}^{T+1} t^{(\tau\mu m_2-2)} dt  \\
	& \leq \,\, \frac{M_2^2\exp(\tau\mu m_2)}{(T+1)^{(\tau\mu m_2)}} \cdot \frac{(T+1)^{(\tau\mu m_2 -1)} - (C_1T^p)^{(\tau\mu m_2 -1)}}{\tau\mu m_2 - 1} \\
	& \leq \,\, \frac{M_2^2\exp(\tau\mu m_2)}{(\tau\mu m_2 - 1)}\cdot\frac{1}{T+1},
	\end{align*}
where the fourth inequality follows from the fact that no matter $\tau\mu m_2 > 2$ or not, we always have $\sum_{t=C_1T^p+1}^{T} \frac{t^{(\tau\mu m_2)}}{t^2} \leq \int_{t=C_1T^p}^{T+1} t^{(\tau\mu m_2-2)} dt $. The last inequality holds since $\kappa := (\tau\mu m_2)(1-p) \geq 1$ and $p\in (0,1)$, we have $\tau\mu m_2 \geq \frac{1}{(1-p)} > 1$.
	Thus 
	\begin{equation}\label{sec2:thm6_gamma2}
	\Gamma_2^T = 2\sigma^2(\Theta_1 + \Theta_2) \leq \frac{2\sigma^2M_1^2\exp(\tau\mu m_2) C_1^{(\tau\mu m_2)}}{\tau\mu m_1T^{\kappa}} + \frac{2\sigma^2M_2^2\exp(\tau\mu m_2)}{(\tau\mu m_2 - 1)} \cdot\frac{1}{T+1}.
	\end{equation}
Hence, combining (\ref{sec2:thm6_gamma1}) and (\ref{sec2:thm6_gamma2}), we obtain the  desired result.
\end{proof}

\section*{Appendix C.} 
\begin{proof}[{\bf  Theorem \ref{sec:thm4}}]	
	In this case, we assume that $\eta(t)$ satisfies the following condition:
	\begin{equation*}
	m \delta(t) \leq \eta(t) \leq M \delta(t),
	\end{equation*}	
	where $\delta(t)$ satisfies (\ref{stepsize_cond3}).	Since  $\frac{d\delta(t)}{dt} \leq 0$, then
\begin{subequations}
\begin{equation}\label{sec:thm5_inequ1}
	\sum_{u=1}^{t} \delta(u)  \geq \int_{u=1}^{t+1}\delta(u) du, 
\end{equation}
\begin{equation}\label{sec:thm5_inequ2}
	\sum_{u=l}^{t} \delta(u)	  \geq \int_{u=l}^{t+1}\delta(u) du.
\end{equation}
\end{subequations}
	
	
	Let $n_0  := \sup \left\lbrace t \in \N^{+}: \eta(t) >  \frac{2-\tau}{2L_f} \right\rbrace $. We assume that $n_0$ is a constant. Thus the conclusion of Lemma \ref{sec:lem2} holds.	Now we invoke (\ref{lem2_equ2}) and incorporate the lower and upper bounds of $\eta(t)$ into (\ref{lem2_equ2}), then apply (\ref{sec:thm5_inequ1}) and (\ref{sec:thm5_inequ2}), consequently, for $t > n_0$, we have
	\begin{align}\label{main_ineq3}
	& \E[  \left\|x_{t+1}-x^{\ast} \right\|^2 ] \notag \\
	\leq & 
	\exp\left( -\tau\mu m \sum_{l=1}^{t}\delta(t)\right)\Delta_{n_0}^0  + 2\sigma^2 M^2\sum_{l=1}^{t}\delta(l)^2\exp\left( -\tau\mu m\sum_{u>l}^{t}\delta(u)\right) \notag \\
	\leq  & 
	\exp\left( -\tau\mu m \sum_{l=1}^{t}\delta(t)\right)\Delta_{n_0}^0  + 2\sigma^2 M^2\sum_{l=1}^{t}\delta(l)^2\exp\left( -\tau\mu m\left( \sum_{u=l}^{t}\delta(u) -\delta(l)\right)\right)  \notag \\
	\leq	& \exp\left( -\tau\mu m\int_{u=1}^{t+1}\delta(u)du\right) \Delta_{n_0}^0  + 2\sigma^2 M^2\sum_{l=1}^{t}\frac{\delta(l)^2\exp(\tau\mu m \delta(l))}{\exp\left( \tau\mu m \int_{u=l}^{t+1}\delta(u) du\right)  } \notag \\ 
	\leq 	&  \frac{ \Delta_{n_0}^0 }{\exp\left(\tau\mu m\int_{u=1}^{t+1}\delta(u)du\right)} + 2\sigma^2 M^2\exp(\tau\mu m \delta(1))\sum_{l=1}^{t}\frac{\delta(l)^2}{\exp\left( \tau\mu m \int_{u=l}^{t+1}\delta(u) du\right) }. 
	\end{align} 	
	
	We consider the following three cases.
	\begin{itemize}		
		\item[{\bf 1.}] $\mathop{\lim}_{t\rightarrow\infty} \delta(t)t = 0$,	 
		that is for all $\epsilon > 0$, there exists an integer constant $t_{\epsilon} > 0 $ such that $\delta(t) t < \epsilon$ for all $t \geq t_{\epsilon}$. To attain such a convergence rate, firstly, we want to prove that for all $t \geq t_{\epsilon}$, there exists $\alpha \in (0, \frac{1}{2}]$ such that the following inequality holds 
		\begin{equation}\label{sec:thm5_inequ3}
		\exp\left( \tau\mu m\int_{t_{\epsilon}}^{t} \delta(l)dl\right)  < t^{\alpha}.
		\end{equation}
		Otherwise, there exists $t_1 \geq t_{\epsilon}$ such that for all $\alpha_1 \in (0,\frac{1}{2}]$ such that 
		\begin{equation*}
		\exp\left( \tau\mu m\int_{t_{\epsilon}}^{t_1} \delta(l)dl\right)  \geq  t_1^{\alpha_1}.
		\end{equation*}
		Thus, we have 
		\begin{equation}\label{A_ineq1}
		\tau\mu m\int_{t_{\epsilon}}^{t_1} \delta(l)dl \geq \alpha_1 \ln(t_1).
		\end{equation} 
		We know that the integral of $\delta(t)$ from $t_{\epsilon}$ to $t_1$ can be rewritten as $$\int_{t_{\epsilon}}^{t} \delta(l)dl = \int_{t_{\epsilon}}^{t} \delta(l)\cdot l \cdot\frac{1}{l}dl.$$
		Since $\delta(t) t < \epsilon $ for $ t \geq t_{\epsilon}$, then $\int_{t_{\epsilon}}^{t} \delta(l)\cdot l \cdot\frac{1}{l}dl < \epsilon \ln(\frac{t_1}{t_{\epsilon}})$. 
		This is contradictory with (\ref{A_ineq1}) for small $\epsilon < \frac{\alpha_1}{\tau\mu m} $. Thus for all $t \geq t_{\epsilon}$, the inequality (\ref{sec:thm5_inequ3}) holds
		for a constant $\alpha \in (0,\frac{1}{2}].$ 		
		Then
		\begin{equation*}
		\begin{split}
		& \sum_{l=1}^{t}\delta(l)^2 \exp\left( -\tau\mu m\int_{u=l}^{t+1} \delta(u)du\right) \\
		= \,\,	&   \sum_{l=1}^{t_{\epsilon}-1}\delta(l)^2 \exp\left( -\tau\mu m\int_{u=l}^{t+1} \delta(u)du\right) dl + \sum_{t_{\epsilon}}^{t}\delta(l)^2 \exp\left( -\tau\mu m\int_{u=l}^{t+1} \delta(u)du\right) \\ 
		\leq\,\, &
		\delta(1)^2\exp\left( -\tau\mu m\int_{u=t_{\epsilon}-1}^{t+1}\delta(u)du\right) (t_{\epsilon}-1) + \frac{\sum_{t_{\epsilon}}^{t}\left( \frac{\epsilon}{l}\right) ^2 \exp\left( \tau\mu m\int_{u=t_{\epsilon}}^{l} \delta(u)du\right)}{\exp\left( \tau\mu m\int_{l=t_{\epsilon}}^{t+1} \delta(l)dl\right) } \\
		\leq\,\, &  \delta(1)^2(t_{\epsilon}-1)\exp\left( -\tau\mu m\int_{u=t_{\epsilon}-1}^{t+1}\delta(u)du\right)  + \frac{\sum_{t_{\epsilon}}^{t}(\frac{\epsilon}{l})^2 (l+1)^{\alpha}}{\exp\left( \tau\mu m\int_{l=t_{\epsilon}}^{t+1} \delta(l)dl\right) }\\
		\leq\,\, &   \left[ \delta(1)^2(t_{\epsilon}-1) + 2\epsilon^2\right] \exp\left( -\tau\mu m\int_{l=t_{\epsilon}}^{t+1}\delta(l)dl\right)  \\
		\leq\,\, &  \frac{\delta(1)^2(t_{\epsilon}-1) + 2\epsilon^2}{\exp\left( -\tau\mu m \int_{l=1}^{t_{\epsilon}}\delta(l)dl\right) }\exp\left( -\tau\mu m\int_{l=1}^{t+1}\delta(l)dl\right) ,
		\end{split}
		\end{equation*}
		where the third inequality follows from the fact that $\sum_{t_{\epsilon}}^{t}(\frac{\epsilon}{l})^2 (l+1)^{\alpha} \leq 2\epsilon^2 $. 
		Thus, in this case,  for $t > n_0$, $\E[ \left\|x_{t+1} -x^{\ast} \right\|^2]  $ is at most
\begin{equation*}
\left(\Delta_{n_0}^0  + 2\sigma^2 M^2\exp(\tau\mu m \delta(1)) \frac{\delta(1)^2(t_{\epsilon}-1) + 2\epsilon^2}{\exp(-\tau\mu  \int_{l=1}^{t_{\epsilon}}\delta(l)dl)}\right) \exp\left( -\tau\mu m\int_{u=1}^{t+1}\delta(u)du\right) . 	
\end{equation*}
		
		\item[{\bf 2.}] $\lim_{t\rightarrow\infty} \delta(t)t = 1 $.
		
		In this case, it is easy to show there exist $m$ and $M$ such that $ \frac{m}{t} \leq \eta(t) \leq \frac{M}{t}$. Hence the theorem follows from Theorem \ref{sec:thm1}.
		\item[{\bf 3.}] $\lim_{t\rightarrow\infty} \delta(t)t = +\infty $, that is for any $M_1 > 0$, there exists a constant $T_M \in \N^{+}$ such that for all $t \geq T_M$, $\delta(t) t > M_1$.

		We suppose that  there exists a constant $c_1 \leq \frac{\tau\mu m}{2}$ such that for all $t \geq T_M$
		\begin{equation}\label{C_ineq}
		- \frac{d\delta(t)}{dt} \leq c_1 \delta(t)^2.
		\end{equation}

		Let $P(l) := \delta(l)^2\exp\left( -\tau\mu m \int_{u=l}^{t+1}\delta(u) du\right) $ for $ 1 \leq l \leq t$, then
		\begin{align}\label{nabla_P}
		 \frac{dP(l)}{dl}		& = 2\delta(l)\frac{d\delta(l)}{dl}\exp\left( -\tau\mu m \int_{u=l}^{t+1}\delta(u) du\right)  + \tau\mu m\delta(l)^3\exp\left( -\tau\mu m \int_{u=l}^{t+1}\delta(u) du\right) \notag \notag \\
	&  	=\delta(l)\exp\left( -\tau\mu m \int_{u=l}^{t+1}\delta(u) du\right) \left[ 2\frac{d\delta(l)}{dl} + \tau\mu m\delta(l)\delta(l)\right].
		\end{align}
		Let $Q(l) := 2\frac{d\delta(l)}{dl} + \tau\mu m\delta(l)\delta(l) $.
		By (\ref{nabla_P}), we know that the sign of $\frac{dP(l)}{dl}$ is determined by the sign of $Q(l)$. If $c_1 \leq \frac{\tau\mu m}{2}$, from (\ref{C_ineq}), we have $Q(l) \geq  0$, then the sequence of $P(l)$ is increasing when $ l \geq T_{M}$.  
		
		If $P(u)$ is increasing for $u\in [l, t]$, then
		\begin{equation}\label{C_ineq2}
		\sum_{u=l}^{t} P(u) \leq \int_{u=l}^{t+1} P(u)du.
		\end{equation}
		Otherwise, if $P(u)$ is decreasing for $u\in [l, t]$, then 
		\begin{equation}\label{C_ineq3}
		\sum_{u=l}^{t} P(u) \leq P(l) + \int_{u=l}^{t} P(u)du.
		\end{equation}
		By (\ref{C_ineq2}), we have
		\begin{equation}\label{sum_integral}	
		\sum_{l=1}^{t} P(l)	 = \sum_{l=1}^{T_M} P(l) + \sum_{l= T_M +1}^{t} P(l)  \leq \sum_{l=1}^{T_M} P(l) + \int_{l=T_M}^{t+1} P(l) dl.
		\end{equation}

		By integration by parts,   $\int_{l=T_M}^{t+1} P(l) dl$ can be written as
		\begin{align*}\label{C_integral}		
		& \tau\mu m\int_{l=T_M}^{t+1} P(l) dl = \tau\mu m\int_{l=T_M}^{t+1}\delta(l)^2\exp\left( -\tau\mu m \int_{u=l}^{t+1}\delta(u) du\right) dl \notag \\
		= &  \delta(t+1) - \delta(T_M)\exp\left( -\tau\mu m \int_{u=T_M}^{t+1}\delta(u) du\right)  - \int_{l=T_M}^{t+1}\frac{d\delta(l)}{dl}\exp\left( -\tau\mu m \int_{u=l}^{t+1}\delta(u) du\right) dl \notag \\
		\leq  &  \delta(t+1) - \delta(T_M)\exp\left( -\tau\mu m \int_{u=T_M}^{t+1}\delta(u) du\right)  + c_1 \int_{l=T_M}^{t+1}\delta(l)^2\exp\left( -\tau\mu m \int_{u=l}^{t+1}\delta(u)du\right) dl,	 
		\end{align*}  
		where the above inequality holds because (\ref{C_ineq}) satisfies. When $c_1 < \tau\mu m $, rearranging the above inequality, we have
		\begin{equation*}
		\int_{l=T_M}^{t+1}\delta(l)^2\exp\left( -\tau\mu m \int_{u=l}^{t+1}\delta(u) du\right) dl \leq \frac{\delta(t+1) - \delta(T_M)\exp\left(-\tau\mu m \int_{u=T_M}^{t+1}\delta(u) du\right)}{(\tau\mu m - c_1)}.
		\end{equation*}   
		Hence,
		\begin{align*}
		&\sum_{l=1}^{t+1}\delta(l)^2\exp\left( -\tau\mu m \int_{u=l}^{t+1}\delta(u) du\right) dl \leq \sum_{l=1}^{T_M} P(l) + \int_{l=T_M}^{t+1} P(l) dl \\
		=  & \sum_{l=1}^{T_M}\delta(l)^2\exp\left( -\tau\mu m \int_{u=l}^{t+1}\delta(u) du\right)  + \int_{l=T_M}^{t+1}\delta(l)^2\exp\left( -\tau\mu m \int_{u=l}^{t+1}\delta(u) du\right) dl \\
		\leq & \frac{\delta(1)^2T_M}{\exp\left( \tau\mu m \int_{u=T_M}^{t+1}\delta(u) du\right) } + \frac{\delta(t+1) - \delta(T_M)\exp(-\tau\mu m \int_{u=T_M}^{t+1}\delta(u) du)}{(\tau\mu m - c_1)}  \\
		= &\frac{\delta(t+1)}{(\tau\mu m - c_1)} + \frac{\delta(1)^2T_M- \frac{\delta(T_M)}{(\tau\mu m - c_1)}}{\exp\left( \tau\mu m \int_{u=T_M}^{t+1}\delta(u) du\right) } \\
		\leq  &\frac{\delta(t+1)}{(\tau\mu m - c_1)} + \frac{\delta(1)^2T_M}{\exp\left( -\tau\mu m \int_{u=1}^{T_M}\delta(u) du\right) } \exp\left( -\tau\mu m \int_{u=1}^{t+1}\delta(u) du\right) .
		\end{align*}

		Finally, incorporating the above inequality into (\ref{main_ineq3}), we can show that $\E[\left\|x_{t+1} -x^{\ast} \right\|^2]$ is bounded by
		\begin{equation*}
		\begin{split}
		 & 
		\frac{\varepsilon_2}{(\tau\mu m - c_1)}
		\delta(t+1)  + \left[\Delta_{n_0}^0  +\frac{\varepsilon_2\delta(1)^2T_M}{\exp(-\tau\mu m \int_{u=1}^{T_M}\delta(u) du)} \right]  \exp\left( -\tau\mu m \int_{u=1}^{t+1}\delta(u) du\right),
		\end{split}
		\end{equation*}
		where $\varepsilon_2 =2\sigma^2M^2 \exp(\tau\mu m \delta(1))$.
		\end{itemize}	
\end{proof}
\begin{proof}[{\bf Lemma \ref{sec:lem4}}]
	Suppose that there exists a constant $c_1 > 0$ such that $$- \frac{d\delta(t)}{dt} \leq c_1 \delta(t)^2.$$ Let $\hat{\delta}(t) = a \delta(t)$ for $a > 0$. Of course, for the new function $\hat{\delta}(t)$, there must be a constant $\hat{c}_1 >  0 $ such that $$- \frac{d\hat{\delta}(t)}{dt} \leq \hat{c}_1 \hat{\delta}(t)^2.$$ Then we have $$ \frac{d\hat{\delta}(t)}{dt} = - a \frac{d\delta(t)}{dt} \leq \hat{c}_1 \hat{\delta}(t)^2 = a^2\hat{c}_1 \delta(t)^2.$$ Thus, $$\frac{d\delta(t)}{dt} \leq a\hat{c}_1 \delta(t)^2.$$ Let $ 0 < a \leq \frac{\tau\mu m}{2\hat{c_1}}$, we have $a\hat{c}_1 \leq \frac{\tau\mu m}{2}$, which shows that there must be a constant $c_1  = a\hat{c}_1 \leq \frac{\tau\mu m}{2}$.

\end{proof}

\section*{Appendix D.}
\begin{proof}[{\bf Theorem \ref{sec4:thm1}}]
	We assume that $\eta(t)$ satisfies the following condition
	\begin{equation*}
	\frac{m}{t+1}\leq \eta(t) \leq \frac{M\ln(t+1)}{t+1},  \, \forall \, 1 \leq t \leq T.
	\end{equation*}
	Let $n_0  := \sup \left\lbrace t \in \N^{+}: \eta(t) >  \frac{2-\tau}{2L_f} \right\rbrace $. For $t\geq (2L_fM/(2-\tau))^2$, we have
	\begin{equation}
	 \eta(t) \leq \frac{M\ln(t+1)}{t+1} \leq \frac{M\sqrt{t+1}}{t+1} \leq \frac{2-\tau}{2L_f}.
	\end{equation}
	Then $n_0$ must exist and is a constant which is independent of $T$. Thus the inequality (\ref{lem2_equ2}) of Lemma \ref{sec:lem2} holds, then we get
	\begin{align}\label{sec4:thm1:ineq1}
	& \E[  \left\|x_{t+1}-x^{\ast} \right\|^2]  \notag \\
	\leq &  \exp\left( -\tau\mu\sum_{l=1}^{t}\eta(l)\right)\Delta_{n_0}^0  + 2\sigma^2\sum_{l=1}^{t}\eta(l)^2\exp\left( -\tau\mu\sum_{u>l}^{t}\eta(u)\right) \notag  \\
	\leq &  \exp\left( -\tau\mu m\sum_{l=1}^{t} \frac{1}{l+1}\right)\Delta_{n_0}^0 + 2\sigma^2M^2 \sum_{l=1}^{t}\frac{\ln^2(l+1)}{(l+1)^2}\exp\left( -\tau\mu m\sum_{u>l}^{t}\frac{1}{u+1}\right)  \notag \\
	\leq & \frac{\Delta_{n_0}^0}{\exp(\tau\mu m (\ln(t+2)-\ln2))} + 2\sigma^2M^2\exp(\tau\mu m) \sum_{l=1}^{t} \frac{\ln^2(l+1)}{(l+1)^2} \cdot \frac{\exp(\tau\mu m \ln(l+1))}{\exp(\tau\mu m \ln(t+2))} \notag \\
	\leq &  \frac{2^{(\tau\mu m)} \Delta_{n_0}^0}{(t+2)^{(\tau\mu m)}} + \frac{2\sigma^2M^2\exp(\tau\mu m)}{ (t+2)^{(\tau\mu m)}} \sum_{l=1}^{t}\frac{\ln^2(l+1)}{(l+1)^2} (l+1)^{(\tau\mu m)} \notag \\
	\leq &  \frac{2^{(\tau\mu m)} \Delta_{n_0}^0}{(t+2)^{(\tau\mu m)}} + \frac{2\sigma^2M^2\exp(\tau\mu m)}{ (t+2)^{(\tau\mu m)}}\left[ \frac{\ln(2)}{2}+ \int_{l=2}^{t+2}\frac{\ln^2(l)}{l^2}\cdot (l)^{(\tau\mu m)}  dl\right], 
	\end{align}
	where the third inequality follows from (\ref{inequ_sum1}) and (\ref{inequ_sum2}), and the last inequality  is obtained from (\ref{sum_integral}).	If $\tau\mu m  = 1$, we have
	\begin{equation*}
	\int_{l=2}^{t+2}\frac{\ln^2(l)}{l^2}\cdot (l)^{(\tau\mu m)} dl = \frac{\ln^3(t+2)}{3} - \ln^32 < \frac{\ln^3(t+2)}{3}.
	\end{equation*}
	Otherwise, if $\tau\mu m  \neq 1$, integrating by parts we get
	\begin{equation*}
	\begin{split}
		&\int_{l=2}^{t+2}\frac{\ln^2(l)}{l^2}\cdot l^{(\tau\mu m)} dl \\ \leq & \, \frac{(t+2)^{(\tau\mu m -1)} \ln^2(t+2) - 2^{(\tau\mu m -1)}\ln^22}{(\tau\mu m-1)} +\frac{2^{(\tau\mu m)}\ln2 }{(\tau\mu m -1)^2}+ \frac{2[(t+2)^{(\tau\mu m -1)} - (2)^{(\tau\mu m -1)}]}{(\tau\mu m -1)^3}.
	\end{split}
	\end{equation*}
	From the above inequality,  we can see that if $\tau\mu m  < 1$, such an integral can be bounded by a scalar
	\begin{equation*}
\int_{l=2}^{t+2}\frac{\ln^2(l)}{l^2}\cdot l^{(\tau\mu m)} dl \leq \frac{2^{(\tau\mu m)}}{(1- \tau\mu m)^3} +\frac{2^{(\tau\mu m)}\ln2 }{(1- \tau\mu m)^2} + \frac{2^{(\tau\mu m-1)}\ln^22 }{(1- \tau\mu m)} \leq \frac{2+2\ln2+\ln^22}{(1-\tau\mu m)^3}.
	\end{equation*}
	While $\tau\mu m > 1$, then
	\begin{equation*}
\int_{l=2}^{t+2}\frac{\ln^2(l)}{l^2}\cdot l^{(\tau\mu m)} dl \leq \left[ \frac{\ln^2(t+2)}{(\tau\mu m -1)} + \frac{2}{(\tau\mu m-1)^3}\right] (t+2)^{(\tau\mu m -1)} +\frac{2^{(\tau\mu m)}\ln2 }{(\tau\mu m -1)^2}.
	\end{equation*}
	
	Thus, collecting the results obtained above, let $t= T$, we can get the result as desired.
	
\end{proof}

\begin{proof}[{\bf Theorem \ref{sec4:thm2}}]	
	In this case, $\eta(t)$ satisfies that
	\begin{equation*}
	\frac{m}{t}\leq \eta(t) \leq \frac{M}{t^{\alpha}},
	\end{equation*} 
	for $\alpha \in (1/2, 1]$. 	Let $n_0  := \sup \left\lbrace t \in \N^{+}: \eta(t) > \frac{2-\tau}{2L_f} \right\rbrace $. For $t\geq (2L_fM/(2-\tau))^{(1/\alpha)}$, we have
		\begin{equation}
		 \eta(t) \leq \frac{M}{t^{\alpha}} \leq \frac{2-\tau}{2L_f}.
		\end{equation}
		Then $n_0$ must exist and is a constant which is independent of $T$. Thus in this case, the inequality (\ref{lem2_equ2}) of Lemma \ref{sec:lem2} holds.	From (\ref{lem2_equ2})  in Lemma \ref{sec:lem2}, we have	
	\begin{equation*}
	\begin{split}
 \E[  \left\|x_{t+1}-x^{\ast} \right\|^2] \leq &\, \exp\left( -\tau\mu\sum_{l=1}^{t}\eta(l)\right)\Delta_{n_0}^0 + 2\sigma^2\sum_{l=1}^{t}\eta(l)^2\exp\left( -\tau\mu\sum_{u>l}^{t}\eta(u)\right)  \\ 
	 \leq &\, \exp\left( -\tau\mu m\sum_{l=1}^{t}\frac{1}{l}\right)\Delta_{n_0}^0 + 2\sigma^2M^2\sum_{l=1}^{t}\frac{1}{l^{2\alpha}}\exp\left( -\tau\mu m\sum_{u>l}^{t}\frac{1}{u}\right)  \\
\leq 	& \,\frac{ \Delta_{n_0}^0}{(t+1)^{(\tau\mu m )}} + \frac{2\sigma^2M^2\exp(\tau\mu m)}{ (t+1)^{(\tau\mu m)}} \sum_{l=1}^{t} l^{(\tau\mu m - 2\alpha)} \\
\leq	& \, \frac{\Delta_{n_0}^0}{(t+1)^{(\tau\mu m )}} + \frac{2\sigma^2M^2\exp(\tau\mu m)}{ (t+1)^{(\tau\mu m)}} \left( \int_{l=1}^{t+1} l^{(\tau\mu m - 2\alpha)}dl + 1\right)  \\
\leq 	& \,\frac{\Delta_{n_0}^0 +2\sigma^2M^2\exp(\tau\mu m)}{(t+1)^{(\tau\mu m)}} + \frac{2\sigma^2M^2\exp(\tau\mu m) }{ (t+1)^{(\tau\mu m)}}\int_{l=1}^{t+1} l^{(\tau\mu m - 2\alpha)}dl.
	\end{split}
	\end{equation*}
	If $\tau\mu m = 2\alpha - 1 > 0$, it follows that  
	\begin{equation*}
	\int_{l=1}^{t+1} l^{(\tau\mu m - 2\alpha)}dl =  \int_{l=1}^{t+1} \frac{dl}{l} = \ln(t+1).
	\end{equation*}	
Consequently,
	\begin{equation*}
	\E[  \left\|x_{t+1}-x^{\ast} \right\|^2]  \leq \frac{\Delta_{n_0}^0 +2\sigma^2M^2\exp(2\alpha-1)}{(t+1)^{(2\alpha-1)}} + \frac{2\sigma^2M^2\exp(2\alpha-1) \ln(t+1)}{ (t+1)^{(2\alpha-1)}}.
	\end{equation*}
	If 	$\tau\mu m \neq  2\alpha -1$, we have
	\begin{equation*}
	\int_{l=1}^{t+1} l^{(\tau\mu m - 2\alpha)}dl =  \frac{(t+1)^{(\tau\mu m - 2\alpha +1)} -1 }{(\tau\mu m - 2\alpha +1)},
	\end{equation*}
then $\E[  \left\|x_{t+1}-x^{\ast} \right\|^2] $ is at most
	\begin{equation*}
\frac{\Delta_{n_0}^0 +2\sigma^2M^2\exp(2\alpha-1)}{(t+1)^{(\tau\mu m)}} + \frac{2\sigma^2M^2\exp(\tau\mu m)  }{(\tau\mu m - 2\alpha +1) }\left[ \frac{1}{(t+1)^{(2\alpha-1)} } - \frac{1}{(t+1)^{(\tau\mu m)}}\right].
	\end{equation*}
	Combing the above results and let $t =T$, we obtain the desired result. 	
\end{proof}

\begin{proof}[{\bf Theorem \ref{sec4:thm3}}]
	In this case, we assume that $\eta(t)$ satisfies that
	\begin{equation*}
	\frac{m}{(t+1)\ln(t+1)}\leq \eta(t) \leq \frac{M}{(t+1)^{\alpha}}
	\end{equation*} 
	for $\alpha \in (1/2, 1]$. 	Let $n_0  := \sup \left\lbrace t \in \N^{+}: \eta(t) >  \frac{2-\tau}{2L_f} \right\rbrace $. For $t\geq (2L_fM/(2-\tau))^{(1/\alpha)}-1$, we have
	\begin{equation}
	\eta(t) \leq \frac{M}{(t+1)^{\alpha}} \leq \frac{2-\tau}{2L_f}.
	\end{equation}
	Therefore $n_0$ must exist and is a constant. In this case, the inequality (\ref{lem2_equ2}) of Lemma \ref{sec:lem2} holds.
	By (\ref{lem2_equ2}), we have	
	\begin{align}\label{sec4:thm3_inequ}
	& \E[  \left\|x_{t+1}-x^{\ast} \right\|^2]  \notag  \\  \leq &  \exp\left( -\tau\mu\sum_{l=1}^{t}\eta(l)\right) \Delta_{n_0}^0 + 2\sigma^2\sum_{l=1}^{t}\eta(l)^2\exp\left( -\tau\mu\sum_{u>l}^{t}\eta(u)\right)  \notag \\
\leq	&  \exp\left( -\tau\mu m\sum_{l=1}^{t}\frac{1}{(l+1)\ln(l+1)}\right) \Delta_{n_0}^0 + 2\sigma^2 M^2\sum_{l=1}^{t}\frac{\exp\left( -\tau\mu m\sum_{u>l}^{t}\frac{1}{(l+1)\ln(l+1)}\right) }{(l+1)^{2\alpha}} \notag \\
 \leq 	&  \frac{(\ln2)^{(\tau\mu m)}\Delta_{n_0}^0}{(\ln(t+2))^{(\tau\mu m)}} + \frac{2\sigma^2 M^2(\ln2)^{(\tau\mu m)}}{(\ln(t+2))^{(\tau\mu m)}}\sum_{l=1}^{t} \frac{(\ln(l+1))^{(\tau\mu m)}}{(\ln(t+2))^{2\alpha}} \notag \\
 \leq 	&\frac{(\ln2)^{(\tau\mu m)}\Delta_{n_0}^0}{(\ln(t+2))^{(\tau\mu m)}} +  \frac{2\sigma^2 M^2(\ln2)^{(\tau\mu m)}}{(\ln(t+2))^{(\tau\mu m)}}\left[ \frac{(\ln2)^{\tau\mu m}}{2^{2\alpha}} + \int_{l=1}^{t+1}\frac{(\ln(l+1))^{(\tau\mu m)}}{(l+1)^{2\alpha}}dl\right], 
\end{align}
	where the third inequality dues to the fact that  $\sum_{l=1}^{t}\frac{1}{l\ln(l)} \geq \int_{l=1}^{t+1}\frac{1}{(l+1)\ln(l+1)}dl = \ln\ln(t+2) -\ln\ln2 $ and the last inequality follows from (\ref{sum_integral}).
	
	We know that for any $\beta \in (0,1)$, there must be a constant $t_{\beta}$ such that $\ln(t+1) \leq (t+1)^\beta$ for all $t \geq t_{\beta}$. Here we choose that $0 < \beta < \frac{2\alpha-1}{\tau\mu m}$. There exists a constant $t_{\beta}$ such that $\ln(t+1) \leq (t+1)^\beta$ for all $t \geq t_{\beta}$. For sufficiently large $t \geq t_{\beta}$, we have
	\begin{align}\label{sec:thm3_inequ1}
	\int_{l=1}^{t+1}\frac{(\ln(l+1))^{(\tau\mu m)}}{(l+1)^{2\alpha}}dl & \leq \int_{l=1}^{t_{\beta}}\frac{(\ln(l+1))^{(\tau\mu m)}}{(l+1)^{2\alpha}}dl +  \int_{t_{\beta}}^{t+1}\frac{(\ln(l+1))^{(\tau\mu m)}}{(l+1)^{2\alpha}}dl \notag \\
	& \leq (\ln(t_{\beta}+1))^{(\tau\mu m)} \int_{l=1}^{t_{\beta}}\frac{dl}{(l+1)^{2\alpha}} + \int_{t_{\beta}}^{t+1} (l+1)^{(\beta\tau\mu m -2\alpha)}dl  \notag \\
	& \leq \frac{2^{(1-2\alpha)}}{2\alpha-1} + \frac{(t+1)^{(\beta\tau\mu m -2\alpha +1)}-(t_{\beta}+1)^{(\beta\tau\mu m -2\alpha +1)}}{(\beta\tau\mu m + 1-2\alpha )}.
	\end{align}
	Thus, applying (\ref{sec:thm3_inequ1}) into (\ref{sec4:thm3_inequ}) and let $t=T$, we can bound $ \E[  \left\|x_{T+1}-x^{\ast} \right\|^2] $ by
	\begin{align*}
	&  \frac{(\ln2)^{(\tau\mu m)}\Delta_{n_0}^0}{(\ln(t+2))^{(\tau\mu m)}} +  \frac{2\sigma^2 M^2(\ln2)^{(\tau\mu m)}}{(\ln(t+2))^{(\tau\mu m)}}\left[ \frac{(\ln2)^{(\tau\mu m)}}{2^{2\alpha}} + \frac{2^{(1-2\alpha)}}{2\alpha-1} + \frac{(t_{\beta}+1)^{(\beta\tau\mu m -2\alpha +1)}}{(2\alpha - 1 -\beta\tau\mu m)}\right].
	\end{align*}
	Therefore, there exists a constant $C_2 > 0$ such that
	\begin{equation*}
	\E[  \left\|x_{T+1}-x^{\ast} \right\|^2]   \leq  \frac{C_2}{(\ln(t+2))^{(\tau\mu m)}}.
	\end{equation*}	
\end{proof}

\bibliography{sgd_sgd.bib}

\end{document}